\documentclass[10pt,letterpaper]{amsart}
\usepackage{mathrsfs}
\usepackage{graphicx}
\usepackage{amssymb}
\usepackage{epstopdf}
\usepackage{amsmath,amscd}
\usepackage{amsthm}
\usepackage{enumerate}
\usepackage{url,verbatim}
\usepackage{subfig}
\usepackage{enumitem}

\usepackage[normalem]{ulem}
\usepackage{fixme}

\RequirePackage[colorlinks,citecolor=blue,urlcolor=blue]{hyperref}
\usepackage{breakurl}
\theoremstyle{plain}
\newtheorem{theorem}{Theorem}
\newtheorem{definition}[theorem]{Definition}
\newtheorem{lemma}[theorem]{Lemma}
\newtheorem{proposition}[theorem]{Proposition}
\newtheorem{corollary}[theorem]{Corollary}

\newtheorem{assumption}[theorem]{Assumption}
\newtheorem{remark}[theorem]{Remark}

\newcommand\ol{\overline}

\newcommand\EE{{\mathbb E}}

\newcommand\RR{{\mathbb R}}
\newcommand\ZZ{{\mathbb Z}}
\newcommand\NN{{\mathbb N}}

\newcommand\HH{{\mathbb H}}

\newcommand\YY{{\mathbb {Y}}}

\renewcommand\a{\alpha}

\renewcommand\ell{l}

\newcommand\CC{\mathbb{C}}

\newcommand\bA{\mathbf{A}}
\newcommand\bB{\mathbf{B}}
\newcommand\bC{\mathbf{C}}
\newcommand\bD{\mathbf{D}}
\newcommand\MP{\mathbb{MP}}

\newcounter{mycount}

\numberwithin{equation}{section}
\numberwithin{theorem}{section}
\numberwithin{figure}{section}

\title[Asymptotics of pure dimer coverings on rail-yard graphs]{Asymptotics of pure dimer coverings on rail-yard graphs}

\date{}

\author{Zhongyang Li}
\address{Department of Mathematics,
University of Connecticut,
Storrs, Connecticut 06269-3009, USA}
\email{zhongyang.li@uconn.edu}
\urladdr{\url{https://mathzhongyangli.wordpress.com}}

\author{Mirjana Vuleti{\'c}}
\address{Department of Mathematics,
University of Massachussetts Boston,
Storrs, Connecticut 06269-3009, USA}
\email{mirjana.vuletic@umb.edu}
\urladdr{\url{https://www.math.umb.edu/~vuletic/Site/Homepage.html}}

\begin{document}
\maketitle

\begin{abstract}We study the asymptotic limit of random pure dimer coverings on rail yard graphs when the mesh sizes of the graphs go to 0. Each pure dimer covering corresponds to a sequence of interlacing partitions starting with an empty partition and ending in an empty partition. Under the assumption that the probability of each dimer covering is proportional to the product of weights of present edges, we obtain the limit shape (law of large numbers) of the rescaled height functions and the convergence of the unrescaled height fluctuations to a diffeomorphic image of Gaussian free field (Central Limit Theorem); answering a question in \cite{bbccr}. Applications include the limit shape and height fluctuations for pure steep tilings (\cite{BCC17}) and pyramid partitions (\cite{Ken05,Sze08,BY09,you10}). The technique to obtain these results is to analyze a class of Macdonald processes which involve dual partitions as well.

\end{abstract}

\section{Introduction}

A dimer cover, or perfect matching on a graph is a subset of edges such that each vertex is incident to exactly one edge in the subset. A dimer model is a probability measure on the set of perfect matchings (see \cite{RK09}).
The dimer model is a natural mathematical model for the structures of matter; for example, each perfect matching on a hexagonal lattice corresponds to a double-bond configuration of a graphite molecule; the dimer model on a Fisher graph has a measure-preserving correspondence with the 2D Ising model (see \cite{Fi66,MW73,ZL12}).

Just as in the structures of matter different molecule configurations have certain probabilities to occur depending on the underlying energy, mathematically we define a probability measure on the set of all perfect matchings of a graph depending on the energy of the dimer configuration, quantified as the product of weights of present edges in the configuration.
The phase transitions and asymptotical behaviors of the dimer model have been an interesting topic for mathematicians and physicists for a long time. A combinatorial argument shows that the total number of perfect matchings on any finite planar graph can be computed by the Pfaffian of the corresponding weighted adjacency matrix (\cite{Kas61,TF61}). The local statistics can be computed by the inverse of the weighted adjacency matrix (\cite{Ken01}); a complete picture of phase transitions was obtained in \cite{KOS}. Empirical results show that in large graphs, there are certain regions where the configurations are almost deterministic, i.e. one type of edges have very high probability to occur in the dimer configuration. These are called ``frozen regions", and their boundaries are called ``frozen boundaries".
When the mesh size of the graph goes to 0 such that the graph approximates a simply-connected region in the plane, the limit shape of the height function can be obtained by a variational principle (\cite{ckp}), and the frozen boundaries are proved to be algebraic curves of a specific type called the cloud curves (\cite{KO07}). It is also known that the fluctuations of (unrescaled) dimer heights converge to the Gaussian free field (GFF) in distribution when the boundary satisfies certain conditions (\cite{Ken01,Li13}).

In this paper, we investigate perfect matchings on a general class of bipartite graphs called rail-yard graphs. The major goal of the paper is to understand the asymptotic behavior of the  model, in particular the limit shape and height fluctuations.

We start with pyramid partitions as an example of rail-yard graphs. They are shown in the figure below.
\begin{figure}[!h]
\includegraphics[width=.250\textwidth]{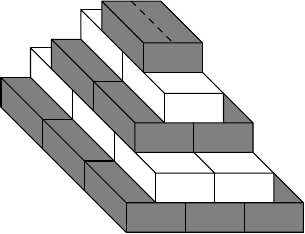}\qquad\qquad \includegraphics[width=.25\textwidth]{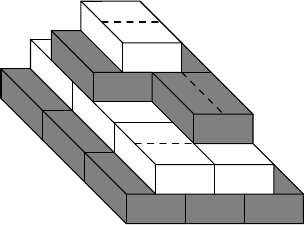}
\end{figure}
These are pyramid shaped objects built out of square bricks. The fundamental pyramid partition is shown on the left and extends infinitely down and any other pyramid partition is obtained by removing finitely many bricks, where we can only remove fully exposed bricks at any given time. In the figure on the right, 3 bricks have been removed. If one looks from above, then the following domino tilings will be observed.
\begin{figure}[!h]
\includegraphics[width=.25\textwidth]{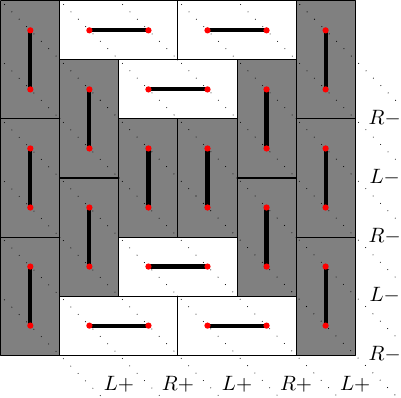}\qquad\qquad\qquad\includegraphics[width=.25\textwidth]{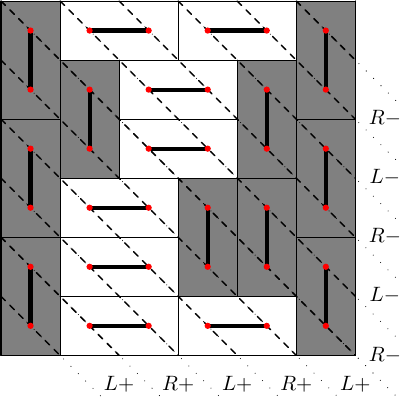}
\end{figure}

In the above figures an edge was drawn for each domino. Now it is possible to slice the tilings diagonally (along dashed lines) and insert new vertices to obtain the following two perfect matchings on the same graph. The graph and the perfect matching for the fundamental pyramid partition is shown below.
\begin{figure}[!h]
\includegraphics[width=.8\textwidth]{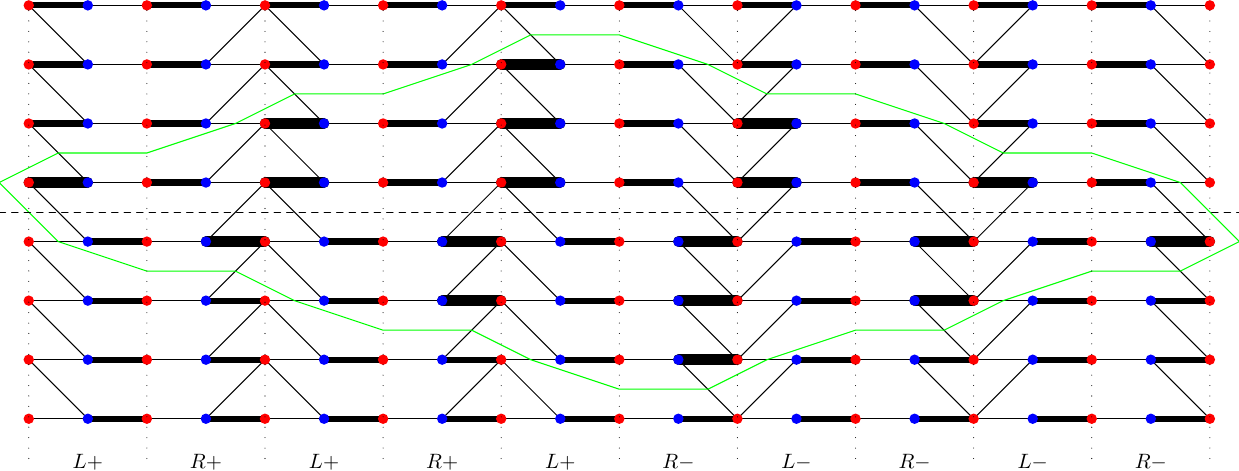}
\end{figure}

The one corresponding to the other pyramid partitions where 3 bricks were removed is given below.

\begin{figure}[!h]
\includegraphics[width=.8\textwidth]{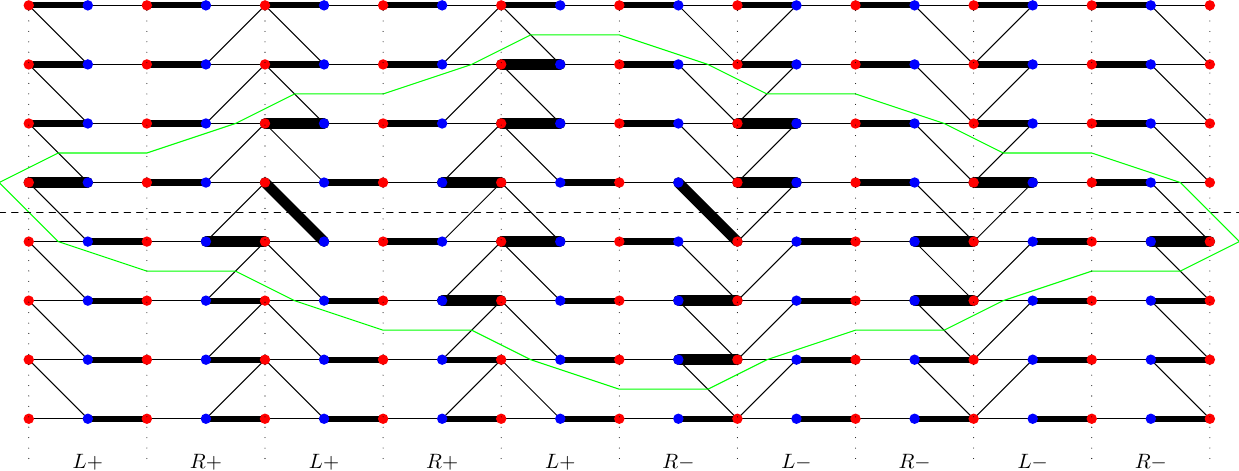}
\end{figure}

This construction is due to \cite{bbccr}. We leave the explanation to the picture. Note that the green curve bounds the region that corresponds to the tiling of the square domain shown above and that edges shown in bold correspond to dominoes. Note also that to represent an arbitrary pyramid partition we might need to enlarge the graph to the left and right if the bricks that were removed fall out of this square domain. This means that if the graph is fixed then only a subset of all pyramid partitions can be represented by its perfect matchings.

The  shown graph is an example of a rail-yard graph. It is a rail-yard graph given by the word $\{L+,R+,\dots, L+,R-,L-,\dots,R-\}$. In general, a rail-yard graph is characterized by a word from the four-letter alphabet $\{L+,L-,R+,R-\}$. 
Each letter corresponds to a building block -- a column consisting of horizontal and diagonal edges. Building blocks differ among themselves  in the location of diagonal edges only. For example, for the $R+$ block diagonal edges are on the right of the column and  going up (from left to right).

Rail-yard graphs were defined in \cite{bbccr}, and the formulas to compute the partition functions of pure dimer coverings on such graphs were also proved in \cite{bbccr}. Special cases of rail-yard graphs include the Aztec diamond (\cite{EKLP92a,EKLP92b,KJ05,bk}), pyramid partition (\cite{Ken05,Sze08,BY09,you10}), steep tiling (\cite{BCC17}), tower graph (\cite{BF15}),  contracting square-hexagon lattice (\cite{BL17,ZL18,ZL20,Li182}), and contracting bipartite graph (\cite{ZL202,ZL211}).

Pure dimer coverings on rail-yard graphs are in one-to-one correspondence with sequences of partitions. To establish the correspondence one recognize each of the red dots in the graph above as either a particle or hole. Red dots that are the left-end points of an edge are holes and the right-end points are particles. That way on each vertical line we obtain a Maya diagram of a partition. The fundamental pyramid partition corresponds to the sequence of empty partitions, and the pyramid partition with 3 bricks removed shown above to $\emptyset, \emptyset,\emptyset,(1),(1),(1),\emptyset,\dots$

Furthermore, using the correspondence, certain random dimer models on rail-yard graphs can be seen as the probability distribution on sequences of partitions known as the Macdonald process. To study dimer configurations on the rail-yard graphs we need to use dual partitions and we will be dealing with a generalized Macdonald process that allow dual interalacing. Such Macdonald process can also be obtained from the Macdonald processes defined in \cite{bc13,bcgs13} by certain specializations,  which are homomorphisms from the algebra of symmetric polynomials to $\CC$, but not function evaluations. In this paper our asymptotics is concerned only with measure that belong to the subclass of Macdonald processes known as Schur processes, hence in later sections we assume $q=t$. 

Our main technical result which allows us to later perform the asymptotic analysis is done in Section \ref{sect:mh}. The height function of a pyramid partition is naturally defined by its 3-D depiction, but a notion of the height function exists for general rail-yard graphs, not just pyramid partitions. It is the Thurston height function, which is well defined for perfect matchings of bipartite planar graphs. In Section \ref{sect:mh} our main result is the formula for the expectation of the moments for quantities associated with the height function. We refer to these quantities as Macdonald observables, as they can be computed within the framework of Macdonald processes. This kind of approach in studying random models using Macdonald processes  was pioneered in \cite{bc13,bcgs13}, and applied to study the asymptotics of lozenge tilings in \cite{DE16,Ah18}. Our method is very similar to \cite{DE16,Ah18} and makes use of Negut operators. 

In Section \ref{sect:as} we study the asymptotics of the moments of the observables in the appropriate scaling limit. The general result is given in Theorem \ref{t57}. We also show the Gaussian fluctuations in Theorem \ref{t58}. In Section \ref{sect:fb}, we prove an integral formula for the Laplace transform of the rescaled height function (see Theorem \ref{t61}), which turns out to be deterministic, as a 2D analog of the law of large numbers. Before we state the result for pyramid partitions here, we explain the limiting regime. 

We take a sequence of pyramid partitions and scale them so that their corresponding rail-yard graphs have a fixed set of transition points $V_0<V_1<V_2$, representing  the abscissas of the vertical lines of the left, mid (transition from $+$s to $-$s), and  right boundary. For example, if $V_0=-2, V_1=0$, and $V_2=1$, then rail-yard graphs associated with $(L+,R+)^{2n}(L-,R-)^n$ have transition points at  $-8n, 0, 4n$ which can be after scaling by $\epsilon=1/(4n)$ brought to $V_0,V_1,V_2$. The random model we study is what we refer to the periodic $q$-volume model which in the 2-periodic case, such as pyramid partitions, depend on parameters $\tau_1$ and $\tau_2$ and $q=e^{-\epsilon}$, where $\epsilon \to 0$. The weights of diagonal edges are products of weights each depending on one on theses parameters. The weight coming from $q$ is the $q$-volume which is the $q$ analog of the uniform measure on plane partitions or more generally sequences of interlacing partitions. Weights coming from $\tau$s are periodic weights which give different weight to diagonal edges in columns associated with $L$ from those associated with $R$. For general rail-yard graphs the precise conditions on the periodicity of the graph and weights is given in Assumptions \ref{ap5}. The asymptotics of the pyramid partitions for uniform weights, i.e. when $\tau_1=\tau_2$,  was studied in \cite{BBV15}.
Theorem \ref{t61} in case of pyramid partitions says:
\begin{theorem} The rescaled random height function of pyramid partitions $\epsilon h\left(\frac{\chi}{\epsilon},\frac{\kappa}{\epsilon}\right)$ converges, as $\epsilon\rightarrow 0$, to a non-random function $\mathcal{H}(\chi,\kappa)$ such that the Laplace transform of $\mathcal{H}(\chi,\cdot)$ is given by
\begin{align*}
\int_{-\infty}^{\infty} e^{-2\alpha \kappa}\mathcal{H}(\chi,\kappa)d\kappa
=\frac{1}{4\alpha^2\pi\mathbf{i}}\oint_{\mathcal{C}} \left[\mathcal{G}_{\chi}(w)\right]^{\alpha}\frac{dw}{w},
\end{align*}
\begin{align*}
\mathcal{G}_{\chi}(w)=\frac{\left(1+e^{-\chi}w \tau_2\right)\left(1-e^{-V_0}w\tau_1\right)\left(1+e^{-V_1}w \tau_2 \right)\left(1-e^{-V_2}w\tau_1\right)}{{(1-e^{-\chi}w \tau_1})\left(1+e^{-V_0}w\tau_2\right)\left(1-e^{-V_1}w\tau_1 \right)\left(1+e^{-V_2}w\tau_2\right)},
\end{align*}
where $\alpha$ is a positive real number and $\mathcal{C}$ is a positively oriented contour that encloses $-e^{-V_0}\tau_2$, 0, and $e^{-V_1}\tau_1$, but no other poles or zeros of $\mathcal{G}_{\chi}$.
\end{theorem}

\begin{figure}[!h]
\includegraphics[trim=150 100 100 100, clip,width=.4\textwidth]{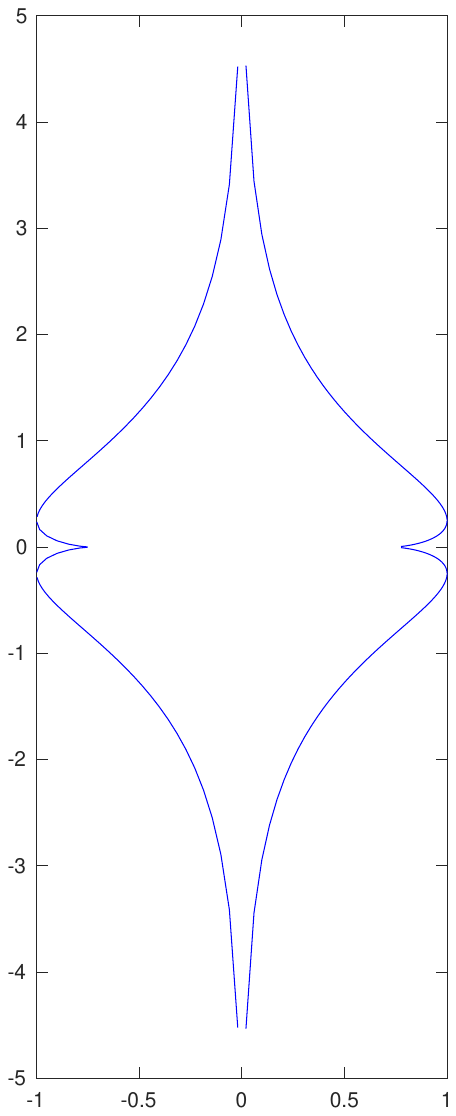}
\includegraphics[trim=160 80 150 120, clip,width=.31\textwidth]{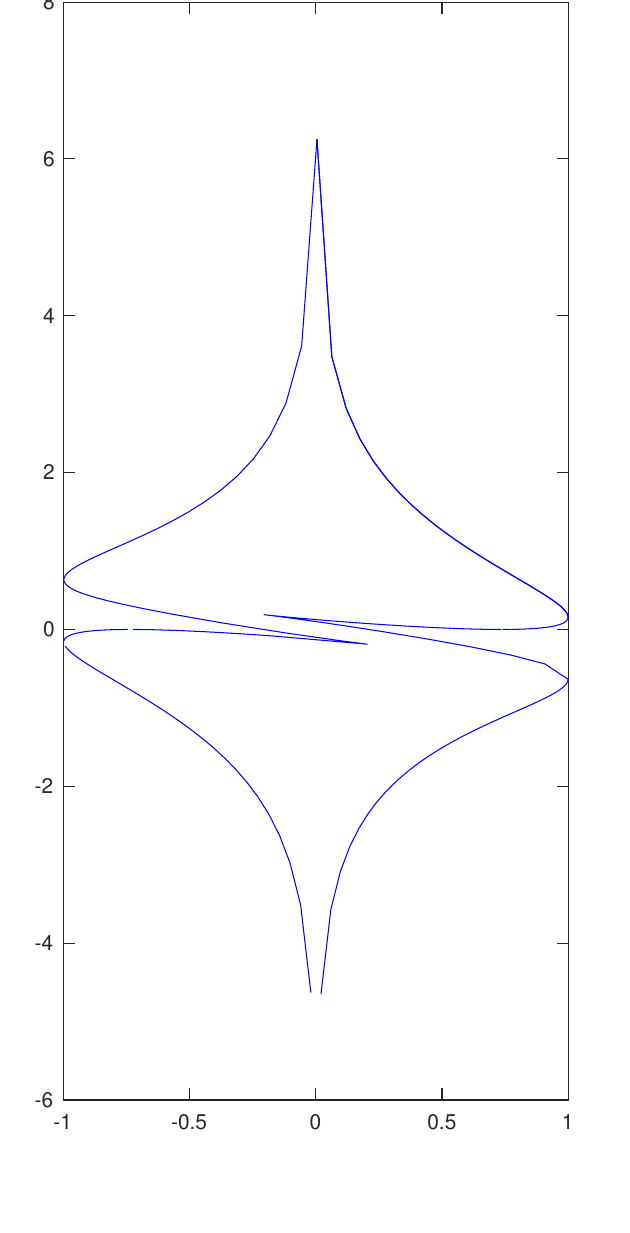}
\caption{Frozen boundary of pyramid partitions with transition points $V_0=-1$, $V_1=0$, $V_2=1$. The left graph has $\tau_1=\tau_2=1$, the right graph has $\tau_1=10,\tau_2=1/10$.}\label{limshapePP}
\end{figure}

The limit shape of pyramid partitions is described as a solution of the parametric equation (parametrized by $w$):
\begin{align*}
\begin{cases}
\mathcal{G}_{\chi}(w)=e^{-2\kappa}\\
\mathcal{G}_{\chi}'(w)=0.
\end{cases}
\end{align*}
The limit shape for pyramid partitions is shown in Figure \ref{limshapePP}. The figure on the left corresponds to the uniform case, and coincides to one obtained by \cite{BBV15}, and the one on the right to a non-uniform case. 

In Section \ref{sect:gff}, we show that the fluctuations of unrescaled height functions converge to the pull-back Gaussian free field (GFF) in the upper half plane under a diffeomorphism from the liquid region to the upper half plane. This result is given in Theorem \ref{t77}, and we state the result for pyramid partitions here. Let $\mathbf{w}_+:\mathcal{L}\rightarrow \HH$ be the diffeomorphism which maps each point $(\chi,\kappa)$ in the liquid region $\mathcal{L}$ to the unique root of $\mathcal{G}_{\chi}(w)=e^{-2\kappa}$ in the upper half plane $\HH$. We discuss in Section \ref{sect:gff} conditions that need to be satisfied so that such map is well defined. Let $\Xi$ be the Gaussian free field (GFF)
on $\HH$ with the zero boundary condition. Then Theorem \ref{t77} for pyramid partitions says:
\begin{theorem}
As $\epsilon\rightarrow 0$, the height function of pyramid partitions converges to the $\mathbf{w}_+$-pullback of GFF in the sense that for any $(\chi,\kappa)\in \mathcal{L}$, $\chi\notin\{V_0,V_1,V_2\}$ and positive real number $\alpha$
\begin{align*}
\int_{-\infty}^{\infty}\left(h\left(\frac{\chi}{\epsilon},\frac{\kappa}{\epsilon}\right)-
\mathbb{E}\left[h\left(\frac{\chi}{\epsilon},\frac{\kappa}{\epsilon}\right)
\right]\right)
e^{-\alpha \kappa}d\kappa\longrightarrow\int_{(\chi,\kappa)\in\mathcal{L}}e^{-\alpha\kappa} \Xi(\mathbf{w}_+(\chi,\kappa)) d\kappa
\end{align*}
in distribution.
\end{theorem}

The organization of the paper is as follows. In Section \ref{sect:bk}, we define the rail-yard graph, the perfect matching and the height function, and review related technical facts. In Section \ref{sect:sc}, we discuss a class of Macdonald processes related to the probability measure of perfect matchings on the rail-yard graphs. In Section \ref{sect:mh}, we compute the moments of height functions of perfect matchings on rail-yard graphs by computing the observables in the generalized Macdonald processes (see Lemma \ref{l36}). In Section \ref{sect:as}, we study the asymptotics of the moments of the random height functions and prove their Gaussian fluctuations in the scaling limit (see Theorems \ref{t57} and \ref{t58}). In Section \ref{sect:fb}, we prove an integral formula for the Laplace transform of the rescaled height function (see Theorem \ref{t61}), which turns out to be deterministic, as a 2D analog of the law of large numbers. We further obtain a parametric equation for the frozen boundary in the scaling limit. In Section \ref{sect:gff}, we prove that the fluctuations of unrescaled height functions converge to the pull-back Gaussian free field (GFF) in the upper half plane under a diffeomorphism from the liquid region to the upper half plane (see Theorem \ref{t77}). In Section \ref{sect:ex}, we discuss specific examples of the rail-yard graphs, where the limit shapes and height fluctuations of perfect matchings can be obtained by the main results in the paper; these examples include the pure steep tilings and pyramid partitions. In Appendix \ref{sc:dmp} we review some facts about Macdonald polynomials and include some known technical results.

\section{Backgrounds}\label{sect:bk}

In this section, we define the rail-yard graph, the perfect matching and the height function, and review related technical facts.

\subsection{Weighted rail-yard graphs}

Let $l,r\in\ZZ$ such that $l\leq r$. Let $[l..r]:=[l,r]\cap\ZZ,$ i.e., $[l..r]$ is the set of integers between $l$ and $r$. For a positive integer $m$, we use $[m]:=\{1,2,\ldots,m\}.$

Consider two binary sequences indexed by integers in $[l..r]$ 
\begin{itemize}
\item the $LR$ sequence $\underline{a}=\{a_l,a_{l+1},\ldots,a_r\}\in\{L,R\}^{[l..r]}$;
\item the sign sequence $\underline{b}=(b_l,b_{l+1},\ldots,b_l)\in\{+,-\}^{[l..r]}$.
\end{itemize}
The rail-yard graph $RYG(l,r,\underline{a},\underline{b})$ with respect to integers $l$ and $r$, the $LR$ sequence $\underline{a}$ and the sign sequence $\underline{b}$, is the bipartite graph with vertex set $[2l-1..2r+1]\times \left\{\ZZ+\frac{1}{2}\right\}$. A vertex is called even (resp.\ odd) if its abscissa is an even (resp.\ odd) integer. Each even vertex $(2m,y)$, $m\in[l..r]$ is incident to 3 edges, two horizontal edges joining it to the odd vertices $(2m-1,y)$ and $(2m+1,y)$ and one diagonal edge joining it to
\begin{itemize}
\item the odd vertex $(2m-1,y+1)$ if $(a_m,b_m)=(L,+)$;
\item the odd vertex $(2m-1,y-1)$ if $(a_m,b_m)=(L,-)$;
\item the odd vertex $(2m+1,y+1)$ if $(a_m,b_m)=(R,+)$;
\item the odd vertex $(2m+1,y-1)$ if $(a_m,b_m)=(R,-)$.
\end{itemize} 

\begin{figure}
\includegraphics[width=.7\textwidth]{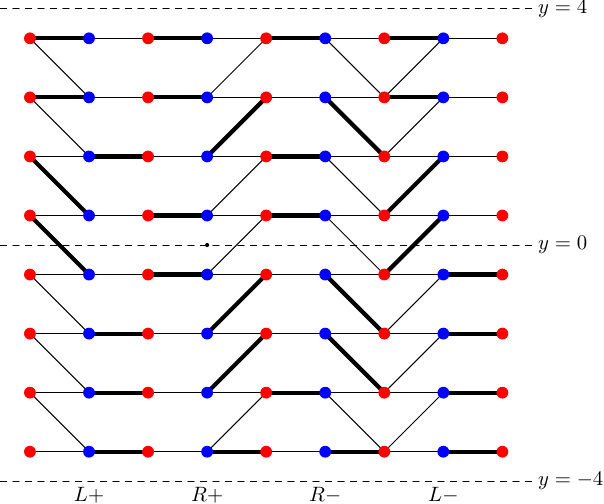}
\caption{A rail yard graph with LR sequence $\underline{a}=\{L,R,R,L\}$and sign sequence $\underline{b}=\{+,+,-,-\}$. Odd vertices are represented by red points, and even vertices are represented by blue points. Dark lines represent a pure dimer covering. Assume that above the horizontal line $y=4$, only horizontal edges with an odd vertex on the left are present in the dimer configuration; and below the horizontal line $y=-4$, only horizontal edges with an even vertex on the left are present in the dimer configuration. The corresponding sequence of partitions (from the left to the right) is given by $\emptyset\prec(2,0,\ldots)\prec' (3,1,1,\ldots)\succ'(2,0,\ldots)\succ \emptyset$.}\label{fig:rye}
\end{figure}

See Figure \ref{fig:rye} for an example of a rail-yard graph. 

The left boundary (resp.\ right boundary) of $RYG(l,r,\underline{a},\underline{b})$ consists of all odd vertices with abscissa $2l-1$ (resp.\ $2r+1$). Vertices which do not belong to the boundaries are called inner. A face of $RYG(l,r,\underline{a},\underline{b})$ is called an inner face if it contains only inner vertices.

We assign edge weights to a rail yard graph $RYG(l,r,\underline{a},\underline{b})$ as follows:
\begin{itemize}
    \item a horizontal edge has weight 1; and
    \item a diagonal edge adjacent to a vertex with abscissa $2i$ has weight $x_i$.
\end{itemize}

\subsection{Dimer coverings and pure dimer coverings}

\begin{definition}
A dimer covering is a subset of edges of $RYG(l,r,\underline{a},\underline{b})$ such that
\begin{enumerate}
\item each inner vertex of $RYG(l,r,\underline{a},\underline{b})$ is incident to exactly one edge in the subset;
\item each left boundary vertex or right boundary vertex is incident to at most one edge in the subset;
\item only a finite number of diagonal edges are present in the subset.
\end{enumerate}
A pure dimer covering of $RYG(l,r,\underline{a},\underline{b})$ is dimer  covering of $RYG(l,r,\underline{a},\underline{b})$ satisfying the following two additional conditions
\begin{itemize}
\item each left boundary vertex $(2l-1,y)$ is incident to exactly one edge (resp.\ no edges) in the subset if $y>0$ (resp.\ $y<0$).
\item each right boundary vertex $(2r+1,y)$ is incident to exactly one edge (resp.\ no edges) in the subset if $y<0$ (resp.\ $y>0$).
\end{itemize}
\end{definition}

See Figure \ref{fig:rye} for an example of pure dimer coverings on a rail yard graph.

For a dimer covering $M$ on the rail-yard graph $RYG(l,r,\underline{a},\underline{b})$, define the associated height function $h_{M}$ on faces of $RYG(l,r,\underline{a},\underline{b})$ as follows. We first define a preliminary height function $\overline{h}_M$ on faces of $RYG(l,r,\underline{a},\underline{b})$. Note that there exists a positive integer $N>0$, such that when $y<-N$, only horizontal edges with even vertices on the left are present.  Fix a face $f_0$ of $RYG(l,r,\underline{a},\underline{b})$ such that the midpoint of $f_0$ is on the horizontal line $y=-N$, and define $\overline{h}_M(f_0)=0.$

For any two adjacent faces $f_1$ and $f_2$ sharing at least one edge, 
\begin{itemize}
\item If moving from $f_1$ to $f_2$ crosses a present (resp.\ absent) horizontal edge in $M$ with odd vertex on the left, then $\ol{h}_M(f_2)-\ol{h}_M(f_1)=1$ (resp.\ $\ol{h}_M(f_2)-\ol{h}_M(f_1)=-1$).
\item If moving from $f_1$ to $f_2$ crosses a present (resp.\ absent) diagonal edge in $M$ with odd vertex on the left, then $\ol{h}_M(f_2)-\ol{h}_M(f_1)=2$ (resp.\ $\ol{h}_M(f_2)-\ol{h}_M(f_1)=0$).
\end{itemize}

Let $\ol{h}_0$ be the preliminary height function associated to the dimer configuration satisfying 
\begin{itemize}
\item no diagonal edge is present; and
\item each present edge is horizontal with an even vertex on the left.
\end{itemize}
The height function $h_{M}$ associated to $M$ is then defined by 
\begin{align}
h_M=\ol{h}_M-\ol{h}_0.\label{dhm}
\end{align}

Let $m\in[l..r]$. Let $x=2m-\frac{1}{2}$ be a vertical line such that all the horizontal edges and diagonal edges of $RYG(l,r,\underline{a},\underline{b})$ crossed by $x=2m-\frac{1}{2}$ have odd vertices on the left. Then for each point $\left(2m-\frac{1}{2},y\right)$ in a face of $RYG(l,r,\underline{a},\underline{b})$, we have
\begin{align}
h_{M}\left(2m-\frac{1}{2},y\right)=2\left[N_{h,M}^{-}\left(2m-\frac{1}{2},y\right)+N_{d,M}^{-}\left(2m-\frac{1}{2},y\right)\right];\label{hm1}
\end{align}
where $N_{h,M}^{-}\left(2m-\frac{1}{2},y\right)$ is the total number of present horizontal edges in $M$ crossed by $x=2m-\frac{1}{2}$ below $y$, and $N_{d,M}^{-}\left(2m-\frac{1}{2},y\right)$ is the total number of present diagonal edges in $M$ crossed by $x=2m-\frac{1}{2}$ below $y$. From the definition of a pure dimer covering we can see that both $N_{h,M}^{-}\left(2m-\frac{1}{2},y\right)$ and $N_{d,M}^{-}\left(2m-\frac{1}{2},y\right)$ are finite for each finite $y$.

Note also that $x=2m+\frac{1}{2}$ is a vertical line such that all the horizontal edges and diagonal edges of $RYG(l,r,\underline{a},\underline{b})$ crossed by $x=2m+\frac{1}{2}$ have even vertices on the left. Then for each point $\left(2m+\frac{1}{2},y\right)$ in a face of $RYG(l,r,\underline{a},\underline{b})$, we have
\begin{align}
h_{M}\left(2m+\frac{1}{2},y\right)=2\left[J_{h,M}^{-}\left(2m+\frac{1}{2},y\right)-N_{d,M}^{-}\left(2m+\frac{1}{2},y\right)\right];\label{hm2}
\end{align}
where $J_{h,M}^{-}\left(2m+\frac{1}{2},y\right)$ is the total number of absent horizontal edges in $M$ crossed by $x=2m+\frac{1}{2}$ below $y$, and $N_{d,M}^{-}\left(2m+\frac{1}{2},y\right)$ is the total number of present diagonal edges in $M$ crossed by $x=2m+\frac{1}{2}$ below $y$. From the definition of a pure dimer covering we can also see that both $J_{h,M}^{-}\left(2m+\frac{1}{2},y\right)$ and $N_{d,M}^{-}\left(2m+\frac{1}{2},y\right)$ are finite for each finite $y$.

\subsection{Partitions}

A partition is a non-increasing sequence $\lambda=(\lambda_i)_{i\geq 0}$ of non-negative integers which vanish eventually. Let $\mathbb{Y}$ be the set of all the partitions. The size of a partition is defined by $|\lambda|=\sum_{i\geq 1}\lambda_i.$
Two partitions $\lambda$ and $\mu$ are called interlaced, and written by $\lambda\succ\mu$ or $\mu\prec \lambda$ if $\lambda_1\geq \mu_1\geq\lambda_2\geq \mu_2\geq \lambda_3\cdots.$ When representing partitions by Young diagrams, this means $\lambda/\mu$ is a horizontal strip. The conjugate partition $\lambda'$ of $\lambda$ is a partition whose Young diagram $Y_{\lambda'}$ is the image of the Young diagram $Y_{\lambda}$ of $\lambda$ by the reflection along the main diagonal. More precisely
\begin{align*}
\lambda_i':=\left|\{j\geq 0: \lambda_j\geq i\}\right|,\qquad \forall i\geq 1.
\end{align*}

The skew Schur functions are defined in Section I.5 of \cite{IGM15}.

\begin{definition}\label{dss}Let $\lambda$, $\mu$ be partitions. Define the skew Schur functions as 
\begin{align*}
s_{\lambda/\mu}=\det\left(h_{\lambda_i-\mu_j-i+j}\right)_{i,j=1}^{l(\lambda)}
\end{align*} 
where for $r<0$, $h_r=0$ and for $r\geq 0$, $h_r$ is the $r$th complete symmetric function defined by the sum of all monomials of total degree $r$ in the variables $x_1,x_2,\ldots$. More precisely,
\begin{align*}
h_r=\sum_{1\leq i_1\leq i_2\leq \ldots\leq i_r} x_{i_1}x_{i_2}\cdots x_{i_r}
\end{align*}
Define the Schur function as $s_{\lambda}=s_{\lambda/\emptyset}.$
\end{definition}

For a dimer covering $M$ of $RYG(l,r,\underline{a},\underline{b})$, we associate a particle-hole configuration to each odd vertex of $RYG(l,r,\underline{a},\underline{b})$ as follows: let $m\in[l..(r+1)]$ and $k\in\ZZ$: if the odd endpoint $\left(2m-1,k+\frac{1}{2}\right)$ is incident to a present edge in $M$ on its right (resp.\ left), then associate a hole (resp.\ particle) to the odd endpoint $\left(2m-1,k+\frac{1}{2}\right)$. When $M$ is a pure dimer covering, it is not hard to check that there exists $N>0$, such that when $y>N$, only holes exist and when $y<-N$, only particles exist.

We associate a partition $\lambda^{(M,m)}$ to the column indexed by $m$ of particle-hole configurations, which corresponds to a pure dimer covering $M$ adjacent to odd vertices with abscissa $(2m-1)$ as follows. Assume
\begin{align*}
\lambda^{(M,m)}=(\lambda^{(M,m)}_1,\lambda^{(M,m)}_2,\ldots),
\end{align*}
Then for $i\geq 1$, $\lambda^{(M,m)}_i$ is the total number of holes in $M$ along the vertical line $x=2m-1$ below the $i$th highest particles. Let $l(\lambda^{(M,m)})$ be the total number of nonzero parts in the partition $\lambda^{(M,m)}$.

We define the charge $c^{(M,m)}$ on column $(2m-1)$ for the configuration $M$ as follows:
\begin{align}
c^{(M,m)}&=&\mathrm{number \ of\ particles\ on\ column\ }(2m-1)\ \mathrm{in\ the\ upper\ half\ plane}\notag\\
&&-\mathrm{number\ of\ holes\ on\ column\ }(2m-1)\ \mathrm{in\ the\ lower\ half\ plane}\label{dcg}
\end{align}

The weight of a dimer covering $M$ of $RYG(l,r,\underline{a},\underline{b})$ is defined as follows
\begin{align*}
w(M):=\prod_{i=l}^{r}x_i^{d_i(M)},
\end{align*}
where $d_i(M)$ is the total number of present diagonal edges of $M$ incident to an even vertex with abscissa $2i$. 

Let $\lambda^{(l)},\lambda^{(r+1)}$ be two partitions. The partition function $Z_{\lambda^{(l)},\lambda^{(r+1)}}(G,\underline{x})$ of dimer coverings on $RYG(l,r,\underline{a},\underline{b})$ whose configurations on the left (resp.\ right) boundary correspond to partition $\lambda^{(l)}$ (resp.\ $\lambda^{(r+1)}$) is the sum of weights of all such dimer coverings on the graph. Given the left and right boundary conditions $\lambda^{(l)}$ and $\lambda^{(r+1)}$, respectively, the probability of a dimer covering $M$ is then defined by
\begin{align}
\mathrm{Pr}(M|\lambda^{(l)},\lambda^{(r+1)}):=\frac{w(M)}{Z_{\lambda^{(l)},\lambda^{(r+1)}}(G,\underline{x})}.\label{ppd}
\end{align}
Note that pure dimer coverings have left and right boundary conditions given by
\begin{align}
\lambda^{(l)}=\lambda^{(r+1)}=\emptyset \label{pbc}.
\end{align}

\begin{figure}
\includegraphics[width=.18\textwidth]{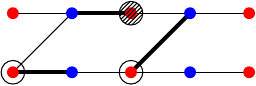}\quad\raisebox{2ex}{$\longleftrightarrow$}\quad\includegraphics[width=.18\textwidth]{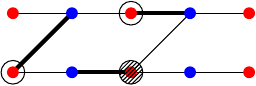} \hfill \includegraphics[width=.18\textwidth]{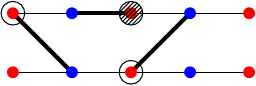}\quad\raisebox{2ex}{$\longleftrightarrow$}\quad\includegraphics[width=.18\textwidth]{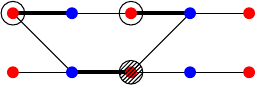}\\
\bigskip
\includegraphics[width=.18\textwidth]{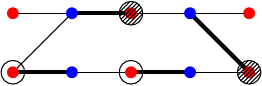}\quad\raisebox{2ex}{$\longleftrightarrow$}\quad\includegraphics[width=.18\textwidth]{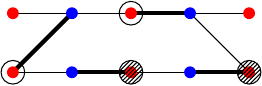}\hfill \includegraphics[width=.18\textwidth]{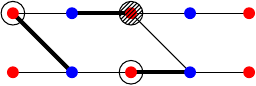}\quad\raisebox{2ex}{$\longleftrightarrow$}\quad\includegraphics[width=.18\textwidth]{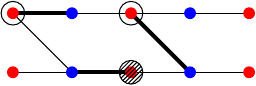}\\
\bigskip
\includegraphics[width=.18\textwidth]{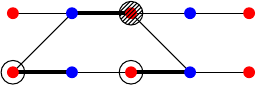}\quad\raisebox{2ex}{$\longleftrightarrow$}\quad\includegraphics[width=.18\textwidth]{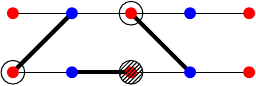}\hfill
\includegraphics[width=.18\textwidth]{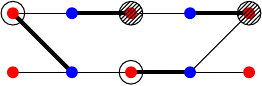}\quad\raisebox{2ex}{$\longleftrightarrow$}\quad\includegraphics[width=.18\textwidth]{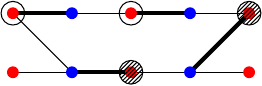}\\
\bigskip
\includegraphics[width=.18\textwidth]{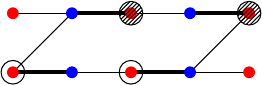}\quad\raisebox{2ex}{$\longleftrightarrow$}\quad\includegraphics[width=.18\textwidth]{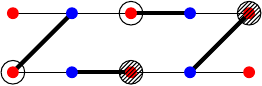}\hfill\includegraphics[width=.18\textwidth]{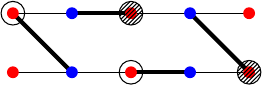}\quad\raisebox{2ex}{$\longleftrightarrow$}\quad\includegraphics[width=.18\textwidth]{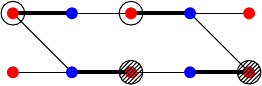}\\
\bigskip
\includegraphics[width=.18\textwidth]{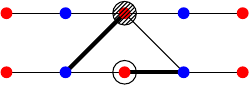}\quad\raisebox{2ex}{$\longleftrightarrow$}\quad\includegraphics[width=.18\textwidth]{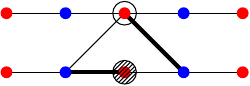}\hfill\includegraphics[width=.18\textwidth]{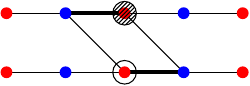}\quad\raisebox{2ex}{$\longleftrightarrow$}\quad\includegraphics[width=.18\textwidth]{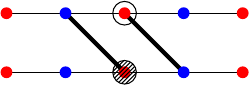}\\
\bigskip
\includegraphics[width=.18\textwidth]{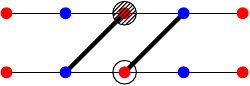}\quad\raisebox{2ex}{$\longleftrightarrow$}\quad\includegraphics[width=.18\textwidth]{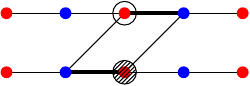}\hfill\includegraphics[width=.18\textwidth]{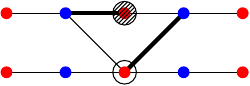}\quad\raisebox{2ex}{$\longleftrightarrow$}\quad\includegraphics[width=.18\textwidth]{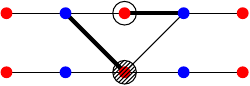}\\
\bigskip
\includegraphics[width=.18\textwidth]{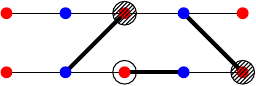}\quad\raisebox{2ex}{$\longleftrightarrow$}\quad\includegraphics[width=.18\textwidth]{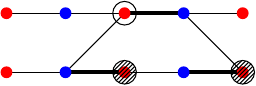}\hfill\includegraphics[width=.18\textwidth]{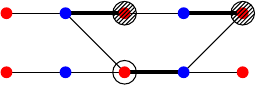}\quad\raisebox{2ex}{$\longleftrightarrow$}\quad\includegraphics[width=.18\textwidth]{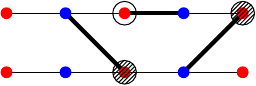}\\
\bigskip
\includegraphics[width=.18\textwidth]{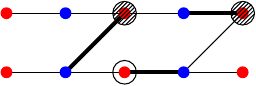}\quad\raisebox{2ex}{$\longleftrightarrow$}\quad\includegraphics[width=.18\textwidth]{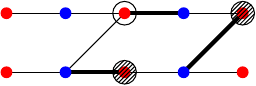}\hfill\includegraphics[width=.18\textwidth]{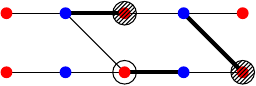}\quad\raisebox{2ex}{$\longleftrightarrow$}\quad\includegraphics[width=.18\textwidth]{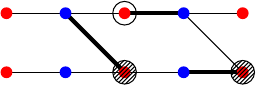}
\caption{Flip of dimer configurations on a face between two columns. Left: $(L-,L-)$, $(L-,R-)$, $(L-,L+)$,$(L-,R+)$ $(R+,L+)$, $(R+,L-)$, $(R+,R-)$, $(R+,R+)$. Right: $(L+,L-)$, $(L+,L+)$, $(L+,R+)$, $(L+,R-)$, $(R-,L+)$, $(R-,L-)$, $(R-,R+)$, $(R-,R-)$. }\label{fig:flip1234}
\end{figure}

Let $f$ be an inner face of $RYG(l,r,\underline{a},\underline{b})$. Let $M$ be a dimer covering of $C$. If exactly half of the edges bordering $f$ are present in $M$, we can obtain another dimer covering $M'$ from $M$, such that $M'$ and $M$ coincide on each edge not bordering $f$; while for an edge bordering $f$, it is present in $M'$ if and only if it is absent in $M$. In particular, $M$ and $M'$ have the same configuration on the left and right boundary. The operation of replacing $M$ by $M'$ is called a flip of $f$; see Figure \ref{fig:flip1234}, where odd vertices are represented by red dots, even vertices are represented by blue dots.

Then we have the following lemma.

\begin{lemma}\label{l01}Let $M$ be a pure dimer covering on the rail-yard graph $RYG(l,r,\underline{a},\underline{b})$. Then 
\begin{align*}
c^{(M,m)}=0,\ \forall m\in[l..(r+1)].
\end{align*}
\end{lemma}

\begin{proof}Let $M_0$ be the pure dimer covering on $RYG(l,r,\underline{a},\underline{b})$ such that
\begin{itemize}
\item all the present edges in the upper half plane are horizontal with odd vertex on the left; and
\item all the present edges in the lower half plane are horizontal with even vertex on the left.
\end{itemize}
It is straightforward to check that in the particle-whole representation for any column in $M_0$, the upper half plane only has holes, while the lower half plane only has particles. By (\ref{dcg}) we obtain
\begin{align*}
c^{(M_0,m)}=0,\qquad \forall m\in[l..(r+1)].
\end{align*}
By Section 2.3 of \cite{bbccr} (see also \cite{jp02}), any pure dimer covering $M$ of $RYG(l,r,\underline{a},\underline{b})$ can be obtained from $M_0$ by finitely many flips.
The particle-hole configuration is associated to each odd vertex.  The particle-hole configuration for each type of a flip is shown in Figure \ref{fig:flip1234}, where particles are represented by hatched circles, while holes are represented by non-hatched circles. Each local particle-hole configuration is lying in two adjacent rows. The following cases might occur
\begin{itemize}
\item both rows are in the upper half plane; or
\item both rows are in the lower half plane; or
\item the top row is in the upper half plane and the bottom row is in the lower half plane.
\end{itemize}
It is straightforward to check that for each one of the 3 cases above, and each type particle-hole configuration the charge $c^{M,m}$ for all $m\in[l..r+1]$ remains unchanged. Then for any pure dimer covering $M$, $c^{(M,m)}=c^{(M_0,m)}=0.$ Then the lemma follows.
\end{proof}

\subsection{Asymptotic height function.}

Let $M$ be a dimer covering of $RYG(l,r,\underline{a},\underline{b})$. Let $Y^{(i,M,m)}$ be the ordinate of the $i$th highest particle along the line $x=2m-1$ for the pure dimer covering $M$. Then by (\ref{dcg}) we obtain
\begin{align}
\lambda_i^{(M,m)}=Y^{(i,M,m)}-c^{(M,m)}+i-\frac{1}{2}.\label{dny}
\end{align}

Assume $k\log t<0$, we have
\begin{align}
\label{ih}\int_{-\infty}^{\infty}h_{M}(x,y)t^{ky}dy&=\frac{1}{k\log t}\int_{-\infty}^{\infty}h_M(x,y)\frac{d e^{ky\log t}}{dy}dy\\
&=-\frac{1}{k\log t}\int_{-\infty}^{\infty}e^{ky\log t}\frac{d h_M(x,y)}{dy}dy.\notag
\end{align}
Let $x=2m-\frac{1}{2}$.
From (\ref{hm1}), we obtain
\begin{align}
\frac{d h_M(2m-\frac{1}{2},y)}{dy}=2\left(1-\sum_{i=1}^{l(\lambda^{(M,m)})}\mathbf{1}_{[Y^{(i,M,m)}-\frac{1}{2},Y^{(i,M,m)}+\frac{1}{2}]}(y)\right).\label{drh}
\end{align}
Here for $A\subseteq \RR$, $\mathbf{1}_{A}(y):\RR\rightarrow \{0,1\}$ is the indicator function for the set $A$, i.e., $\mathbf{1}_{A}(y)=1$ if $y\in A$ and 0 otherwise.

By (\ref{dny}) we obtain for $1\leq i\leq l(\lambda^{(M,m)})$
\begin{align}
Y^{(i,M,m)}=\frac{1}{2}+\lambda_i^{(M,m)}-i+c^{(M,m)}.\label{dyi}
\end{align}
Let
\begin{align*}
B_M(m):=Y^{(l(\lambda^{(M,m)})+1)}+\frac{1}{2}.
\end{align*}
Note that below $B_M(m)$, only particles are present along the vertical line $y=2m-1$, hence we have
\begin{align*}
\frac{d h_M(2m-\frac{1}{2},y)}{dy}=0,\qquad \forall y<B_M(m).
\end{align*} 
Moreover, since the charge $c^{(M,m)}=0$, there are exactly the same number of particles on the upper half plane and holes in the lower half plane along the line $x=2m-1$, we obtain
\begin{align}
-B_M(m)&=\mathrm{number\ of\ particles\ at\ } \left(2m-1,y\right)\ \mathrm{with\ }B_M(m)<y<0\notag\\&+
\mathrm{number\ of\ holes\ at\ }\left(2m-1,y\right)\ \mathrm{with}\ B_M(m)<y<0\notag
\\
&=\mathrm{number\ of\ particles\ at\ }(2m-1,y)\ \mathrm{with\ }B_M(m)<y<0\notag\\&+
\mathrm{number\ of\ particles\ at\ } (2m-1,y)\ \mathrm{with}\ y>0\notag\\
&=
l(\lambda^{(M,m)}).\label{b0}
\end{align}
Then from (\ref{ih}) and (\ref{drh}) we obtain
\begin{align*}
&\int_{-\infty}^{\infty}h_{M}(x,y)t^{ky}dy=\\
&=\frac{2}{k\log t}\left[-\int_{B_M(m)}^{\infty}e^{ky\log t}dy+\int_{B_M(m)}^{\infty}\sum_{i=1}^{l(\lambda^{(M,m)})}\mathbf{1}_{[Y^{(i,M,m)}-\frac{1}{2},Y^{(i,M,m)}+\frac{1}{2}]}(y) e^{ky\log t}dy\right]\\
&=\frac{2 t^{kB_M(m)}}{(k\log t)^2}+\frac{2}{(k\log t)^2}\sum_{i=1}^{l(\lambda^{(M,m)})}\left(e^{k(Y^{(i,M,m)}+\frac{1}{2})\log t}-e^{k(Y^{(i,M,m)}-\frac{1}{2})\log t}\right).
\end{align*}
By (\ref{dyi}), we obtain
\begin{align}
&\int_{-\infty}^{\infty}h_{M}(x,y)t^{ky}dy=\notag\\
&=\frac{2t^{k(B_M(m)+l(\lambda^{(M,m)}))}}{(k\log t)^2}\left[t^{-kl(\lambda^{(M,m)})}+(1-t^{-k})\sum_{i=1}^{l(\lambda^{(M,m)})}t^{k(\lambda_i^{(M,m)}+c^{(M,m)}-i+1)}\right]\notag\\
&=\frac{2}{(k\log t)^2}\left[t^{-kl(\lambda^{(M,m)})}+(1-t^{-k})\sum_{i=1}^{l(\lambda^{(M,m)})}t^{k(\lambda_i^{(M,m)}+c^{(M,m)}-i+1)}\right]\label{hme},
\end{align}
where the last identity follows from (\ref{b0}). In particular if $M$ is a pure dimer covering, we have
\begin{align*}
\int_{-\infty}^{\infty}h_{M}(x,y)t^{ky}dy
=\frac{2}{(k\log t)^2}\left[t^{-kl(\lambda^{(M,m)})}+(1-t^{-k})\sum_{i=1}^{l(\lambda^{(M,m)})}t^{k(\lambda_i^{(M,m)}-i+1)}\right].
\end{align*}

The bosonic Fock space $\mathcal{B}$ is the infinite dimensional Hilbert space spanned by the orthonormal basis vectors $|\lambda\rangle$, where $\lambda$ runs over all the partitions. Let $\langle \lambda |$ denote the dual basis vector. Let $x$ be a formal or a complex variable. Introduce the operators $\Gamma_{L+}(x)$, $\Gamma_{L-}(x)$, $\Gamma_{R+}(x)$, $\Gamma_{R-}(x)$ from $\mathcal{B}$  to $\mathcal{B}$ as follows
\begin{align*}
\Gamma_{L+}(x)|\lambda\rangle=\sum_{\mu\prec \lambda}x^{|\lambda|-|\mu|}|\mu\rangle;\qquad \Gamma_{R+}(x)|\lambda\rangle=\sum_{\mu'\prec \lambda'}x^{|\lambda|-|\mu|}|\mu\rangle;\\
\Gamma_{L-}(x)|\lambda\rangle=\sum_{\mu\succ \lambda}x^{|\mu|-|\lambda|}|\mu\rangle;\qquad \Gamma_{R-}(x)|\lambda\rangle=\sum_{\mu'\succ \lambda'}x^{|\mu|-|\lambda|}|\mu\rangle.
\end{align*}
These operators were first introduced by the Kyoto school, and used to study random partitions in \cite{oko01}.

\begin{lemma}\label{l12}Let $a_1,a_2\in \{L,R\}$. We have the following commutation relations for the operators $\Gamma_{a_1,\pm}$, $\Gamma_{a_2,\pm}$
\begin{align*}
\Gamma_{a_1,+}(x_1)\Gamma_{a_2,-}(x_2)=\begin{cases}\frac{\Gamma_{a_2,-}(x_2)\Gamma_{a_1,+}(x_1)}{1-x_1x_2}&\mathrm{if}\ a_1=a_2\\(1+x_1x_2)\Gamma_{a_2,-}(x_2)\Gamma_{a_1,+}(x_1)&\mathrm{if}\ a_1\neq a_2\end{cases}.
\end{align*}
Moreover,
\begin{align*}
\Gamma_{a_1,b}(x_1)\Gamma_{a_2,b}(x_2)=\Gamma_{a_2,b}(x_2)\Gamma_{a_1,b}(x_1);
\end{align*}
for all $a_1,a_2\in\{L,R\}$ and $b\in\{+,-\}$.
\end{lemma}

\begin{proof}See Proposition 7 of \cite{bbccr}; see also \cite{you10,bbb14}.
\end{proof}

Given the definitions of the operators $\Gamma_{a,b}(x)$ with $a\in\{L,R\}$, $b\in\{+,-\}$, it is straightforward to check the following lemma.

\begin{lemma}\label{l13}The partition function of dimer coverings on a rail yard graph $G=RYG(l,r,\underline{a},\underline{b})$ with left and right boundary conditions given by $\lambda^{(l)},\lambda^{(r+1)}$, respectively,  is 
\begin{align}
Z_{\lambda^{(l)},\lambda^{(r+1)}}(G;\underline{x})=\langle\lambda^{(l)}| \Gamma_{a_lb_l}(x_l)\Gamma_{a_{l+1}b_{l+1}}(x_{l+1})\cdots \Gamma_{a_rb_r}(x_r)|\lambda^{(r+1)} \rangle.\label{bf}
\end{align}
\end{lemma}

\begin{corollary}
The partition function of pure dimer coverings can be computed as follows:
\begin{align}
Z_{\emptyset,\emptyset}(G;\underline{x})=\prod_{l\leq i<j\leq r;b_i=+,b_j=-}z_{i,j},\label{fp}
\end{align}
where
\begin{align}
z_{ij}=\begin{cases}1+x_ix_j&\mathrm{if}\ a_i\neq a_j\\\frac{1}{1-x_ix_j}&\mathrm{if}\ a_i=a_j\end{cases}.\label{dzij}
\end{align}
\end{corollary}

\begin{proof}The corollary follows from Lemma \ref{l13} by letting $\lambda^{(l)}=\lambda^{(r+1)}=\emptyset$; it also appears in Proposition 8 of \cite{bbccr} for (\ref{bf}) and Theorem 1 of \cite{bbccr} for (\ref{fp}).
\end{proof}

\noindent{\textbf{Remark.}} The partition function $Z(G;\underline{x})$ is always well-defined as a power series in $\underline{x}$. When we consider the edge weights $x_i$'s to be positive numbers, to make sure the convergence of the power series representing the partition function, we need to assume that for any $i,j\in[l..r]$, $i<j$, $a_i=a_j$ and $b_i=+$, $b_j=-$ we have $x_ix_j<1$. However, when considering the corresponding probability measure, we do not necessarily need this assumption.

\section{Macdonald Processes}\label{sect:sc}

In this section, we discuss a class of Macdonald processes related to the probability measure of perfect matchings on the rail-yard graphs. The major characterstic of the processes defined here is that the processes involve dual partitions as well, which, as we will see, can also be obtained from certain non-function-evaluation specializations of the Macdonald processes defined without dual partitions (see \cite{bc13}), when the parameters satisfy $q=t$.

Let $G=RYG(l,r,\underline{a},\underline{b})$ be a rail-yard graph. Let $(\lambda^{(M,l)},\lambda^{(M,l+1)},\ldots,\lambda^{(M,r+1)})$
be the sequence of partitions corresponding to a dimer covering $M$ on $G$. By  Lemmas \ref{l12} and \ref{l13}, we obtain for $i\in[l..r]$
\begin{enumerate}
\item If $(a_i,b_i)=(L,-)$, $\lambda^{(M,i+1)}\prec \lambda^{(M,i)}$;
\item If $(a_i,b_i)=(L,+)$, $\lambda^{(M,i+1)}\succ \lambda^{(M,i)}$;
\item If $(a_i,b_i)=(R,-)$, $[\lambda^{(M,i+1)}]'\prec [\lambda^{(M,i)}]'$;
\item If $(a_i,b_i)=(R,+)$, $[\lambda^{(M,i+1)}]'\succ [\lambda^{(M,i)}]'$.
\end{enumerate}

Given Definition \ref{dss}, we can express the probability of a dimer covering $M$ conditional on the left and right boundary conditions $\lambda^{(l)}$ and $\lambda^{(r+1)}$ respectively, as defined by (\ref{ppd}), as follows:
\begin{align}
\mathrm{Pr}&(M|\lambda^{(l)},\lambda^{(r+1)}):=\frac{1}{Z_{\lambda^{(l)},\lambda^{(r+1)}}(G,\underline{x})}\notag\\
&\prod_{\substack{i\in[l..r]\\(a_i,b_i)=(L,-)}} s_{\lambda^{(M,i)}/\lambda^{(M,i+1)}}(x_i)\prod_{\substack{j\in[l..r]\\(a_i,b_i)=(L,+)}} s_{\lambda^{(M,j+1)}/\lambda^{(M,j)}}(x_j)\notag\\
&\prod_{\substack{i\in[l..r]\\(a_i,b_i)=(R,-)}} s_{[\lambda^{(M,i)}]'/[\lambda^{(M,i+1)}]'}(x_i)\prod_{\substack{j\in[l..r]\\(a_i,b_i)=(R,+)}} s_{[\lambda^{(M,j+1)}]'/[\lambda^{(M,j)}]'}(x_j).\label{pm}
\end{align}

Now we define a generalized Macdonald process, which is a formal probability measure on sequences of partitions such that the probability of each sequence of partitions is proportional to a sum of products of skew Macdonald polynomials. See Section \ref{sc:dmp} for definitions of Macdonald polynomials $P_{\lambda}$, $Q_{\lambda}$, $P_{\lambda/\mu}$, $Q_{\lambda/\mu}$.
\begin{definition}\label{df22}Let $\bA=(A^{(l)},\ldots,A^{(r+1)})$ and $\bB=(B^{(l+1)},\ldots,B^{(r+1)})$ be $2(r-l+1)$ set of variables, in which each $A^{(i)}$ or $B^{(j)}$ consists of countably many variables.
Let $\mathcal{P}=\{\mathcal{L},\mathcal{R}\}$ be a partition of the set $[l..r]$, i.e. $\mathcal{L}\cup\mathcal{R}=[l..r]$ and $\mathcal{L}\cap \mathcal{R}=\emptyset$.

Define a formal probability measure on the set of sequences of $(r-l+2)$ partitions $(\lambda^{(l)},\lambda^{(r+1)},\ldots,\lambda^{(r+1)})$ with respect to $\mathcal{P}$, $\bA$ and $\bB$ and parameters $q,t\in(0,1)$ by
\begin{align}
\MP_{\bA,\bB,\mathcal{P},q,t}(\lambda^{(l)},\ldots,\lambda^{(r+1)})&\propto \left[\prod_{i\in\mathcal{L}}\Psi_{\lambda^{(i)},\lambda^{(i+1)}}(A^{(i)},B^{(i+1)};q,t)\right] \notag\\
&\times\left[\prod_{j\in\mathcal{R}}\Phi_{[\lambda^{(j)}]',[\lambda^{(j+1)}]'}(A^{(j)},B^{(j+1)};q,t)\right], \label{pmd}
\end{align}
where for two partitions $\lambda,\mu\in \YY$, and two countable set of variables $A$, $B$,
\begin{align*}
\Psi_{\lambda,\mu}(A,B;q,t)=\sum_{\nu\in \YY}P_{\lambda/\nu}(A;q,t)Q_{\mu/\nu}(B;q,t),\\
\Phi_{\lambda,\mu}(A,B;q,t)=\sum_{\nu\in \YY}Q_{\lambda/\nu}(A;t,q)P_{\mu/\nu}(B;t,q).
\end{align*}
\end{definition}

\begin{remark}In terms of the scalar product as defined in (\ref{dsp})
\begin{align*}
\Psi_{\lambda,\mu}(A,B;q,t)&=\langle P_{\lambda}(A,Y;q,t),Q_{\mu}(Y,B;q,t) \rangle_Y, \\
\Phi_{\lambda,\mu}(A,B)&=\langle P_{\mu}(Y,B;t,q),Q_{\lambda}(A,Y;t,q) \rangle_Y,
\end{align*}
where $Y$ is a countable set of variables.
\end{remark}

\begin{lemma}\label{l23}Consider dimer coverings on the rail-yard graph with probability measure conditional on left and right boundary conditions $\lambda^{(l)}$ and $\lambda^{(r+1)}$, respectively, given by (\ref{pm}). Then the corresponding sequences of partitions form a generalized Macdonald process as in Definition \ref{df22} with 
\begin{enumerate}
\item $\mathcal{L}=\{i\in[l..r]:a_i=L\}$ and $\mathcal{R}=\{j\in[l..r]:a_j=R\}$; and
\item For $i\in [l..r]$, 
\begin{enumerate}
\item if $b_i=-$, then $A^{(i)}=\{x_i\},B^{(i+1)}=\{0\}$;
\item if $b_i=+$, then $A^{(i)}=\{0\},B^{(i+1)}=\{x_i\}$;
\end{enumerate}
\item $q=t$;
\end{enumerate}
conditional on fixed $\lambda^{(l)}$ and $\lambda^{(r+1)}$ on the left and right boundaries, respectively.
\end{lemma}
\begin{proof}Note that $q=t$ implies  $\Psi=\Phi$ and 
\begin{align*}
\Psi_{\lambda^{(i)},\lambda^{(i+1)}}(A^{(i)},B^{(i+1)};t,t)=\sum_{\nu\in\YY}s_{\lambda^{(i)}/\nu}(A^{(i)})s_{\lambda^{(i+1)}/\nu}(B^{(i+1)}).
\end{align*}
When $b_i=-$,
\begin{align*}
s_{\lambda^{(i+1)}/\nu}(0)=\begin{cases}1&\mathrm{if}\ \nu=\lambda^{(i+1)}\\0&\mathrm{otherwise}\end{cases};
\end{align*}
and therefore
\begin{align*}
\Psi_{\lambda^{(i)},\lambda^{(i+1)}}(x_i,0;t,t)=s_{\lambda^{(i)}/\lambda^{(i+1)}}(x_i)=\Phi_{\lambda^{(i)},\lambda^{(i+1)}}(x_i,0;t,t).
\end{align*}
Similarly, when $b_i=+$,
\begin{align*}
\Psi_{\lambda^{(i)},\lambda^{(i+1)}}(x_i,0;t,t)=s_{\lambda^{(i+1)}/\lambda^{(i)}}(x_i)=
\Phi_{\lambda^{(i)},\lambda^{(i+1)}}(x_i,0;t,t).
\end{align*}
\end{proof}

\section{Moments of Random Height Functions}\label{sect:mh}

In this section, we compute the moments of height functions of perfect matchings on rail-yard graphs by computing the observables in the generalized Macdonald processes. The main result is Lemma \ref{l31}, which implies the formula for the moments given in Lemma \ref{l36}.  

Let $\lambda\in \YY$ be a partition and $q,t\in(0,1)$ be parameters. Let
\begin{align}
\gamma_k(\lambda;q,t)=(1-t^{-k})\sum_{i=1}^{l(\lambda)}q^{k\lambda_i}t^{k(-i+1)}+t^{-kl(\lambda)}.\label{dgm}
\end{align}

\begin{lemma}\label{l41} For $\lambda\in \YY$ and $q,t\in(0,1)$
\begin{align*}
\gamma_k(\lambda';t,q)=\gamma_k\left(\lambda;
\frac{1}{q},\frac{1}{t}\right).
\end{align*}
\end{lemma}

\begin{proof} Let $f_{\lambda}(q,t):=(1-t)\sum_{i\geq 1}(q^{\lambda_i}-1)t^{i-1}$. Then $f_{\lambda}(q,t)=f_{\lambda'}(t,q)$ (see Example 1 in Sect. VI 5 of \cite{IGM15}). Also, note that $f(q,t)=\gamma_1(\lambda;q,\frac{1}{t})-1$. Finally,
\begin{align*}
\gamma_k(\lambda';t,q)=1+f_{\lambda'}\left(t^k,\frac{1}{q^k}\right)=1+f_{\lambda}\left(\frac{1}{q^k},t^k\right)=\gamma_k\left(\lambda;\frac{1}{q},\frac{1}{t}\right).
\end{align*}
Then the lemma follows.
\end{proof}

Let
\begin{align}
H(W,X;q,t)=\prod_{i=1}^k\prod_{x_j\in X}\frac{w_i-\frac{q x_j}{t}}{w_i-q x_j}\label{dh}
\end{align}
and
\begin{align}
\Pi_{L,L}(X,Y)&=\Pi(X,Y;q,t)\notag\\
\Pi_{R,R}(X,Y)&=\Pi(X,Y;t,q)\label{PiLR}\\
\Pi_{L,R}(X,Y)&=\Pi_{R,L}(X,Y)=\Pi'(X,Y),\notag
\end{align}
where $\Pi$ and $\Pi'$ are defined by (\ref{defPPprime}).

\begin{lemma}\label{l31} Let $\mathrm{Pr}$ be the probability measure on pure dimer coverings of the rail-yard graph $RYG(l,r,\underline{a},\underline{b})$ as defined by (\ref{ppd})and (\ref{pbc}); and let
\begin{align*}
\Lambda=\{\lambda^{(i)}\}_{i\in[l+1..r]}
\end{align*}
be the corresponding sequence of partitions. Let $l_i$ be non-negative integers for $i\in[l+1..r]$. Then 
\begin{align*}
&\EE_{\mathrm{Pr}}\left[\prod_{i\in[l+1..r]} \gamma_{l_i}(\lambda^{(i)};t,t)\right]=\oint\ldots\oint 
\prod_{i=[l+1..r]}D(W^{(i)};\omega(t,t,a_i))\\
&\times \prod_{i\leq j;i,j\in[l+1..r]}\left(H(W^{(i)},(-1)^{\delta_{a_i,a_j}-1}A^{(j)};\omega(t,t,a_i))\right)^{(-1)^{\delta_{a_i,a_j}-1}}
\\
&\times \prod_{i<j;i,j\in[l..r]}\frac{\Pi_{a_i,a_j}(B^{(i+1)},W^{(j)})}{\Pi_{a_i,a_j}(B^{(i+1)},\xi(t,t,a_j)W^{(j)})} \times \prod_{i<j;i,j\in[l+1..r]}T_{a_i,a_j}(W^{(i)},W^{(j)}),
\end{align*}
where
\begin{align}
\omega(q,t,a_{i})=\begin{cases}
(q,t)&\mathrm{if}\ a_{i}=L\\
\left(\frac{1}{t},\frac{1}{q}\right)&\mathrm{if}\ a_{i}=R
\end{cases}, \quad \xi(q,t,a_j)=\begin{cases}q^{-1}&\mathrm{If}\ a_j=L\\ t&\mathrm{If}\ a_j=R \end{cases},\label{dox}
\end{align}
$D(W;q,t)$, $H(W,X;q,t)$, and $\Pi_{c,d}(X,Y)$ are given by (\ref{ddf}), (\ref{dh}), and (\ref{PiLR}), and
\begin{align*}
T_{c,d}(Z,W):=\begin{cases}
\prod_{z_i\in Z}\prod_{w_j\in W}\frac{(1-w_jz_i^{-1})^2}{(1-t^{-1}w_j z_i^{-1})
(1-tw_j z_i^{-1})
}&\mathrm{if}\ c=d\\
\prod_{z_i\in Z}\prod_{w_j\in W}\frac{(1+tw_jz_i^{-1})^2}{(1+t^2w_j z_i^{-1})
(1+w_j z_i^{-1})}
&\mathrm{if}\ c=L\ \mathrm{and}\ d=R\\
\prod_{z_i\in Z}\prod_{w_j\in W}\frac{(1+t^{-1}w_jz_i^{-1})^2}{(1+t^{-2}w_j z_i^{-1})
(1+w_j z_i^{-1})
}&\mathrm{if}\ c=R\ \mathrm{and}\ d=L
\end{cases}.
\end{align*}
Furthermore, $|W^{(i)}|=l_i$ and the integral contours are given by $\{\mathcal{C}_{i,j}\}_{i\in [l+1..r],s\in[l_i]}$ such that
\begin{enumerate}
    \item $\mathcal{C}_{i,s}$ is the integral contour for the variable $w^{(i)}_s\in W^{(i)}$;
    \item $\mathcal{C}_{i,s}$ encloses 0 and every singular point of 
    \begin{align*}
\prod_{j\in[i..r]}\left(H(W^{(i)},(-1)^{\delta_{a_i,a_j}-1}A^{(j)};\omega(t,t,a_i)\right)^{(-1)^{\delta_{a_i,a_j}-1}},
\end{align*}
but no other singular points of the integrand;
\item the contour $\mathcal{C}_{i,j}$ is contained in the domain bounded by $t\mathcal{C}_{i',j'}$ whenever $(i,j)<(i',j')$ in lexicographical ordering.
\end{enumerate}
\end{lemma}

\begin{proof}By Lemma \ref{l41}, we obtain
\begin{small}
\begin{align*}
&\left.\EE_{\mathrm{Pr}}\left[\prod_{i\in[l+1..r]} \gamma_{l_i}(\lambda^{(i)};q,t)\right]\right|_{q=t}\\
&=\left.\EE_{\mathrm{Pr}}\left[\prod_{i\in[l+1..r]\cap \mathcal{L}} \gamma_{l_i}(\lambda^{(i)};q,t)\right] \left[\prod_{i\in[l+1..r]\cap \mathcal{R}} \gamma_{l_i}\left(\left[\lambda^{(i)}\right]';\frac{1}{t},\frac{1}{q}\right)\right]
\right|_{q=t}.
\end{align*}
\end{small}
Recall that the Macdonald polynomials satisfy (See Page 324 of \cite{IGM15})
\begin{align}
P_{\lambda}(X;q,t)=P_{\lambda}\left(X;\frac{1}{q},\frac{1}{t}\right);\quad
Q_{\lambda}(X;q,t)=\left(\frac{t}{q}\right)^{|\lambda|}Q_{\lambda}\left(X;\frac{1}{q},\frac{1}{t}\right).\label{pqr}
\end{align}
We obtain
\begin{small}
\begin{align*}
&\left.\EE_{\mathrm{Pr}}\left[\prod_{i=[l+1..r]\cap \mathcal{L}} \gamma_{l_i}(\lambda^{(i)};q,t)\right] \left[\prod_{i=[l+1..r]\cap \mathcal{R}} \gamma_{l_i}\left(\left[\lambda^{(i)}\right]';\frac{1}{t},\frac{1}{q}\right)\right]
\right|_{q=t}\\
&=\sum_{\lambda^{(l)},\ldots,\lambda^{(r+1)}\in\YY}
\left[\prod_{i=[l+1..r]\cap \mathcal{L}} \gamma_{l_i}(\lambda^{(i)};q,t)\right] \left[\prod_{i=[l+1..r]\cap \mathcal{R}} \gamma_{l_i}\left(\left[\lambda^{(i)}\right]';\frac{1}{t},\frac{1}{q}\right)\right]\\
&\times\left.\mathrm{Pr}(\lambda^{(l)},\ldots,\lambda^{(r+1)}|\lambda^{(l)}=\lambda^{(r+1)}=\emptyset)\right|_{q=t}\\
&=\frac{1}{\mathcal{Z}}\sum_{\lambda^{(l+1)},\ldots,\lambda^{(r)}\in\YY}
\left[\prod_{i=[l+1..r]\cap \mathcal{L}} \gamma_{l_i}(\lambda^{(i)};q,t)\right] \left[\prod_{i=[l+1..r]\cap \mathcal{R}} \gamma_{l_i}\left(\left[\lambda^{(i)}\right]';\frac{1}{t},\frac{1}{q}\right)\right]\\
&\times\left[\prod_{i\in\mathcal{L}}\langle P_{\lambda^{(i)}}(A^{(i)},Y^{(i)};q,t),Q_{\lambda^{(i+1)}}(Y^{(i)},B^{(i+1)};q,t) \rangle_{Y^{(i)}}\right]\\
&\left.\left[\prod_{i\in\mathcal{R}}\langle Q_{[\lambda^{(i)}]'}\left(A^{(i)},Y^{(i)};t,q\right),P_{[\lambda^{(i+1)}]'}\left(Y^{(i)},B^{(i+1)};t,q\right) \rangle_{Y^{(i)}}\right]\right|_{q=t,\lambda^{(l)}=\lambda^{(r+1)}=\emptyset}
\end{align*}
\end{small}
where for each $i$, $Y^{(i)}$ is a countable collection of variables; and
\begin{small}
\begin{align*}
\mathcal{Z}&=\sum_{\lambda^{(l+1)},\ldots,\lambda^{(r)}\in\YY}\left[\prod_{i\in\mathcal{L}}\langle P_{\lambda^{(i)}}(A^{(i)},Y^{(i)};q,t),Q_{\lambda^{(i+1)}}(Y^{(i)},B^{(i+1)};q,t) \rangle_{Y^{(i)}}\right]\\
&\left.\left[\prod_{i\in\mathcal{R}}\langle Q_{[\lambda^{(i)}]'}(A^{(i)},Y^{(i)};t,q),P_{[\lambda^{(i+1)}]'}(Y^{(i)},B^{(i+1)};t,q) \rangle_{Y^{(i)}}\right]\right|_{q=t,\lambda^{(l)}=\lambda^{(r+1)}=\emptyset}.
\end{align*}
\end{small}
For $i\in[l+1..r]$, let 
\begin{small}
\begin{align*}
\mathbf{E}_i=\begin{cases}{\displaystyle\sum_{\lambda^{(i)}\in \YY}}\gamma_{l_i}(\lambda^{(i)};q,t)P_{\lambda^{(i)}}(A^{(i)},Y^{(i)};q,t)Q_{\lambda^{(i)}}(Y^{(i-1)},B^{(i)};q,t)&\substack{(i-1,i)\\ \in \mathcal{L}\times\mathcal{L}}\\ 
\displaystyle\sum_{\lambda^{(i)}\in\YY}\gamma_{l_i}([\lambda^{(i)}]';\frac{1}{t},\frac{1}{q})
Q_{\lambda^{(i)}}(Y^{(i-1)},B^{(i)};q,t)Q_{[\lambda^{(i)}]'}(A^{(i)},Y^{(i)};t,q)&\substack{(i-1,i)\\ \in \mathcal{L}\times\mathcal{R}}\\
\displaystyle\sum_{\lambda^{(i)}\in\YY}\gamma_{l_i}(\lambda^{(i)};q,t)P_{[\lambda^{(i)}]'}(Y^{(i-1)},B^{(i)};t,q)P_{\lambda^{(i)}}(A^{(i)},Y^{(i)};q,t)&\substack{(i-1,i)\\ \in \mathcal{R}\times\mathcal{L}}\\
\displaystyle\sum_{\lambda^{(i)}\in\YY}\gamma_{l_i}([\lambda^{(i)}]';\frac{1}{t},\frac{1}{q})P_{[\lambda^{(i)}]'}(Y^{(i-1)},B^{(i)};t,q)Q_{[\lambda^{(i)}]'}(A^{(i)},Y^{(i)};t,q)&\substack{(i-1,i)\\ \in \mathcal{R}\times\mathcal{L}}
\end{cases}
\end{align*}
\end{small}
and
\begin{align*}
\mathbf{E}_l&=\begin{cases}P_{\lambda^{(l)}}(A^{(l)},Y^{(l)};q,t) &l\in \mathcal{L}\\
Q_{[\lambda^{(l)}]'}(A^{(l)},Y^{(l)};t,q)& l\in \mathcal{R}
\end{cases}\\
\mathbf{E}_{r+1}&=\begin{cases}
P_{[\lambda^{(r+1)}]'}(Y^{(r)},B^{(r+1)};t,q)&\ r\in \mathcal{R}\\
Q_{\lambda^{(r+1)}}(Y^{(r)},B^{(r+1)};q,t) & r\in \mathcal{L}
\end{cases}.
\end{align*}
When $\lambda^{(l)}=\lambda^{(r+1)}=\emptyset$ and $q=t$, we have $\mathbf{E}_l=\mathbf{E}_{r+1}=1.$
Then
\begin{align*}
&\left.\EE_{\mathrm{Pr}}
\left[\prod_{i=[l+1..r]\cap \mathcal{L}} \gamma_{l_i}(\lambda^{(i)};q,t)\right] \left[\prod_{i=[l+1..r]\cap \mathcal{R}} \gamma_{l_i}\left(\left[\lambda^{(i)}\right]';t,q\right)\right]
\right|_{q=t}\\
&=\left.\frac{1}{\mathcal{Z}}\langle \mathbf{E}_{l}\langle \mathbf{E}_{l+1}\ldots\langle \mathbf{E}_{r},\mathbf{E}_{r+1}\rangle_{Y^{((r))}}\ldots \rangle_{Y^{(l+1)}} \rangle_{Y^{(l)}}\right|_{q=t,\ \lambda^{(l)}=\lambda^{(r+1)}=\emptyset}.
\end{align*}
Observe that for $i\in[l+1..r]$,
\begin{align*}
&\mathbf{E}_i=
\begin{cases}
D_{-l_i,(A^{(i)},Y^{(i)});q,t}\Pi_{a_{i-1},a_i}((A^{(i)},Y^{(i)}),(Y^{(i-1)},B^{(i)}))&\mathrm{If}\ a_i=L\\
D_{-l_i,(A^{(i)},Y^{(i)});\frac{1}{t},\frac{1}{q}}\Pi_{a_{i-1},a_i}((A^{(i)},Y^{(i)}),(Y^{(i-1)},B^{(i)}))&\mathrm{If}\ a_i=R
\end{cases},
\end{align*}
where $D_{-l_i,(A^{(i)},Y^{(i-1)})}$ is the operator defined as in (\ref{ngt}).

By Proposition \ref{pa2}, we obtain for $i\in[l+1..r]$ and $q=t$,
\begin{itemize}
\item If $a_i=L$,
\begin{align*}
\mathbf{E}_i&=\Pi_{a_{i-1},a_i}((A^{(i)},Y^{(i)}),(Y^{(i-1)},B^{(i)}))\\
&\oint\cdots\oint D(W^{(i)};q,t)H(W^{(i)},(A^{(i)},Y^{(i)});q,t)\frac{\Pi_{a_{i-1},a_i}((Y^{(i-1)},B^{(i)}),W^{(i)})}{\Pi_{a_{i-1},a_i}((Y^{(i-1)},B^{(i)}),q^{-1}W^{(i)})};
\end{align*}
\item If $a_i=R$,
\begin{align*}
\mathbf{E}_i&=\Pi_{a_{i-1},a_i}((A^{(i)},Y^{(i)}),(Y^{(i-1)},B^{(i)}))\\
&\oint\cdots\oint D\left(W^{(i)};\frac{1}{t},\frac{1}{q}\right)H\left(W^{(i)},(A^{(i)},Y^{(i)});\frac{1}{t},\frac{1}{q}\right)\frac{\Pi_{a_{i-1},a_i}((Y^{(i-1)},B^{(i)}),W^{(i)})}{\Pi_{a_{i-1},a_i}((Y^{(i-1)},B^{(i)}),tW^{(i)})};
\end{align*}
\end{itemize}
where each integral contour encloses $0$ and all poles of $H(W^{(i)},(A^{(i)},Y^{(i)});\omega(q,t;a_i))$; moreover, if 
$W^{(i)}=(w^{(i)}_1,w^{(i)}_2,\ldots, w^{(i)}_{l_i})$,
then along the integral contours $|w_j^{(i)}|\leq |\min\{q,t\}w_{j+1}^{(i)}|$ for each $i\in[l_i-1]$.

Since the integrand in each $\mathbf{E}_i$ is $\Lambda_{Y^{(i)}}$-projective, by Lemma \ref{A12} we can interchange the order of the residue and Macdonald scalar product and obtain
\begin{small}
\begin{align*}
&\left.\EE_{\mathrm{Pr}}\left[\prod_{i=[l+1..r]} \gamma_{l_i}(\lambda^{(i)};q,t)\right] \right|_{q=t}
\\&=\frac{1}{\mathcal{Z}}\left(\prod_{i=l+1}^r \Pi_{a_{i-1},a_i}(A^{(i)},B^{(i)})\right)\cdot\oint\langle F_{l}\langle F_{l+1}\ldots\langle F_{r},F_{r+1} \rangle_{Y^{(r)}} \rangle_{Y^{(l+1)}} \rangle_{Y^{(l)}}\\
&\times
\left(\prod_{i=[l+1..r]\cap \mathcal{L}}D(W^{(i)};q,t)H(W^{(i)},A^{(i)};q,t)\frac{\Pi_{a_{i-1},a_i}(B^{(i)},W^{(i)})}{\Pi_{a_{i-1},a_i}(B^{(i)},q^{-1}W^{(i)})}\right)\\
&\left.\times
\left(\prod_{i=[l+1..r]\cap \mathcal{R}}D\left(W^{(i)};
\frac{1}{t},\frac{1}{q}\right)H\left(W^{(i)},A^{(i)};\frac{1}{t},\frac{1}{q}\right)\frac{\Pi_{a_{i-1},a_i}(B^{(i)},W^{(i)})}{\Pi_{a_{i-1},a_i}(B^{(i)},tW^{(i)})}\right)\right|_{q=t}.
\end{align*}
\end{small}
Moreover for $i\in[l+1..r]$,
\begin{itemize}
\item If $a_i=L$,
\begin{align*}
F_i&=\Pi_{a_{i-1},a_i}(A^{(i)},Y^{(i-1)})\cdot\Pi_{a_{i-1},a_i}(Y^{(i-1)},Y^{(i)})\cdot
\Pi_{a_{i-1},a_i}(Y^{(i)},B^{(i)})\\
&\times H(W^{(i)},Y^{(i)};q,t)\frac{\Pi_{a_{i-1},a_i}(Y^{(i-1)},W^{(i)})}{\Pi_{a_{i-1},a_i}(Y^{(i-1)},q^{-1}W^{(i)})};
\end{align*}
\item If $a_i=R$,
\begin{align*}
F_i&=\Pi_{a_{i-1},a_i}(A^{(i)},Y^{(i-1)})\cdot\Pi_{a_{i-1},a_i}(Y^{(i-1)},Y^{(i)})\cdot
\Pi_{a_{i-1},a_i}(Y^{(i)},B^{(i)})\\
&\times H\left(W^{(i)},Y^{(i)};\frac{1}{t},\frac{1}{q}\right)\frac{\Pi_{a_{i-1},a_i}(Y^{(i-1)},W^{(i)})}{\Pi_{a_{i-1},a_i}(Y^{(i-1)},tW^{(i)})};
\end{align*}
\end{itemize}
and $F_l=F_{r+1}=1.$

By Lemmas \ref{la5} and \ref{la6}, we obtain
\begin{align*}
\mathcal{F}_r:=
\langle F_r,F_{r+1} \rangle_{Y^{(r)}}=
\begin{cases}
\frac{\Pi_{a_{r-1},a_r}((A^{(r)},W^{(r)}),Y^{(r-1)})}{\Pi_{a_{r-1},a_r}(q^{-1}W^{(r)},Y^{(r-1)})}&\mathrm{If}\ a_r=L\\
\frac{\Pi_{a_{r-1},a_r}((A^{(r)},W^{(r)}),Y^{(r-1)})}{\Pi_{a_{r-1},a_r}(tW^{(r)},Y^{(r-1)})}&\mathrm{If}\ a_r=R
\end{cases}.
\end{align*}
Then the lemma follows by inductively computing the scalar product 
\begin{align*}
\langle F_{l}\langle F_{l+1}\ldots\langle F_{r},F_{r+1} \rangle_{Y^{(r)}} \rangle_{Y^{(l+1)}} \rangle_{Y^{(l)}};
\end{align*}
and applying Lemmas \ref{l32} and \ref{l33}. Note that factors $\Pi_{a_i,a_j}(B^{(i+1)},A^{(j)})$ cancel out with $1/\mathcal{Z}$ since the partition function for the Macdonald process with empty partitions on both ends is
\begin{align*}
\prod_{i<j;i,j\in[l,..,r]}\Pi_{a_i,a_j}(B^{(i+1)},A^{(j)}).
\end{align*}
\end{proof}

\begin{remark}Using similar arguments, we can also obtain the following formula.
\begin{align*}
&\EE_{\mathrm{Pr}}\left[\prod_{i\in[l+1..r]\cap \mathcal{L}} \gamma_{l_i}(\lambda^{(i)};t,t) \prod_{i\in[l+1..r]\cap \mathcal{R}} \gamma_{l_i}\left([\lambda^{(i)}]';t,t\right)\right]\\
&=\oint\ldots\oint \prod_{i \in [l+1..r]}D(W^{(i)};t,t)\\
&\times \prod_{i\leq j;i,j\in[l+1..r]}\left(H(W^{(i)},(-1)^{\delta_{a_i,a_j}-1}A^{(j)};t,t)\right)^{(-1)^{\delta_{a_i,a_j}-1}}\\
&\times \prod_{i<j;i,j\in[l..r]}\frac{\Pi_{a_i,a_j}(B^{(i+1)},W^{(j)})}{\Pi_{a_i,a_j}(B^{(i+1)},t^{-1}W^{(j)})}\times \prod_{i< j;i,j\in[l+1..r]}S_{a_i,a_j}(W^{(i)},W^{(j)}),
\end{align*}
where $D(W;q,t)$, $H(W,X;q,t)$, and $\Pi_{c,d}(X,Y)$ are given by (\ref{ddf}), (\ref{dh}), and (\ref{PiLR}) and
\begin{align*}
S_{c,d}(Z,W):=\left[\prod_{z_i\in Z}\prod_{w_j\in W}\frac{(1+(-1)^{\delta_{c,d}}t^{-1}w_j z_i^{-1})
(1+(-1)^{\delta_{c,d}}tw_j z_i^{-1})
}{(1+(-1)^{\delta_{c,d}}w_jz_i^{-1})^2}\right]^{(-1)^{\delta_{c,d}}}.
\end{align*}
Furthermore, $|W^{(i)}|=l_i$ and the integral contours are given by $\{\mathcal{C}_{i,j}\}_{i\in [l+1..r],s\in[l_i]}$ such that
\begin{enumerate}
    \item $\mathcal{C}_{i,s}$ is the integral contour for the variable $w^{(i)}_s\in W^{(i)}$;
    \item $\mathcal{C}_{i,s}$ encloses 0 and every singular point of 
    \begin{align*}
\prod_{j\in[i..r]}\left(H(W^{(i)},(-1)^{\delta_{a_i,a_j}-1}A^{(j)};t,t)\right)^{(-1)^{\delta_{a_i,a_j}-1}},
\end{align*}
but no other singular points of the integrand;
    \item the contour $\mathcal{C}_{i,j}$ is contained in the domain bounded by $t\mathcal{C}_{i',j'}$ whenever $(i,j)<(i',j')$ in lexicographical ordering.
\end{enumerate}
\end{remark}

\begin{lemma}\label{l32}Let $c_1,c_2,c_3\in\{L,R\}$. Let $A,B,Y$ be 3 collections of countably many variables. 
Then we have
\begin{align*}
\langle \Pi_{c_1,c_2}(A,Y),\Pi_{c_2,c_3}(Y,B) \rangle_Y
=\Pi_{c_1,c_3}(A,B).
\end{align*}
where if $c_2=L$ (resp.\ $c_2=R$), the scalar product $\langle \cdot,\cdot\rangle$ is with respect to $(q,t)$ (resp.\ $(t,q)$).
\end{lemma}

\begin{lemma}\label{l33}For $i\in[l+2..r]$, and $j\in[i..r]$,
let
\begin{align*}
\mathcal{G}:&=\left\langle \Pi_{a_{i-2},a_{i-1}}\left(\left(Y^{(i-2)},B^{(i-1)}\right),Y^{(i-1)}\right)H\left(W^{(i-1)},Y^{(i-1)};\omega(q,t;a_{i-1})\right), \right.\\
&\left.\frac{\Pi_{a_{i-1},a_j}((A^{(j)},W^{(j)}),Y^{(i-1)})}{\Pi_{a_{i-1},a_j}(\xi(q,t,a_j)W^{(j)},Y^{(i-1)})}\right\rangle_{Y^{(i-1)}},
\end{align*}
where in $\mathcal{G}$, the scalar product $\langle\cdot,\cdot \rangle$ is with respect to $(q,t)$ if $a_{i-1}=L$ and with respect to $(t,q)$ if $a_{i-1}=R$; and $\omega$, $\xi$ are defined as in (\ref{dox}).
Assume $q=t$, then $\mathcal{G}$ is equal to
\begin{align*}
\mathcal{G}&=\frac{\Pi_{a_{i-2},a_j}((Y^{(i-2)},B^{(i-1)}),(A^{(j)},W^{(j)}))}
{\Pi_{a_{i-2},a_j}(Y^{(i-2)},B^{(i-1)}),\xi(t,t,a_j)W^{(j)})
}\\
&\times \left[\frac{H(W^{(i-1)},(A^{(j)},W^{(j)}),\omega(t,t;a_{i-1}))}{
H(W^{(i-1)},(-1)^{\delta_{a_{i-1},a_j}-1}\xi(t,t,a_j)W^{(j)},\omega(t,t;a_{i-1}))
}\right]^{(-1)^{\delta_{a_{i-1},a_j}-1}}.
\end{align*}
\end{lemma}

\begin{proof}The lemma follows from Lemma \ref{la5} with
\begin{align*}
d_n=\begin{cases}
(-1)^{(n+1)(\delta_{a_{i-2},a_{i-1}}-1)}p_n(Y^{(i-2)},B^{(i-1)})+(1-t^{-n})p_n(\frac{q}{W^{(i-1)}}),&a_{i-1}=L\\
(-1)^{(n+1)(\delta_{a_{i-2},a_{i-1}}-1)}p_n(Y^{(i-2)},B^{(i-1)})+(1-q^{n})p_n(\frac{1}{tW^{(i-1)}}),&a_{i-1}=R
\end{cases}
\end{align*}
\begin{align*}
u_n&=\begin{cases}
(-1)^{(n+1)}\left[p_n(A^{(j)},W^{(j)})-p_n\left(W^{(j)}\xi(a_j)\right)\right]&\mathrm{If}\ a_{i-1}\neq a_j\\
\frac{1-t^n}{1-q^n}\left[p_n(A^{(j)},W^{(j)})-p_n\left(W^{(j)}\xi(a_j)\right)\right]&\mathrm{If}\ a_{i-1}=a_j=L\\
\frac{1-q^{-n}}{1-t^{-n}}\left[p_n(A^{(j)},W^{(j)})-p_n\left(W^{(j)}\xi(a_j)\right)\right]&\mathrm{If}\ a_{i-1}=a_j=R
\end{cases}.
\end{align*}
\end{proof}

\begin{lemma}\label{l35}Let $\mathrm{Pr}$ be the probability measure on pure dimer coverings of the rail-yard graph $RYG(l,r,\underline{a},\underline{b})$ as defined by (\ref{ppd}) and (\ref{pbc}); and let
\begin{align*}
\Lambda=\{\lambda^{(i)}\}_{i\in[l+1..r]}
\end{align*}
be the corresponding sequence of partitions. Let $l_i$ be non-negative integers for $i\in[l+1..r]$. Then
\begin{align*}
&\mathbb{E}_{\mathrm{Pr}}\left[\prod_{i\in[l+1..r]} \gamma_{l_i}(\lambda^{(i)};t,t)\right] 
=\oint\ldots\oint  \prod_{i=[l+1..r]}D(W^{(i)};\omega(t,t;a_i)) \\
&\times \prod_{\substack{i,j\in[l+1..r] \\i \leq j, b_j=-}}\left(H(W^{(i)},(-1)^{\delta_{a_i,a_j}-1}\{x_j\};\omega(t,t;a_i))\right)^{(-1)^{\delta_{a_i,a_j}-1}}\\
&\times \prod_{\substack{i,j\in[l..r]\\i<j,b_i=+}}\frac{\Pi_{a_i,a_j}(\{x_i\},W^{(j)})}{\Pi_{a_i,a_j}(\{x_i\},\xi(t,t;a_j)W^{(j)})} \times \prod_{\substack{i,j\in[l+1..r]\\i<j}}T_{a_i,a_j}(W^{(i)},W^{(j)}),
\end{align*}
where $|W^{(i)}|=l_i$, and the integral contours are given by $\{\mathcal{C}_{i,j}\}_{i\in [l+1..r],s\in[l_i]}$ satisfying the condition as described in Lemma \ref{l31}.
\end{lemma}

\begin{proof} Follows from Lemmas \ref{l31} and \ref{l23}(2) and $H(X,\{0\};q,t)=1$ and $\Pi_{a,b}(X,\{0\})=\Pi_{a,b}(\{0\},X)=1$.
\end{proof}

To compute the moments of $\gamma$s we need to generalize Lemma \ref{l35} and find expectation of powers of $\gamma$s. We use the following auxiliary lemma, which is straightforward from the properties of Macdonald functions.  \begin{lemma}\label{MacPrRest}Let $\bA$, $\bB$ be defined as in Definition \ref{df22}. For each $i\in[l..r]$, let $\Lambda_{A^{(i)}}$, (resp. $\Lambda_{B^{(i+1)}}$) be the algebra of symmetric functions on $A^{(i)}$ (resp. $B^{(i+1)}$) over $\CC$. Define a map 
\begin{align*}
\phi_0^{(i)}:\Lambda_{A^{(i)}}\otimes \Lambda_{B^{(i+1)}}\rightarrow \CC
\end{align*}
by
\begin{align*}
\phi_0^{(i)}(f\otimes g)=f(\mathbf{0})g(\mathbf{0}).
\end{align*}
Define a formal measure
\begin{align*}
M^{(i)}:=\phi_0^{(i)}\left(\MP_{\bA,\bB,\mathcal{P},q,t}\right).
\end{align*}
Then for any sequence of partitions $(\lambda^{(l)},\lambda^{(l+1)},\ldots,\lambda^{(r+1)})\in \YY^{r-l+2}$ in the support of $M^{(i)}$, we have $\lambda^{(i)}=\lambda^{(i+1)}$. Moreover, for all $i\in [l..r]$ let
\begin{align*}
\widehat{\bA}^{(i)}=\bA\setminus\{A^{(i)}\},\qquad
\widehat{\bB}^{(i+1)}=\bB\setminus\{B^{(i+1)}\},
\end{align*}
and
\begin{align*}
\widehat{\mathcal{L}}^{(i)}=\mathcal{L}\setminus\{i\},\qquad
\widehat{\mathcal{R}}^{(i)}=\mathcal{R}\setminus\{i\},
\end{align*}
so that $\{\widehat{\mathcal{L}}^{(i)},\widehat{\mathcal{R}}^{(i)}\}$ form a partition $\mathcal{P}'$ of $[l..r]\setminus\{i\}$. Then the restriction of $M^{(i)}$ to 
\begin{align*}
(\lambda^{(l)},\ldots,\lambda^{(i-1)},\lambda^{(i+1)},\ldots,\lambda^{(r+1)})\in \YY^{(r-l+1)}
\end{align*}
is the formal Macdonald process $\MP_{\widehat{A}^{(i)},\widehat{B}^{(i)},\mathcal{P}',q,t}$.
\end{lemma}

\begin{lemma}\label{l36}Let $i_1\leq i_2\leq \ldots\leq i_m\in[l+1..r]$, and let $l_1,\ldots,l_m>0$ be integers. Let $I:=\{i_1,i_2,\ldots,i_m\}$. Then
\begin{align*}
&\mathbb{E}_{\mathrm{Pr}}\left[\prod_{j\in [m]} \gamma_{l_j}(\lambda^{(i_j)};t,t)\right]=\oint\ldots\oint \prod_{i\in I}D(W^{(i)};\omega(t,t;a_i))\\
&\times\prod_{\substack{i \in I, j\in[l+1..r]\\ i \leq j, b_j=-.}}\left(H(W^{(i)},(-1)^{\delta_{a_i,a_j}-1}\{x_j\};\omega(t,t;a_i))\right)^{(-1)^{\delta_{a_i,a_j}-1}}\\
&\times \prod_{\substack{j \in I,i\in[l..r]\\i<j,b_i=+}}\frac{\Pi_{a_i,a_j}(\{x_i\},W^{(j)})}{\Pi_{a_i,a_j}(\{x_i\},\xi(t,t;a_j)W^{(j)})} \times \prod_{1\leq s<j\leq m}T_{a_{i_s},a_{i_j}}\left(W^{(i_s)},W^{(i_j)}\right),
\end{align*}
where $|W^{(i)}|=l_i$, and the integral contours are given by $\{\mathcal{C}_{i,j}\}_{i\in I,s\in[l_i]}$ satisfying the condition as described in Lemma \ref{l31}.
\end{lemma}

\begin{proof} The proof is the same as in \cite{bcgs13}, where Corollary 3.11 is derived from Theorem 3.10. 
The idea is to apply Lemma \ref{l31} to an auxiliary Macdonald process $
\left.\MP_{\bC,\bD,\mathcal{P}_*,q,t}\right|_{\mu^{(0)}=\mu^{(r-l+m+1)}=\emptyset}$, where parameters $\bC$ and $\bD$ consist of $r-l+m$ variables and the additional indices correspond to copies of $i_j$ for $j\in [m]$. The expectation formula from Lemma \ref{l31} restricted to the original process, as in Lemma \ref{MacPrRest}, gives the desired formula. For details see the proof of Corollary 3.11 in
\cite{bcgs13}.
\end{proof}

\section{Asymptotics}\label{sect:as}
In this section, we study the asymptotics of the moments of the random height functions and prove its Gaussian fluctuation in the scaling limit. More precisely, we study the limit of the moments of the observables obtained in the previous section. The main results are given in Theorem \ref{t57} and Theorem \ref{t58}. We first specify the conditions under which the limit is taken.


We consider a sequence of rail-yard graphs that depend on $\epsilon$ and study the limit when $\epsilon \to 0$. The conditions can be split in three groups: an assumption on piecewise periodicity of the graphs, which are described by the sequences $\underline{a}$ and $\underline{b}$, assumption on periodicity of weights, which allows for periodic non-uniform weights on diagonal edges in the $q$-volume analog of the uniform model, and the limit regime, which contains further assumptions under which the limit is taken.

\begin{assumption}\label{ap5} Let $\{RYG(l^{(\epsilon)},r^{(\epsilon)},\underline{a}^{(\epsilon)},\underline{b}^{(\epsilon)})\}_{\epsilon>0}$  be a sequence of rail-yard graphs with the weights of diagonal edges incident with $x=2i$ given with $x_i^{(\epsilon)}$.

\begin{enumerate}
\item \label{ap41} \textbf{Piecewise periodicity of the graph.} For a positive integer $n$ and real numbers $V_0<V_1<\ldots<V_m$, we say that a given sequence of rail-yard graphs is $n$-periodic with transition points $V_0, V_1, \dots,V_m$ as $\epsilon\rightarrow 0$ if
\begin{enumerate}
\item For each $\epsilon>0$, there exist integer multiples of $n$
\begin{align*}
l^{(\epsilon)}=v_0^{(\epsilon)}<v_1^{(\epsilon)}<\ldots<v_m^{(\epsilon)}=r^{(\epsilon)},
\end{align*}
such that $\lim_{\epsilon\rightarrow 0}\epsilon v_p^{(\epsilon)}=V_p,\ \forall p\in\{0\}\cup[m].$
\item The sequence $\underline{a}^{(\epsilon)}$ is $n$-periodic on $\left[l^{(\epsilon)},r^{(\epsilon)}\right]$ and does not depend on $\epsilon$. More precisely, there exist $a_1, a_2, \dots, a_n \in \{L,R\}$ such that
\begin{align*}
a_i^{(\epsilon)}=a_{i_{\equiv_n}},
\end{align*}
where $i_{\equiv_n}\in[n]$ is in the same congruence class modulo $n$ as $i$.
\item For each $p\in[m]$, the sequence
$\underline{b}^{(\epsilon)}$ is $n$-periodic on $(v_{p-1}^{(\epsilon)},v_{p}^{(\epsilon)})$ and does not depend on $\epsilon$, but it may depend on $p$. More precisely, there exist $b_{p,1}, b_{p,2}, \dots, b_{p,n} \in \{+,-\}$ such that for $i\in (v_{p-1}^{(\epsilon)},v_{p}^{(\epsilon)})$
\begin{align*}
b_i^{(\epsilon)}=b_{p,i_{\equiv_n}}. 
\end{align*}
\end{enumerate}

\item \label{ap43} \textbf{Periodicity of weights.} The weights $x_i^{(\epsilon)}$ are periodic q-volume weights. Precisely they are given with
\begin{align*}
x_i^{(\epsilon)}=\begin{cases}
e^{-\epsilon(i-i_{\equiv_n})}\tau_k,& b_i^{(\epsilon)}=+\\
e^{\epsilon(i-i_{\equiv_n})}\tau^{-1}_k,& b_i^{(\epsilon)}=-
\end{cases}.
\end{align*}

\item \label{ap52} \textbf{Limit regime.} 
\begin{enumerate}
    \item Assume that $\lim_{\epsilon\rightarrow 0}t=1$ in a such a way that 
\begin{align*}
\lim_{\epsilon\rightarrow 0}-\frac{\log t}{n \epsilon}=\beta,
\end{align*}
where $\beta$ is a positive real number independent of $\epsilon$.
\item Let $s$ be a positive integer and assume that for $d \in [s]$ the sequence $i_d^{(\epsilon)}\in[l^{(\epsilon)}..r^{(\epsilon)}]$ satisfies 
\begin{align*}
\lim_{\epsilon\rightarrow 0}\epsilon i_d^{(\epsilon)}=\chi_d, 
\end{align*}
for $\chi_1\leq \chi_2\leq \ldots\leq \chi_s$ and that $i_d^{(\epsilon)}\,\mathrm{mod}\, n$ does not depend on $\epsilon$, i.e. that there exist $i^*_1, i^*_2, \dots, i^*_s \in [n]$ such that
\begin{align*}
(i_d^{(\epsilon)})_{\equiv_n} =  i^*_d.
\end{align*}
If $s=1$ we drop the index, i.e. we assume that $\lim_{\epsilon\rightarrow 0}\epsilon i^{(\epsilon)}=\chi$ and $(i^{(\epsilon)})_{\equiv_n} =  i^*$.
\end{enumerate}
\end{enumerate}
\end{assumption}

\begin{lemma} \label{lLeftMoments}Under the same conditions as in Lemma \ref{l36} and assuming the index set $I$ is a subset of $\mathcal{L}$
\begin{align*}
\mathbb{E}&_{\mathrm{Pr}}\left[\prod_{j\in [m]} \gamma_{l_j}(\lambda^{(i_j)};t,t)\right]=\oint\ldots\oint \prod_{i\in I}D(W^{(i)};t,t)\notag\\
&\times\prod_{\substack{i\in I,j\in[l+1..r],j\geq i\\b_j=-,a_j=a_i}}G_{1,>}(W^{(i)},x_j,t)
\prod_{\substack{i\in I,j\in[l+1..r],j \geq i\\b_j=-,a_j \neq a_i}}G_{0,>}(W^{(i)},x_j,t)\notag \\
&\times\prod_{\substack{i\in I,j\in[l..r],j<i\\b_j=+,a_j= a_i}}G_{1,<}(W^{(i)},x_j,t)
\prod_{\substack{i\in I,j\in[l..r],j<i\\b_j=+,a_j\neq  a_i}}G_{0,<}(W^{(i)},x_j,t)\\
&\times \prod_{1\leq s<j\leq m}T_{L,L}(W^{(i_s)},W^{(i_j)}),
\end{align*}
where $|W^{(i)}|=l_i$, and the integral contours are given by $\{\mathcal{C}_{i,j}\}_{i\in I,s\in[l_i]}$ satisfying the condition as described in Lemma \ref{l31}, where 
\begin{align}
G_{1,>}(W,x,t)&=\prod_{w_s\in W}\frac{w_s-x_j}{w_s-tx_j},\quad G_{1,<}(W,x,t)=\prod_{w_s\in W}\frac{t-w_sx}{t(1-w_s x)},\label{dg1}\\
G_{0,>}(W,x,t)&=\prod_{w_s\in W}\frac{w_s+tx_j}{w_s+x_j},\quad G_{0,<}(W,x,t)=\prod_{w_s \in W}\frac{t(1+w_sx)}{t+w_sx},\label{dg2}
\end{align}
\begin{align}
T_{L,L}(Z,W):=\prod_{z_i\in Z}\prod_{w_j\in W}\frac{\left(z_i-w_j\right)^2}{\left(z_i-t^{-1}w_j \right)
\left(z_i-tw_j \right)}.\label{tll}
\end{align}
\end{lemma}
We apply this lemma to the sequence of rail-yard graphs depending on $\epsilon$ and consider the limit of it under the assumption stated above. To shorten the notation we use the following abbreviations
\begin{align*}
G^{(\epsilon)}_{1,>i}(W)&=\prod_{\substack{j\in[(l+1)^{(\epsilon)}..r^{(\epsilon)}] ,j\geq i\\b^{(\epsilon)}_j=-,a_j=a_i}}G_{1,>}(W,x_j^{(\epsilon)},t),\\
G^{(\epsilon)}_{1,<i}(W)&=\prod_{\substack{j\in[l^{(\epsilon)}..r^{(\epsilon)}],j<i\\b^{(\epsilon)}_j=+,a_j= a_i}}G_{1,<}(W,x^{(\epsilon)}_j,t),\\
G^{(\epsilon)}_{0,>i}(W)&=\prod_{\substack{j\in[(l+1)^{(\epsilon)}..r^{(\epsilon)}],j \geq i\\b^{(\epsilon)}_j=-,a_j \neq a_i}}G_{0,>}(W,x^{(\epsilon)}_j,t),\\
G^{(\epsilon)}_{0,<i}(W)&=\prod_{\substack{j\in[l^{(\epsilon)}..r^{(\epsilon)}],j<i\\b^{(\epsilon)}_j=+,a_j \neq  a_i}}G_{0,<}(W,x^{(\epsilon)}_j,t).
\end{align*}

The limit of these functions is given by the following lemma.

\begin{lemma}\label{l55}Suppose Assumption \ref{ap5} holds. Then
\begin{align*}
 \lim_{\epsilon\rightarrow 0}G^{(\epsilon)}_{1,>i^{(\epsilon)}}(W)
 &=\left[\prod_{w_s\in W}\mathcal{G}_{1,>\chi}(w_s)\right]^{\beta}, \quad \lim_{\epsilon\rightarrow 0}G^{(\epsilon)}_{1<i^{(\epsilon)}}(W)
= \left[\prod_{w_s\in W}\mathcal{G}_{1,<\chi}(w_s)\right]^{\beta},\\
 \lim_{\epsilon\rightarrow 0}G^{(\epsilon)}_{0,>i^{(\epsilon)}}(W)
 &=\left[\prod_{w_s\in W}\mathcal{G}_{0,>\chi}(w_s)\right]^{\beta},\quad
 \lim_{\epsilon\rightarrow 0}G^{(\epsilon)}_{0,<i^{(\epsilon)}}(W)=\left[\prod_{w_s\in W}\mathcal{G}_{0,<\chi}(w_s)\right]^{\beta},
 \end{align*}
where
\begin{align}
\mathcal{G}_{1,>\chi}(w):&=\prod_{(p\in [m],V_p>\chi)}\prod_{\substack{j\in[n]\\b_{p,j}=-,a_j=a_{i^*}}}
\frac{1-\left[w\tau_j\right]^{-1}e^{V_p}}{1-e^{\max\{V_{p-1},\chi\}}\left[w\tau_j\right]^{-1}}\label{dgs1},\\
\mathcal{G}_{1,<\chi}(w):&=\prod_{(p\in[m],V_{p-1}<\chi)}\prod_{\substack{j\in[n]\\b_{p,j}=+,a_j=a_{i^*}}}
\frac{1-w e^{-V_{p-1}}\tau_j}{1-e^{-\min\{V_p,\chi\}}w \tau_j}\label{dgs2},\\
\mathcal{G}_{0,>\chi}(w):&=\prod_{(p\in[m],V_p>\chi)}\prod_{\substack{j\in[n]\\b_{p,j}=-,a_j\neq a_{i^*}}}
\frac{1+e^{\max\{V_{p-1},\chi\}}\left[w\tau_j\right]^{-1}}{1+\left[w\tau_j\right]^{-1}e^{V_p}}\label{dgs3},\\
\mathcal{G}_{0,<\chi}(w):&=\prod_{(p\in[m],V_{p-1}<\chi)}\prod_{\substack{j\in[n]\\b_{p,j}=+,a_j\neq a_{i^*}}}
\frac{1+e^{-\min\{V_p,\chi\}}w \tau_j}{1+w e^{-V_{p-1}}\tau_j}.\label{dgs4}
\end{align}

Here the logarithmic branches for $\mathcal{G}_{1,>\chi}$, $\mathcal{G}_{1,<\chi}$, $\mathcal{G}_{0,>\chi}$, $\mathcal{G}_{0,<\chi}$ are chosen so that when $z$ approaches the positive real axis, the imaginary part of $\log z$ approaches 0.
\end{lemma}

\begin{proof}
For $p\in [m]$, $j\in[n]$ and $i\in[l^{(\epsilon)}..r^{(\epsilon)}]$, let
\begin{align*}
I^{(\epsilon)}_{j,p,>i,1}&=\left\{u\in\left[v_{p-1}^{(\epsilon)}+1,v_p^{(\epsilon)}\right]\cap \{n\ZZ+j\}\cap[i..r^{(\epsilon)}]:b^{(\epsilon)}_u=-,a_u=a_i\right\},\\
I^{(\epsilon)}_{j,p,<i,1}&=\left\{u\in\left[v_{p-1}^{(\epsilon)}+1,v_p^{(\epsilon)}\right]\cap \{n\ZZ+j\}\cap[l^{(\epsilon)}..i-1]:b^{(\epsilon)}_u=+,a_u=a_i\right\},\\
I^{(\epsilon)}_{j,p,>i,0}&=\left\{u\in\left[v_{p-1}^{(\epsilon)}+1,v_p^{(\epsilon)}\right]\cap \{n\ZZ+j\}\cap[i..r^{(\epsilon)}]:b^{(\epsilon)}_u=-,a_u\neq a_i\right\},\\
I^{(\epsilon)}_{j,p,<i,0}&=\left\{u\in\left[v_{p-1}^{(\epsilon)}+1,v_p^{(\epsilon)}\right]\cap \{n\ZZ+j\}\cap[l^{(\epsilon)}..i-1]:b^{(\epsilon)}_u=+,a_u\neq a_i\right\}.
\end{align*}
Let $N_{j,p,>i,1}^{(\epsilon)}$, $N_{j,p,<i,1}^{(\epsilon)}$, $N_{j,p,>i,0}^{(\epsilon)}$, and $N_{j,p,<i,0}^{(\epsilon)}$ be their cardinalities, respectively, and
\begin{align*}
&q_{j,p,>i,1}=\max I^{(\epsilon)}_{j,p,>i,1},\qquad
q_{j,p,<i,1}=\min I^{(\epsilon)}_{j,p,<i,1},\\
&q_{j,p,>i,0}=\max I^{(\epsilon)}_{j,p,>i,0},\qquad
q_{j,p,<i,0}=\min I^{(\epsilon)}_{j,p,<i,0},
\end{align*}
where we take the convention that the minimum (resp.\ maximum) of an empty set is $\infty$ (resp.\ $-\infty$), and define $x_{-\infty}^{(\epsilon)}=x_{\infty}^{(\epsilon)}:=0$ for all $\epsilon\geq 0.$

Then using the above notation we can rewrite $G^{(\epsilon)}_{1,>i}(W)$ as
\begin{align*}
G^{(\epsilon)}_{1,>i}(W)=\prod_{p=1}^{m}\prod_{j=1}^{n}\prod_{c\in I_{j,p,>i,1}^{(\epsilon)}} \left[G_{1,>}
(W,e^{-n\epsilon(q_{j,p,>i,1}-c)}x_{q_{j,p,>i,1}}^{(\epsilon)},t)\right],
\end{align*}
which can be further rewritten as 
\begin{align*}
G^{(\epsilon)}_{1,>i}(W)=\prod_{p=1}^{m}\prod_{j=1}^{n}\prod_{w_s\in W}
\frac{\left(w_s^{-1}x_{q_{j,p,>i,1}}^{(\epsilon)};e^{-n\epsilon}\right)_{N_{j,p,>i,1}^{(\epsilon)}}}{
\left(tw_s^{-1}x_{q_{j,p,>i,1}}^{(\epsilon)};e^{-n\epsilon}\right)_{N_{j,p,>i,1}^{(\epsilon)}}
},
\end{align*}
where $(a;q)_{N}=\prod_{i=0}^{N-1}(1-aq^{i})$. Similarly, we have
\begin{align*}
G^{(\epsilon)}_{1,<i}(W)
&=\prod_{p=1}^{m}\prod_{j=1}^{n}\prod_{w_s \in W}
\frac{\left( t^{-1}w_sx_{q_{j,p,<i,1}}^{(\epsilon)};e^{-n\epsilon}\right)_{N_{j,p,<i,1}^{(\epsilon)}}}{
\left(w_sx_{q_{j,p,<i,1}}^{(\epsilon)};e^{-n\epsilon}\right)_{N_{j,p,<i,1}^{(\epsilon)}}
},
\end{align*}
\begin{align*}
G^{(\epsilon)}_{0,>i}(W)=\prod_{p=1}^{m}\prod_{j=1}^{n}\prod_{w_s\in W}
\frac{\left( -t w_s^{-1}x_{q_{j,p,>i,0}}^{(\epsilon)};e^{-n\epsilon}\right)_{N_{j,p,>i,0}^{(\epsilon)}}}{
\left(-w_s^{-1}x_{q_{j,p,>i,0}}^{(\epsilon)};e^{-n\epsilon}\right)_{N_{j,p,>i,0}^{(\epsilon)}}
},
\end{align*}
\begin{align*}
G^{(\epsilon)}_{0,<i}=\prod_{p=1}^{m}\prod_{j=1}^{n}\prod_{w_s\in W}
\frac{\left( -w_sx_{q_{j,p,<i,0}}^{(\epsilon)},e^{-n\epsilon}\right)_{N_{j,p,<i,0}^{(\epsilon)}}}{
\left(-t^{-1}w_sx_{q_{j,p,<i,0}}^{(\epsilon)},e^{-n\epsilon}\right)_{N_{j,p,<i,0}^{(\epsilon)}}
}.
\end{align*}

Assume now that $i=i^{(\epsilon)}$ vary with $\epsilon$ and satisfies the assumptions. For convenience we continue writing $i$ instead of  $i^{(\epsilon)}$. By Lemma \ref{al51}, we obtain as $\epsilon\to 0$.
\begin{align*}
G^{(\epsilon)}_{1,>i}(W)\sim\prod_{p=1}^{m}\prod_{j=1}^{n}\prod_{w_s\in W}
\left(\frac{1-w_s^{-1}x_{q_{j,p,>i,1}}^{(\epsilon)}}{1-e^{-n\epsilon N_{j,p,>i,1}^{(\epsilon)}}w_s^{-1}x_{q_{j,p,>i,1}}^{(\epsilon)}}\right)^{-\frac{\log t}{n\epsilon}},
\end{align*}

\begin{align*}
G^{(\epsilon)}_{1<i}(W)\sim \prod_{p=1}^{m}\prod_{j=1}^{n}\prod_{w_s\in W}
\left(\frac{1-t^{-1}w_sx_{q_{j,p,<i,1}}^{(\epsilon)}}{1-e^{-n\epsilon N_{j,p,<i,1}^{(\epsilon)}}t^{-1}w_sx_{q_{j,p,<i,1}}^{(\epsilon)}}\right)^{-\frac{\log t}{n\epsilon}},
\end{align*}

\begin{align*}
G^{(\epsilon)}_{0,>i}(W)\sim\prod_{p=1}^{m}\prod_{j=1}^{n}\prod_{w_s\in W}
\left(\frac{1+w_s^{-1}x_{q_{j,p,>i,0}}^{(\epsilon)}}{1+e^{-n\epsilon N_{j,p,>i,0}^{(\epsilon)}}w_s^{-1}x_{q_{j,p,>i,0}}^{(\epsilon)}}\right)^{\frac{\log t}{n\epsilon}},
\end{align*}

\begin{align*}
G^{(\epsilon)}_{0,<i}(W)\sim\prod_{p=1}^{m}\prod_{j=1}^{n}\prod_{w_s\in W}
\left(\frac{1+t^{-1}w_sx_{q_{j,p,<i,0}}^{(\epsilon)}}{1+e^{-n\epsilon N_{j,p,<i,0}^{(\epsilon)}}w_sx_{q_{j,p,<i,0}}^{(\epsilon)}}\right)^{\frac{\log t}{n\epsilon}},
\end{align*}
where when $W\in \CC^k$, we choose the branch such that when a complex number approaches the positive real line, its argument approaches 0.

Note that for nonempty $I^{(\epsilon)}_{j,p,>i,1}$, $I^{(\epsilon)}_{j,p,<i,1}$, $I^{(\epsilon)}_{j,p,>i,0}$, and $I^{(\epsilon)}_{j,p,<i,0}$ we have
\begin{align*}
\lim_{\epsilon\rightarrow 0} n\epsilon N_{j,p,>i,1}^{(\epsilon)}
&
= V_p-\max\{V_{p-1},\chi\},\\
\lim_{\epsilon\rightarrow 0} n\epsilon N_{j,p,<i,1}^{(\epsilon)}&
=\min\{V_p,\chi\}-V_{p-1},\\
\lim_{\epsilon\rightarrow 0} n\epsilon N_{j,p,>i,0}^{(\epsilon)}
&
=V_p-\max\{V_{p-1},\chi\},\\
\lim_{\epsilon\rightarrow 0} n\epsilon N_{j,p,<i,0}^{(\epsilon)}
&
=\min\{V_p,\chi\}-V_{p-1}.
\end{align*}

Also note that
\begin{align*}
\lim_{\epsilon\rightarrow 0}x_{q_{j,p,>i,1}}^{(\epsilon)}&
=\tau_j^{-1}e^{-\max\{V_{p-1},\chi\}}\mathbf{1}_{\{V_p>\chi,b_{p,j}=-,a_j=a_{i^*}\}},\\
\lim_{\epsilon\rightarrow 0}x_{q_{j,p,<i,1}}^{(\epsilon)}&
=\tau_je^{-V_{p-1}}\mathbf{1}_{\{V_{p-1}<\chi,b_{p,j}=+,a_j=a_{i^*}\}},\\
\lim_{\epsilon\rightarrow 0}x_{q_{j,p,>i,0}}^{(\epsilon)}&
=\tau_j^{-1}e^{-\max\{V_{p-1},\chi\}}\mathbf{1}_{\{V_{p}>\chi,b_{p,j}=-,a_j \neq a_{i^*}\}},\\
\lim_{\epsilon\rightarrow 0}x_{q_{j,p,<i,0}}^{(\epsilon)}&
=\tau_je^{-V_{p-1}}\mathbf{1}_{\{V_{p-1}<\chi,b_{p,j}=+,a_j \neq a_{i^*}\}}.
\end{align*}
This proves the lemma.
\end{proof}

For $w\in \CC$ and $\chi\in \RR$, define
\begin{align}
\mathcal{G}_{\chi}(w):=\mathcal{G}_{1,>{\chi}}(w)\cdot\mathcal{G}_{1,<{\chi}}(w)\cdot\mathcal{G}_{0,>\chi}(w)\cdot
\mathcal{G}_{0,<\chi}(w).\label{dgc}
\end{align}

We have the following asymptotic result.

\begin{theorem}\label{t57}Suppose Assumption \ref{ap5} holds. Assume $a_{i^{(\epsilon)}}=L$ for all $\epsilon>0$. Let $\mathrm{Pr}^{(\epsilon)}$ be the corresponding probability measure. Then
\begin{align}
\lim_{\epsilon\rightarrow 0}\mathbb{E}_{\mathrm{Pr}^{(\epsilon)}}\left[ \gamma_{k}(\lambda^{(i^{(\epsilon)})};t,t)\right]=\frac{1}{2\pi\mathbf{i}}\oint_{\mathcal{C}} \left[\mathcal{G}_{\chi}(w)\right]^{k\beta}\frac{dw}{w}\label{skm},
\end{align}
where the contour is positively oriented (which may be a union of disjoint simple closed curves) enclosing $0$ and every pole of $\mathcal{G}_{1,>{\chi}}$ and $\mathcal{G}_{0,>\chi}$, but does not enclose any other poles or zeros of 
$\mathcal{G}_{\chi}$; the expression
\begin{align*}
F^{k\beta}=e^{k\beta\log F}
\end{align*}
where the branch of $\log F$ is the one which takes positive real values when $F$ is positive and real.
\end{theorem}

\begin{proof}
By Lemma  \ref{lLeftMoments},
\begin{align*}
&\mathbb{E}_{\mathrm{Pr}^{(\epsilon)}}\left[ \gamma_{k}(\lambda^{(i^{(\epsilon)})};t,t)\right]
=\frac{1}{(2\pi\mathbf{i})^{k}}\oint_{\mathcal{C}_1^{(\epsilon)}}\ldots\oint_{\mathcal{C}_k^{(\epsilon)}} \frac{\sum_{i=1}^k\frac{1 }{w_i}}{\left(w_2-w_1\right)\ldots\left(w_k-w_{k-1}\right)}\prod_{i=1}^kdw_i \\
&\prod_{\substack{j\in[(l+1)^{(\epsilon)}..r^{(\epsilon)}],j\geq i^{(\epsilon)}\\b_j=-,a_j=a_{i^*}}}G_{1,>i^{(\epsilon)}}(W,x_j^{(\epsilon)},t)
\prod_{\substack{j\in[(l+1)^{(\epsilon)}..r^{(\epsilon)}],j\geq i^{(\epsilon)}\\b_j=-,a_j \neq a_{i^*}}}G_{0,>i^{(\epsilon)}}(W,x_j^{(\epsilon)},t)\\
&\times \prod_{\substack{j\in[l^{(\epsilon)}..r^{(\epsilon)}],j<i^{(\epsilon)}\\b_j=+,a_j= a_{i^*}}}G_{1,<i^{(\epsilon)}}(W,x_j^{(\epsilon)},t)
\prod_{\substack{j\in[l^{(\epsilon)}..r^{(\epsilon)}],j<i^{(\epsilon)}\\b_j=+;a_j\neq a_{i^*}}}G_{0,<i^{(\epsilon)}}(W,x_j^{(\epsilon)},t),
\end{align*}
where for $1\leq i\leq k$, $\mathcal{C}_i^{(\epsilon)}$ is  the integral contour for $w_i$, and $\mathcal{C}_1^{(\epsilon)},\ldots,\mathcal{C}_k^{(\epsilon)}$ satisfy the conditions as described in Lemma \ref{l31}. As $\epsilon\rightarrow 0$, assume that 
$\mathcal{C}_1^{(\epsilon)},\ldots,\mathcal{C}_k^{(\epsilon)}$ converge to contours $\mathcal{C}_1,\ldots,\mathcal{C}_k$, respectively, such that $\mathcal{C}_1,\ldots\mathcal{C}_k$ are separated from one another and do not cross any of the singularities of the integrand, by Lemma \ref{l55}, we obtain
\begin{align*}
&\lim_{\epsilon\rightarrow 0}\mathbb{E}_{\mathrm{Pr}^{(\epsilon)}}[ \gamma_{k}(\lambda^{(i^{(\epsilon)})};t,t)]=\frac{1}{(2\pi\mathbf{i})^{k}}\oint_{\mathcal{C}_1}\ldots\oint_{\mathcal{C}_k} \frac{\sum_{i=1}^k\frac{1 }{w_i}}{\left(w_2-w_1\right)\ldots\left(w_k-w_{k-1}\right)}\prod_{i=1}^kdw_i \\
&\times \prod_{w_s\in W}\left[\mathcal{G}_{1,>}\chi(w_x)\mathcal{G}_{1,<\chi}(w_s)\mathcal{G}_{0,>\chi}(w_s)\mathcal{G}_{0,<\chi}(w_s)\right]^{\beta}.
\end{align*}
Then (\ref{skm}) follows from Lemmas \ref{lb2} and $\ref{l55}$.
\end{proof}

\begin{theorem}\label{t58}Suppose Assumption \ref{ap5} holds. Assume that
\begin{align}
a_{i_1^{(\epsilon)}}=a_{i_2^{(\epsilon)}}=\ldots=a_{i_s^{(\epsilon)}}=L.\label{ael}
\end{align}
Let $\mathrm{Pr}^{(\epsilon)}$ be the corresponding probability measure and 
\begin{align*}
Q_{k_d}^{(\epsilon)}(\epsilon i_d^{(\epsilon)}):=\frac{1}{\epsilon}
\left(\gamma_{k_d}(\lambda^{(i_d^{(\epsilon)})};t,t)-\mathbb{E}_{\mathrm{Pr}^{(\epsilon)}}(\gamma_{k_d}(\lambda^{(i_d^{(\epsilon)})};t,t)\right).
\end{align*}
Then
$(Q_{k_1}^{(\epsilon)}(\epsilon i_1^{(\epsilon)}),\dots,Q_{k_s}^{(\epsilon)}(\epsilon i_s^{(\epsilon)}))
$ converges in distribution to the centered Gaussian vector 
\begin{align*}
(Q_{k_1}(\chi_1),\ldots,Q_{k_s}(\chi))
\end{align*}
as $\epsilon\rightarrow 0$, whose covariances are
\begin{align*}
\mathrm{Cov}\left[Q_{k_d}(\chi_d),Q_{k_h}(\chi_h)\right]
=\frac{k_dk_hn^2\beta^2}{(2\pi\mathbf{i})^2}\oint_{\mathcal{C}_d}\oint_{\mathcal{C}_h} \frac{\left[\mathcal{G}_{\chi_d}(z)\right]^{k_d\beta}\left[\mathcal{G}_{\chi_h}(w)\right]^{k_h\beta}}{(z-w)^2}dzdw,
\end{align*}
where
\begin{itemize}
    \item the $z$-contour ${\mathcal{C}_d}$ is positively oriented enclosing $0$ and every pole of $\mathcal{G}_{1,>\chi_d}\cup\mathcal{G}_{0,>\chi_d}$, but does not enclose any other poles or zeros of 
$\mathcal{G}_{\chi_d}(z)$;
    \item the $w$-contour ${\mathcal{C}_h}$ is positively oriented  enclosing $0$ and every pole of $\mathcal{G}_{1,>\chi_h}\cup\mathcal{G}_{0,>\chi_h}$, but does not enclose any other poles or zeros of $\mathcal{G}_{\chi_h}(w)$;
    \item the $z$-contour ${\mathcal{C}_d}$ and the $w$-contour ${\mathcal{C}_h}$ are disjoint;
    \item the branch of logarithmic function is chosen to take positive real values along the positive real axis.
  \end{itemize}
\end{theorem}

To prove Theorem \ref{t58}, we shall compute the moments of $Q_{k_d}^{(\epsilon)}(\epsilon i_d^{(\epsilon)})$ and show that these moments satisfy Wick's formula in the limit as $\epsilon\rightarrow 0$. We start with the following lemma about covariance.

\begin{lemma}\label{le59}Let $d,h\in[s]$. Under the assumptions of Theorem \ref{t58}, we have
\begin{align*}
&\lim_{\epsilon\rightarrow 0}\mathrm{Cov}\left[Q_{k_d}^{(\epsilon)}(\epsilon i_d^{(\epsilon)}),Q_{k_h}^{(\epsilon)}(\epsilon i_h^{(\epsilon)})\right]
=\frac{n^2\beta^2 k_dk_h}{(2\pi\mathbf{i})^2}\oint_{\mathcal{C}_d}\oint_{\mathcal{C}_h}
\frac{\left[\mathcal{G}_{\chi_d}(z)\right]^{k_d\beta}\left[\mathcal{G}_{\chi_h}(w)\right]^{k_h\beta}}{(z-w)^2}dzdw,
\end{align*}
where the $z$-contour $\mathcal{C}_d$ and the $w$-contour $\mathcal{C}_h$ satisfy the same conditions as in Theorem \ref{t58}.
\end{lemma}

\begin{proof}
Note that when $\mathrm{Cov}\left[Q_{k_d}^{(\epsilon)}(\epsilon i_d^{(\epsilon)})
Q_{k_h}^{(\epsilon)}(\epsilon i_h^{(\epsilon)})\right]$ is equal to
\begin{align*}
\frac{\mathbb{E}\left[\gamma_{k_d}(\lambda^{(i_d^{(\epsilon)})};t,t)
\gamma_{k_h}(\lambda^{(i_h^{(\epsilon)})};t,t)\right]-\mathbb{E}\left[\gamma_{k_d}(\lambda^{(i_d^{(\epsilon)})};t,t)\right]\mathbb{E}_{\mathrm{Pr}}
\left[\gamma_{k_h}(\lambda^{(i_h^{(\epsilon)})};t,t)\right]}{\epsilon^2}.
\end{align*}
For $W=(w_1,\dots,w_k)$ we use abbreviation
\begin{align*}
F_i^{(\epsilon)}(W)&=\prod_{\substack{j\in[(l+1)^{(\epsilon)}..r^{(\epsilon)}]\\j\geq i,b_j=-,a_j=a_i}}G_{1,>}(W,x_j^{(\epsilon)},t)\prod_{\substack{j\in[(l+1)^{(\epsilon)}..r^{(\epsilon)}]\\j\geq i,b_j=-,a_j \neq a_i}}G_{0,>}(W,x_j^{(\epsilon)},t)\\
&\times\prod_{\substack{j\in[l^{(\epsilon)}..r^{(\epsilon)}]\\j<i,b_j=+;a_i= a_j}}G_{1,<}(W,x_j^{(\epsilon)},t)\prod_{\substack{j\in[l^{(\epsilon)}..r^{(\epsilon)}]\\j<i,b_j=+;a_i \neq a_j}}G_{0,<}(W,x_j^{(\epsilon)},t)\\
&\times \frac{\sum_{j=1}^{k}\frac{1}{w_j}}{\left(w_2-w_1\right)\ldots\left(w_{k}-w_{k-1}\right)}\prod_{1\leq i<j\leq k }\frac{\left(1-\frac{w_i}{w_j}\right)^2}{\left(1-\frac{w_i}{tw_j}\right)\left(1-\frac{tw_i}{w_j}\right)}.
\end{align*}

Note that for $\epsilon i^{(\epsilon)}\to \chi$ and $t\to 1$
\begin{align}
\lim_{\epsilon \to 0}F_{i^{(\epsilon)}}^{(\epsilon)}(W)=\prod_{i=1}^k\mathcal{G}_\chi(w_i) \frac{\sum_{j=1}^{k}\frac{1}{w_j}}{\left(w_2-w_1\right)\ldots\left(w_{k}-w_{k-1}\right)}\label{limF}.
\end{align}

By Lemma \ref{lLeftMoments}, we obtain that  $\mathrm{Cov}\left[Q_{k_d}^{(\epsilon)}(\epsilon i_d^{(\epsilon)}),Q_{k_h}^{(\epsilon)}(\epsilon i_h^{(\epsilon)})\right]$ is equal to
\begin{align*}
\frac{1}{\epsilon^2(2\pi \mathbf{i})^{k_h+k_d}}\oint_{\mathcal{C}_{1,1}^{(\epsilon)}}\ldots\oint_{\mathcal{C}_{1,k_d}^{(\epsilon)}} 
\oint_{\mathcal{C}_{2,1}^{(\epsilon)}}\ldots\oint_{\mathcal{C}_{2,k_h}^{(\epsilon)}}& \prod_{\xi \in\{ {d,h}\}}\prod_{i=1}^{k_\xi}dw^{(i_\xi^{(\epsilon)})}_iF^{(\epsilon)}_{i_\xi^{(\epsilon)}}(W^{(i_\xi^{(\epsilon)})})\\
&\times\left[T_{L,L}(W^{(i_d^{(\epsilon)})},W^{(i_h^{(\epsilon)})})-1\right],
\end{align*}
where $T_{L,L}(Z,W)$ is given by (\ref{tll}), $|W^{(i_d^{(\epsilon)})}|=k_d$, $|W^{(i_h^{(\epsilon)}}|=k_h$ and for $1\leq i\leq k_d$ (resp.\ $1\leq j\leq k_h$), $\mathcal{C}_{1,i}^{(\epsilon)}$ (resp.\ $\mathcal{C}_{2,j}^{(\epsilon)}$) is  the integral contour for $w_i^{(i_d^{(\epsilon)})}$ (resp.\ 
$w_j^{(i_h^{(\epsilon)})}$
), and $\mathcal{C}_{1,1}^{(\epsilon)},\ldots,\mathcal{C}_{1,k_d}^{(\epsilon)}$, 
$\mathcal{C}_{2,1}^{(\epsilon)},\ldots,\mathcal{C}_{2,k_h}^{(\epsilon)}$
satisfy the conditions as described in Lemma \ref{l31}. As $\epsilon\rightarrow 0$, assume that 
$\mathcal{C}_1^{(\epsilon)},\ldots,\mathcal{C}_k^{(\epsilon)}$ converge to contours $\mathcal{C}_1,\ldots,\mathcal{C}_k$, respectively, such that $\mathcal{C}_1,\ldots\mathcal{C}_k$ are separated from one another and do not cross any of the singularities of the integrand.

Note that
\begin{align*}
T_{L,L}(Z,W)-1=\sum_{\substack{\emptyset\neq S\\S \subset [k_d]\times[k_h]}}\prod_{(i,j)\in S}\frac{(1-t)(t^{-1}-1)z_iw_j}{\left(z_i-t^{-1}w_j \right)
\left(z_i-tw_j \right)}.
\end{align*}
Under Assumption \ref{ap5}, we obtain
\begin{align}
\frac{1}{\epsilon^2}(1-t)(t^{-1}-1)=n^2\beta^2+O(\epsilon)\label{lmt}.
\end{align}
Therefore 
\begin{align*}
&\frac{1}{\epsilon^2}\left[T_{L,L}(W^{(i_d^{(\epsilon)})},W^{(i_h^{(\epsilon)})})-1\right]\\
&=n^2\beta^2\left[\sum_{(u,v)\in [k_d]\times  [k_h] }\frac{w_{u}^{(i_d^{(\epsilon)})}
w_{v}^{(i_h^{(\epsilon)})}
}{(w_{u}^{(i_d^{(\epsilon)})}-t^{-1}
w_{v}^{(i_h^{(\epsilon)})}
)
(w_{u}^{(i_d^{(\epsilon)})}-t
w_{v}^{(i_h^{(\epsilon)})}
)}\right]+o(\epsilon).
\end{align*}
In the above formula the main contribution comes from one-element subsets $S$, while the others have a negligible $o(\epsilon)$ contribution. This together with (\ref{limF}) gives that in the limit  $\mathrm{Cov}\left[Q_{k_d}^{(\epsilon)}(\epsilon i_d^{(\epsilon)}),
Q_{k_h}^{(\epsilon)}(\epsilon i_h^{(\epsilon)})\right]$ is equal to the limit of
\begin{align*}
&\frac{n^2\beta^2}{(2\pi\mathbf{i})^{k_d+k_h}}\oint_{\mathcal{C}_{1,1}^{(\epsilon)}}\ldots\oint_{\mathcal{C}_{1,k_d}^{(\epsilon)}} 
\oint_{\mathcal{C}_{2,1}^{(\epsilon)}}\ldots\oint_{\mathcal{C}_{2,k_h}^{(\epsilon)}} \prod_{\xi \in\{ {d,h}\}}\prod_{i=1}^{k_\xi}dw^{(i_\xi^{(\epsilon)})}_i\\
&\times \prod_{\xi \in \{d,h\}}\left[\prod_{i=1}^{k_{\xi}}\mathcal{G}_{\chi_{\xi}}(w_i^{(i^{(\epsilon)}_\xi)})\right] \frac{\sum_{j=1}^{k_{\xi}}{\Bigl[w_j^{(i^{(\epsilon)}_\xi)}\Bigr]}^{-1}}{(w_2^{(i^{(\epsilon)}_\xi)}-w_1^{(i^{(\epsilon)}_\xi)})\ldots(w_{k_{\xi}}^{(i^{(\epsilon)}_\xi)}-w_{k_{\xi}-1}^{(i^{(\epsilon)}_\xi)})}\\
&\times\left[\sum_{(u,v)\in [k_d]\times  [k_h] }\frac{w_{u}^{(i_d^{(\epsilon)})}
w_{v}^{(i_h^{(\epsilon)})}
}{(w_{u}^{(i_d^{(\epsilon)})}-t^{-1}
w_{v}^{(i_h^{(\epsilon)})}
)
(w_{u}^{(i_d^{(\epsilon)})}-t
w_{v}^{(i_h^{(\epsilon)})}
)}+o(\epsilon)\right].
\end{align*}

Then by Lemma \ref{lb2}, we have
\begin{align*}
&\lim_{\epsilon\rightarrow 0}\frac{1}{(2\pi\mathbf{i})^{k_h}} 
\oint_{\mathcal{C}_{2,1}^{(\epsilon)}}\ldots\oint_{\mathcal{C}_{2,k_h}^{(\epsilon)}} \prod_{i=1}^{k_h}dw^{(i_h^{(\epsilon)})}_i\\
&\times \left[\prod_{i=1}^{k_{h}}\mathcal{G}_{\chi_h}(w_i^{(i^{(\epsilon)}_h)})\right] \frac{\sum_{j=1}^{k_h}{\Bigl[w_j^{(i^{(\epsilon)}_h)}\Bigr]}^{-1}}{(w_2^{(i^{(\epsilon)}_h)}-w_1^{(i^{(\epsilon)}_h)})\ldots(w_{k_{h}}^{(i^{(\epsilon)}_h)}-w_{k_{h}-1}^{(i^{(\epsilon)}_h)})}\\
&\times \left[\sum_{(u,v)\in [k_d]\times  [k_h] }\frac{w_{u}^{(i_d^{(\epsilon)})}
w_{v}^{(i_h^{(\epsilon)})}
}{(w_{u}^{(i_d^{(\epsilon)})}-t^{-1}
w_{v}^{(i_h^{(\epsilon)})}
)
(w_{u}^{(i_d^{(\epsilon)})}-t
w_{v}^{(i_h^{(\epsilon)})}
)}+o(\epsilon)\right]\\
&=\frac{k_h}{2\pi \mathbf{i}}\oint_{\mathcal{C}_h}\left[
\mathcal{G}_{\chi_h}(w)\
\right]^{k_h\beta}\times\lim_{\epsilon\rightarrow 0}\left[\sum_{u\in [k_d]}\frac{w_{u}^{(i_d^{(\epsilon)})}
w
}{(w_{u}^{(i_d^{(\epsilon)})}-w
)
(w_{u}^{(i_d^{(\epsilon)})}-w
)}\right]  \frac{dw}{w}.
\end{align*}

Applying Lemma \ref{lb2} again to integrals over $\mathcal{C}_{1,1}^{(\epsilon)}$, \ldots, $\mathcal{C}_{1,k_d}^{(\epsilon)}$, we obtain the result.
\end{proof}

\begin{lemma}\label{l510}Assume the assumptions of Theorem \ref{t58} hold. 
\begin{enumerate}
\item Let $s\in \NN$ be odd, and $s\geq 3$. Then
\begin{align*}
&\lim_{\epsilon\rightarrow 0}\mathbb{E}_{\mathrm{Pr}^{(\epsilon)}}\left[\prod_{u=1}^{s}Q_{k_{u}}^{(\epsilon)}(\epsilon i_u^{(\epsilon)})\right]=0.
\end{align*}
\item If $s\in \NN$ is even, then
\begin{align*}
&\lim_{\epsilon\rightarrow 0}\mathbb{E}_{\mathrm{Pr}^{(\epsilon)}}\left[\prod_{u=1}^{s}Q_{k_{u}}^{(\epsilon)}(\epsilon i_u^{(\epsilon)})\right]=\lim_{\epsilon\rightarrow 0}\sum_{P\in\mathcal{P}_{s}^2}\prod_{u,v\in P}\mathbb{E}_{\mathrm{Pr}^{(\epsilon)}}\left[
Q_{k_{u}}^{(\epsilon)}(\epsilon i_u^{(\epsilon)})
Q_{k_{v}}^{(\epsilon)}(\epsilon i_v^{(\epsilon)})
\right],
\end{align*}
where the sum runs over all pairings of $[s]$.
\end{enumerate}
\end{lemma}

\begin{proof}
Note that
\begin{align*}
&\mathbb{E}_{\mathrm{Pr}^{(\epsilon)}}\left[\prod_{u=1}^{s}Q_{k_{u}}^{(\epsilon)}(\epsilon i_u^{(\epsilon)})\right]\\
&=\frac{1}{\epsilon^s}\sum_{J\subseteq [s]}(-1)^{\left|[s]\setminus J\right|}\left[\mathbb{E}_{\mathrm{Pr}^{(\epsilon)}} \prod_{j\in J}\gamma_{k_j}(\lambda^{(i_j^{(\epsilon)})})\right]\prod_{u\in[s]\setminus J}\mathbb{E}_{\mathrm{Pr}^{(\epsilon)}}
\gamma_{k_u}(\lambda^{(i_u^{(\epsilon)})}).
\end{align*}
When computing $\mathbb{E}\left[\prod_{u=1}^{s}Q_{k_{u}}^{(\epsilon)}(\epsilon i_u^{(\epsilon)})\right]$ by Lemma \ref{lLeftMoments}, in the integrand there is a factor
\begin{align}
\sum_{J\subseteq [s]}(-1)^{\left|[s]\setminus J\right|}\prod_{u<v;u,v\in[s]\setminus J}T_{L,L}(W^{(i_u^{(\epsilon)})},W^{(i_v^{(\epsilon)})})\label{fct},
\end{align}
which is by (\ref{tll}) equal to
\begin{align*}
\sum_{J\subseteq [s]}(-1)^{\left|[s]\setminus J\right|}\prod_{\substack{u,v\in[s]\setminus J\\u<v}}\prod_{\substack{w_j^{(u)}\in W^{(i_u^{(\epsilon)})}\\w_f^{(v)}\in 
W^{(i_v^{(\epsilon)})}}}\left[1+\frac{(1-t)(t^{-1}-1)w_j^{(u)}w_f^{(v)}}{(w_j^{(u)}-t^{-1}w_f^{(v)})
(w_j^{(u)}-tw_f^{(v)})}\right].
\end{align*}

Under the assumptions, (\ref{lmt}) is true. Let
\begin{align*}
\frac{1}{\epsilon^2}\frac{(1-t)(t^{-1}-1)w_j^{(u)}w_f^{(v)}}{\left(w_j^{(u)}-t^{-1}w_f^{(v)}\right)
(w_j^{(u)}-tw_f^{(v)})}=C_{u,v,j,f}^{(\epsilon)},
\end{align*}
where $C^{(\epsilon)}_{u,v,j,f}$ tends to a constant as $\epsilon \to 0$. Let
\begin{align*}
K_J:=\left\{(u,v,j,f):u<v;u,v\in[s]\setminus J, w_j^{(u)}\in W^{(i_u^{(\epsilon)})},w_f^{(v)}\in 
W^{(i_v^{(\epsilon)})}\right\}.
\end{align*}
Then (\ref{fct}) is equal to
\begin{align}
\sum_{J\subseteq [s]}(-1)^{\left|[s]\setminus J\right|}+
\sum_{J\subseteq [s]}(-1)^{\left|[s]\setminus J\right|}\left(\sum_{\emptyset\neq H\subseteq K_J}\epsilon^{2|H|}\prod_{(u,v,j,f)\in H}C_{u,v,j,f}^{(\epsilon)}\right).\label{ep5}
\end{align}
Note that $\sum_{J\subseteq [s]}(-1)^{\left|[s]\setminus J\right|}=0$. For each fixed $H\subseteq K_{\emptyset}$, if $H\neq K_{\emptyset}$, let
\begin{align*}
H_0:=\{u\in[s]: \exists v\in[s]\ \mathrm{and } j,f,\ \mathrm{s.t.} (u,v,j,f)\in H,\ \mathrm{or}\ (v,u,j,f)\in H\}.
\end{align*}
The sum of terms with $\prod_{(u,v,j,f)\in H}C_{u,v,j,f}^{(\epsilon)}$ in (\ref{ep5}) is
\begin{align}
\epsilon^{2|H|}\prod_{(u,v,j,f)\in H}C_{u,v,j,f}^{(\epsilon)}\sum_{J\in [s]:J\cap H_0=\emptyset}(-1)^{|[s]\setminus J|}.\label{smh}
\end{align}
As long as $H_0\neq S$, the sum of $(-1)^{|[s]\setminus J|}$ over all the subsets of $[s]\setminus H_0$ is 0. Therefore (\ref{ep5}) is equal to 
\begin{align*}
\sum_{\emptyset\neq H\subseteq K_{\emptyset}, H_0=[s]}\epsilon^{2|H|}\prod_{(u,v,j,f)\in H}C_{u,v,j,f}^{(\epsilon)}.
\end{align*}
\begin{enumerate}
\item If $s$ is odd, as $\epsilon\rightarrow 0$,
\begin{align*}
\sum_{\emptyset\neq H\subseteq K_{\emptyset}, H_0=[s]}\epsilon^{2|H|}\prod_{(u,v,j,f)\in H}C_{u,v,j,f}^{(\epsilon)}=O(\epsilon^{s+1});
\end{align*}
therefore
\begin{align*}
\lim_{\epsilon\rightarrow 0}\frac{1}{\epsilon^s}\sum_{J\subseteq [s]}(-1)^{\left|[s]\setminus J\right|}\prod_{u<v;u,v\in[s]\setminus J}T_{L,L}(W^{(i_u^{(\epsilon)})},W^{(i_v^{(\epsilon)})})=0.
\end{align*}
\item If $s$ is even, as $\epsilon\rightarrow 0$,
\begin{align*}
&\sum_{\emptyset\neq H\subseteq K_{\emptyset}, H_0=[s]}\epsilon^{2|H|}\prod_{(u,v,j,f)\in H}C_{u,v,j,f}^{(\epsilon)}\\
&=\epsilon^{s}\sum_{P\in \mathcal{P}_s^2}\prod_{(u,v)\in \mathcal{P}}\prod_{{(j,f):w_j^{(u)}\in W^{(i_u^{(\epsilon)})},
w_j^{(v)}\in W^{(i_v^{(\epsilon)})}}}C_{u,v,j,f}^{(\epsilon)}
+ O(\epsilon^{s+1})\\
&=\epsilon^{s}\sum_{P\in \mathcal{P}_s^2}\prod_{(u,v)\in \mathcal{P}}
\left(\frac{1}{\epsilon^2}T_{L,L}(W^{(i_u^{(\epsilon)})},W^{(i_v^{(\epsilon)})})\right).
\end{align*}
\end{enumerate}
Then the lemma follows.
\end{proof}

\noindent{\textbf{Proof of Theorem \ref{t58}}.} The theorem follows from Lemmas \ref{le59} and \ref{l510}.

\section{Frozen Boundary}\label{sect:fb}

In this section, we prove an integral formula for the Laplace transform of the rescaled height function (see Theorem \ref{t61}), which turns out to be deterministic, as a 2D analog of the law of large numbers. We further obtain an explicit formula for the frozen boundary in the scaling limit.

\begin{theorem}\label{t61}
Let $M$ be a random pure dimer covering on the rail yard graph $RYG(l,r,\underline{a},\underline{b})$ with probability distribution given by (\ref{ppd}) and (\ref{pbc}). Let $h_M$ be the height function associated to $M$ as defined in (\ref{dhm}).
Suppose Assumption \ref{ap5} holds. 
Then the rescaled random height function $\epsilon h_M\left(\frac{\chi}{\epsilon},\frac{\kappa}{\epsilon}\right)$ converges, as $\epsilon\rightarrow 0$, to a non-random function $\mathcal{H}(\chi,\kappa)$ such that the Laplace transform of $\mathcal{H}(\chi,\cdot)$ is given by
\begin{align}
\int_{-\infty}^{\infty} e^{-n\alpha \kappa}\mathcal{H}(\chi,\kappa)d\kappa
=\frac{1}{n^2\alpha^2\pi\mathbf{i}}\oint_{\mathcal{C}} \left[\mathcal{G}_{\chi}(w)\right]^{\alpha}\frac{dw}{w},\label{lph}
\end{align}
where $\alpha$ is a positive real number and the contour $\mathcal{C}$ satisfies the conditions of Theorem \ref{t57}. Here 
\begin{align*}
\left[\mathcal{G}_{\chi}(w)\right]^{\alpha}=e^{\alpha \log [\mathcal{G}_{\chi}(w)]}
\end{align*}
and the branch of $\log(\zeta) $ is chosen to be real positive when $\zeta$ is real positive.
Note that the right hand side is non-random. 
\end{theorem}

\begin{proof}
By Theorem \ref{t57}, (\ref{hme}) and (\ref{dgm}), let $\frac{\kappa}{\epsilon}=y$, we obtain
\begin{align*}
&\lim_{\epsilon\rightarrow 0}\mathbb{E}_{\mathrm{Pr}^{(\epsilon)}}\int_{-\infty}^{\infty} e^{-n\beta \kappa k}\epsilon h_{M}\left(\frac{\chi}{\epsilon},\frac{\kappa}{\epsilon}\right)d\kappa=\lim_{\epsilon\rightarrow 0}\epsilon^2\mathbb{E}_{\mathrm{Pr}^{(\epsilon)}}\int_{-\infty}^{\infty} h_{M}\left(\frac{\chi}{\epsilon},y\right)t^{ky}dy
\\
&=\lim_{\epsilon\rightarrow 0}\mathbb{E}_{\mathrm{Pr}^{(\epsilon)}}\frac{2\epsilon^2}{(k\log t)^2}\gamma_{k}(\lambda^{(i^{(\epsilon)})};t,t)=\frac{1}{k^2n^2\beta^2\pi\mathbf{i}}\oint_{\mathcal{C}} \left[\mathcal{G}_{x}(w)\right]^{k\beta}\frac{dw}{w}.
\end{align*}
To show that the limit, as $\epsilon\rightarrow 0$, of $\int_{-\infty}^{\infty} e^{-n\beta \kappa k}\epsilon h_{M}\left(\frac{\chi}{\epsilon},\frac{\kappa}{\epsilon}\right)d\kappa$, is non-random, it suffices to show that the limit of its variance is 0. Note that
\begin{align*}
&\lim_{\epsilon\rightarrow 0}\mathbb{E}_{\mathrm{Pr}^{(\epsilon)}}\left[\int_{-\infty}^{\infty} e^{-n\beta \kappa k}\epsilon h_{M}\left(\frac{\chi}{\epsilon},\frac{\kappa}{\epsilon}\right)d\kappa-\mathbb{E}_{\mathrm{Pr}^{(\epsilon)}}\int_{-\infty}^{\infty} e^{-n\beta \kappa k}\epsilon h_{M}\left(\frac{\chi}{\epsilon},\frac{\kappa}{\epsilon}\right)d\kappa\right]^2\\
&=\lim_{\epsilon\rightarrow 0}\epsilon^4\mathbb{E}_{\mathrm{Pr}^{(\epsilon)}}\left[\int_{-\infty}^{\infty} e^{-n\beta \kappa k} h_{M}\left(\frac{\chi}{\epsilon},y\right)d y-\mathbb{E}_{\mathrm{Pr}^{(\epsilon)}}\int_{-\infty}^{\infty} e^{-n\beta \kappa k}\epsilon h_{M}\left(\frac{\chi}{\epsilon},y\right)dy\right]^2\\
&=\lim_{\epsilon\rightarrow 0}\frac{\epsilon^6}{(k\log t)^4}\mathrm{Var}\left[Q_k^{(\epsilon)}(\chi)\right]=0,
\end{align*}
where the last identity follows from Lemma \ref{le59} and the limit regime stated in Assumption \ref{ap5}. Let $\alpha=k\beta$ and consider analytic continuation if necessary, then the lemma follows.
\end{proof}

By (\ref{lph}), for $\alpha>0$ we obtain 
\begin{align}
\int_{-\infty}^{\infty} e^{-n\alpha \kappa}\frac{\partial\mathcal{H}(\chi,\kappa)}{\partial \kappa}d\kappa
&=n\alpha\int_{-\infty}^{\infty} e^{-n\alpha \kappa}\mathcal{H}(\chi,\kappa)d\kappa
=
\frac{1}{n\alpha\pi\mathbf{i}}\oint_{\mathcal{C}} \left[\mathcal{G}_{\chi}(w)\right]^{\alpha}\frac{dw}{w}.\label{dedh}
\end{align}

Let $\mathbf{m}_{\chi}$ be the measure on $(0,\infty)$ defined by 
\begin{align*}
\mathbf{m}_{\chi}(ds)=\left.e^{-\kappa}\frac{\partial\mathcal{H}(\chi,\kappa)}{\partial \kappa}|d\kappa|\right|_{\kappa=-\ln s}.
\end{align*}
We are particularly interested in the measure $\mathbf{m}_{\chi}$ because its density with respect to the Lebesgue measure on $\RR$ is given by
\begin{align}
\frac{\mathbf{m}_{\chi}(ds)}{ds}=\left.\frac{\partial\mathcal{H}(\chi,\kappa)}{\partial \kappa}\right|_{\kappa=-\ln s};\label{slht}
\end{align}
which is exactly the slope of the limiting rescaled height function in the $\kappa$-direction when $s=e^{-\kappa}$.

By (\ref{dedh})
we deduce that for any $\chi\in (l^{(0)},r^{(0)})$
$\int_{0}^{\infty}\mathbf{m}_{\chi}(d s)<\infty$, i.e. $\mathbf{m}_{\chi}(d s)$ is a measure on $\RR$ with finite total mass. Note also that for any positive integer $j$, by (\ref{dedh}) we obtain
\begin{align*}
\int_0^{\infty} s^{j-1}\mathbf{m}_{\chi}(ds)=\int_{-\infty}^{\infty}e^{-\kappa j}\frac{\partial\mathcal{H}(\chi,\kappa)}{\partial \kappa}d\kappa=\frac{1}{j\pi \mathbf{i}}\oint_{\mathcal{C}}[\mathcal{G}_{\chi}(w)]^{\frac{j}{n}}\frac{dw}{w}
\leq C^j,
\end{align*}
where $C>0$ is a positive constant independent of $j$. Hence we obtain
\begin{align*}
\int_{2C}^{\infty}\mathbf{m}_{\chi}(ds)\leq \frac{\int_{2C}^{\infty}s^{j-1}\mathbf{m}_{\chi}(ds)}{(2C)^{j-1}}\leq \frac{1}{2C}\left(\frac{1}{2}\right)^j\longrightarrow 0
\end{align*}
as $j\rightarrow\infty$. Hence we obtain that $\mathbf{m}_{\chi}(ds)$ has compact support in $(0,\infty)$. 

We shall now compute the density of the measure $\mathbf{m}_{\chi}(ds)$ with respect to the Lebesgue measure on $\RR$. It is a classical fact about Stieltjes transform that
\begin{align*}
\frac{\mathbf{m}_{\chi}(ds)}{ds}
=-\lim_{\epsilon\rightarrow 0+}\frac{1}{\pi}\Im\left(\mathrm{St}_{\mathbf{m}_{\chi}}(s+\mathbf{i}\epsilon)\right),
\end{align*}
where $\Im$ denotes the imaginary part of a complex number and $\mathrm{St}_{\mathbf{m}_{\chi}}$ is the Stieltjes transform of the measure $\mathbf{m}_{\chi}$, which can be computed as follows: 
for $\zeta\in \CC\setminus\mathrm{supp}(\mathbf{m}_{\chi})$,
\begin{align}
\mathrm{St}_{\mathbf{m}_{\chi}}(\zeta)=\int_{0}^{\infty}\frac{\mathbf{m}_{\chi}(ds)}{\zeta-s}
=\sum_{i=0}^{\infty}\int_{0}^{\infty}\frac{s^i\mathbf{m}_{\chi}(ds)}{\zeta^{i+1}}
=\sum_{j=1}^{\infty}\frac{1}{\zeta^j}\int_0^{\infty}e^{-j\kappa}\frac{\partial\mathcal{H}(\chi,\kappa)}{\partial\kappa}d\kappa.\label{stc}
\end{align}
Again by (\ref{dedh}) we obtain
\begin{align*}
\mathrm{St}_{\mathbf{m}_{\chi}}(\zeta)=\sum_{j=1}^{\infty}\frac{1}{\zeta^j j\pi\mathbf{i}}
\oint_{\mathcal{C}}[\mathcal{G}_{\chi}(w)]^{\frac{j}{n}}\frac{dw}{w}.
\end{align*}
Let $\mathcal{R}_\chi$ denotes the set of all poles of $\mathcal{G}_{1,>\chi}$ and $\mathcal{G}_{0,>\chi}$. When the contour $\mathcal{C}$ satisfies the conditions given as in Theorem \ref{t57}, we can split $\mathcal{C}$ into a positively oriented simple closed curve $\mathcal{C}_0$ enclosing only 0, and a union $\mathcal{C}_1$ of positively oriented simple closed curves enclosing every point in $\mathcal{R}_\chi$. By the residue theorem, we obtain that 
\begin{align*}
\sum_{j=1}^{\infty}\frac{1}{\zeta^j j\pi\mathbf{i}}
\oint_{\mathcal{C}_0}[\mathcal{G}_{\chi}(w)]^{\frac{j}{n}}\frac{dw}{w}=
\sum_{j=1}^{\infty}\frac{2}{\zeta^j j}
[\mathcal{G}_{\chi}(0)]^{\frac{j}{n}}=-2\log\left(1-\frac{[\mathcal{G}_{\chi}(0)]^{\frac{1}{n}}}{\zeta}\right).
\end{align*}
Moreover,
\begin{align*}
\sum_{j=1}^{\infty}\frac{1}{\zeta^j j\pi\mathbf{i}}
\oint_{\mathcal{C}_1}[\mathcal{G}_{\chi}(w)]^{\frac{j}{n}}\frac{dw}{w}
=-\frac{1}{\pi \mathbf{i}}\oint_{\mathcal{C}_1}\log\left(1-\frac{[\mathcal{G}_{\chi}(w)]^{\frac{1}{n}}}{\zeta}\right)\frac{dw}{w},
\end{align*}
where $\left|\frac{[\mathcal{G}_{\chi}(0)]^{\frac{1}{n}}}{\zeta}\right|<1$ and 
$\max_{\zeta\in\mathcal{C}_1}\left|\frac{[\mathcal{G}_{\chi}(\zeta)]^{\frac{1}{n}}}{\zeta}\right|<1$ to ensure the convergence of Maclaurin series.
Hence we have
\begin{align*}
\mathrm{St}_{\mathbf{m}_{\chi}}(\zeta)=-2\log\left(1-\frac{[\mathcal{G}_{\chi}(0)]^{\frac{1}{n}}}{\zeta}\right)-\frac{1}{\pi \mathbf{i}}\oint_{\mathcal{C}_1}\log\left(1-\frac{[\mathcal{G}_{\chi}(w)]^{\frac{1}{n}}}{\zeta}\right)\frac{dw}{w}.
\end{align*}

We would like to get rid of the fractal exponent for the simplicity of computing complex integrals. To that end, we define another function
\begin{align}
\Theta_{\chi}(\zeta):=-2\log\left(1-\frac{\mathcal{G}_{\chi}(0)}{\zeta^n}\right)-\frac{1}{\pi \mathbf{i}}\oint_{\mathcal{C}_1}\log\left(1-\frac{\mathcal{G}_{\chi}(w)}{\zeta^n}\right)\frac{dw}{w}.\label{dtx}
\end{align}
Let $\omega=e^{\frac{2\pi\mathbf{i}}{n}}$, then it is straightforward to check that
$
\Theta_{\chi}(\zeta)=\sum_{i=0}^{n-1}\mathrm{St}_{\mathbf{m}_{
\chi}}(\omega^{-i}\zeta)$. Then we obtain
\begin{align*}
\frac{1}{\pi}\lim_{\epsilon\rightarrow 0+}\Im \Theta_{\chi}(s+\mathbf{i}\epsilon)
=\frac{1}{\pi}\lim_{\epsilon\rightarrow 0+}\sum_{i=0}^{n-1}\Im \mathrm{St}_{\mathbf{m}_{\chi}}(\omega^{-i}(s+\mathbf{i}\epsilon)).
\end{align*}
Since $\mathrm{St}_{\mathbf{m}_{\chi}}(\zeta)$ is continuous in $\zeta$ when $\zeta\in \CC\setminus \mathrm{supp}(\mathbf{m}_{\chi})$, $\mathrm{supp}(\mathbf{m}_{\chi})\in(0,\infty)$, and $\mathrm{St}_{\mathbf{m}_{\chi}}(\overline{\zeta})=
\overline{\mathrm{St}_{\mathbf{m}_{\chi}}(\zeta)}$, we obtain that when $s\in \mathrm{supp}(\mathbf{m}_{\chi})$,
\begin{align}
\label{dsm2}\frac{1}{\pi}\lim_{\epsilon\rightarrow 0+}\Im \Theta_{\chi}(s+\mathbf{i}\epsilon)
&=\frac{1}{\pi}\sum_{i=1}^{n-1}\Im \mathrm{St}_{\mathbf{m}_{\chi}}(\omega^{-i}s)+
\frac{1}{\pi}\lim_{\epsilon\rightarrow 0+}\Im\mathrm{St} _{\mathbf{m}_{\chi}}(s+\mathbf{i}\epsilon)\\
&=\frac{1}{\pi}\lim_{\epsilon\rightarrow 0+}\Im\mathrm{St} _{\mathbf{m}_{\chi}}(s+\mathbf{i}\epsilon)=-\frac{\mathbf{m}_{\chi}(ds)}{ds}.\notag
\end{align}

Hence by (\ref{slht}), to compute the slope of the limiting rescaled height function in the $\kappa$-direction, it suffices to compute 
$-\frac{1}{\pi}\lim_{\epsilon\rightarrow 0+}\Im \Theta_{\chi}(s+\mathbf{i}\epsilon)$ when $s=e^{-\kappa}$.

By (\ref{dtx}) we obtain
\begin{align}
\label{ctx}\Theta_{\chi}(\zeta):&=-2\log\left(1-\frac{\mathcal{G}_{\chi}(0)}{\zeta^n}\right)-\frac{1}{\pi \mathbf{i}}\oint_{\mathcal{C}_1}\log\left(1-\frac{\mathcal{G}_{\chi}(w)}{\zeta^n}\right)d\log w\\
&=-2\log\left(1-\frac{\mathcal{G}_{\chi}(0)}{\zeta^n}\right)-\frac{1}{\pi \mathbf{i}}\oint_{\mathcal{C}_1}d\left[\log\left(1-\frac{\mathcal{G}_{\chi}(w)}{\zeta^n}\right)\log w\right]\notag\\
&-\frac{1}{\pi \mathbf{i}}\oint_{\mathcal{C}_1}
\frac{\mathcal{G}'_{\chi}(w)\log w}{\zeta^n-\mathcal{G}_{\chi}(w)}dw.\notag
\end{align}

To compute the contour integral above, we need to consider the root of the following equation in $w$:
\begin{align}
\mathcal{G}_{\chi}(w)=\zeta^n;\label{eeq}
\end{align}
in particular, the roots of (\ref{eeq}) that are enclosed by the contour $\mathcal{C}_1$. Recall that $\mathcal{C}_1$ is the union of positively oriented simple closed curves enclosing every point in $\mathcal{R}_{\chi}$, but no other poles or zeros $\mathcal{G}_{\chi}$. We may assume
\begin{align*}
\mathcal{C}_1:=\cup_{\xi\in \mathcal{R}}\mathcal{C}_{\xi};
\end{align*}
where $\mathcal{C}_{\xi}$ is a positively oriented simple closed curve enclosing $\xi$ but no other poles or zeros of $\mathcal{G}_{\chi}$.

When $\zeta\rightarrow \infty$, zeros of (\ref{eeq}) will approach poles of $\mathcal{G}_{\chi}$. For each $\xi\in \mathcal{R}_{\chi}$, let $w_{\xi,\chi}(\zeta)$ be a root of (\ref{eeq}) such that $\lim_{\zeta\rightarrow\infty}w_{\xi,\chi}(\zeta)=\xi$.

When $|\zeta|$ is sufficiently large, $w_{\xi,\chi}(\zeta)$ is enclosed by $\mathcal{C}_{\xi}$. Enclosed by each $\mathcal{C}_{\xi}$, there is exactly one zero and one pole for $1-\frac{\mathcal{G}_{\chi}(w)}{\zeta^n}$, hence
\begin{align}
\frac{1}{\pi \mathbf{i}}\oint_{\mathcal{C}_1}d\left[\log\left(1-\frac{\mathcal{G}_{\chi}(w)}{\zeta^n}\right)\log w\right]=0.\label{ctx1}
\end{align}
By computing residues at each $w_{\xi,\chi}(\zeta)$ and $\xi$, we obtain 
\begin{align}
-\frac{1}{\pi \mathbf{i}}\oint_{\mathcal{C}_1}
\frac{\mathcal{G}'_{\chi}(w)\log w}{\zeta^n-\mathcal{G}_{\chi}(w)}dw
=2\sum_{\xi\in \mathcal{R}_{\chi}}\left[\log w_{\xi,\chi}(\zeta)-\log\xi\right].\label{ctx2}
\end{align}

\begin{figure}[!h]
\includegraphics[width=1\textwidth]{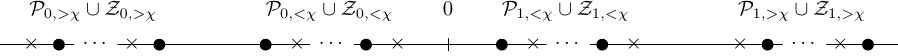}\\
\caption{Poles, represented with crosses, and zeros, represented with dots, of $\mathcal{G}_\chi$. Satisfied for Assumption \ref{ap65} where (\ref{sep1}) holds.} \label{figzp1} 
\end{figure}

\begin{figure}[!h]
\includegraphics[width=1\textwidth]{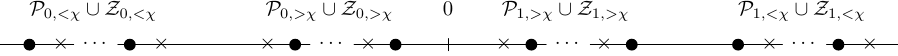}\\
\caption{Poles, represented with crosses, and zeros, represented with dots, of $\mathcal{G}_\chi$. Satisfied for Assumption \ref{ap65} where (\ref{sep2}) holds.} \label{figzp2}
\end{figure}

We now want to establish conditions under which (\ref{eeq}) has at most one pair of complex conjugate roots. For that we need to consider zeros and poles of $\mathcal{G}_{1,>\chi}(w)$, $\mathcal{G}_{1,<\chi}(w)$,
$\mathcal{G}_{0,>\chi}(w)$, and $\mathcal{G}_{0,<\chi}(w)$. Our goal will be to fully separate zeros and poles of each function from the zeros and poles of the others. By this we mean that all zeros and poles of of one function are either all to the left or all to the right of all zeros and poles of the other. We will further require that zeros and poles for each function alternate, i.e. that sorted from smallest to largest, we have that a zero is followed by a pole and vice versa. More precisely, we will look for conditions so that  $\mathcal{G}_\chi$ has poles and zeros positioned as in Figure \ref{figzp1} or as in Figure \ref{figzp2}. 

Conditions needed to satisfy the separating property and alternating property are stated in Assumption \ref{ap65}. The analysis in Remark \ref{remzp} justifies the conditions from Assumption \ref{ap65}.

\begin{assumption}\label{ap65}
Let $i,j\in [n]$ and $p_1,p_2 \in [m]$.
\begin{enumerate}
\item
\textbf{Separating condition.} For $b_{p_1,i}=-$, $b_{p_2,j}=+$ and $a_j=a_i$ holds that
\begin{align}\label{sep1}
\tau_{i}^{-1}\tau_{j}&
\begin{cases}\geq e^{V_{p_2}-V_{p_1-1}},&{if} \ p_1>p_2\\
>1,&{if} \ p_1=p_2
\end{cases}.
\end{align}

\item \textbf{Alternating condition.} For $b_{p_1,i}=b_{p_2,j}$, $a_i=a_j$, and $\tau_{i}> \tau_{j}$ holds that
\begin{align*}
\tau_{i}^{-1}\tau_{j}< e^{V_{p_2-1}-V_{p_1}}.
\end{align*}
\item Alternatively, we can assume that instead of (\ref{sep1}) it holds that
\begin{align}\label{sep2}
\tau_{i}^{-1}\tau_{j}< e^{V_{p_2-1}-V_{p_1}},\quad {if} \ p_1\geq p_2.
\end{align}
\end{enumerate}
\end{assumption}

\begin{remark} \label{uniformzp} Note that in the uniform case Assumption \ref{ap65} is satisfied if between any two transition points  all Ls have the same sign and all Rs have the same sign (not necessarily the same as Ls). More precisely, if for all $i,j\in [n]$ and $p \in [m]$,  $a_i=a_j$ implies $b_{p,i}=b_{p,j}$.

When we say Assumption \ref{ap65} holds, we mean either Assumption \ref{ap65}(1)(2) hold, or Assumption \ref{ap65}(2)(3) hold.
\end{remark}

In the remark below we justify Assumption \ref{ap65}.

\begin{remark} \label{remzp}Let $\mathrm{x}$ stand for one of the following 4 cases: $1,>\chi$; $1,<\chi$; $0,>\chi$; or $0,<\chi$. We denote the set of poles and zeros of $\mathcal{G}_\mathrm{x}$ with $\mathcal{P}_{\mathrm{x}}$ and $\mathcal{Z}_\mathrm{x}$, respectively. We have that $\mathcal{Z}_{\mathrm{x}}=\mathcal{N} _{\mathrm{x}}\setminus \mathcal{D}_{\mathrm{x}}$ and $\mathcal{P}_{\mathrm{x}}=\mathcal{D} _{\mathrm{x}}\setminus \mathcal{N}_{\mathrm{x}}$, where $\mathcal{D} _{\mathrm{x}}$, respectively $\mathcal{N}_{\mathrm{x}}$, denotes the set of points where the denominator, respectively numerator, of $\mathcal{G}_\mathrm{x}$ vanishes. Further,  $\mathcal{D}_{\mathrm{x}}= \bigcup_{j\in [n]}\mathcal{D}_{j,\mathrm{x}}$ and $\mathcal{N}_{\mathrm{x}}= \bigcup_{j\in [n]}\mathcal{N}_{j,\mathrm{x}}$, where

\begin{align*}
\mathcal{D}_{j,1,>\chi}&=\left\{e^{\max\{V_{p-1},\chi\}}\tau_j^{-1}|V_p>\chi, b_{p,j}=-, a_j=a_{i^*}\right\}\\
\mathcal{N}_{j,1,>\chi}&=\left\{e^{V_p}\tau_j^{-1}|V_p>\chi, b_{p,j}=-, a_j=a_{i^*}\right\}
\end{align*}

\begin{align*}
\mathcal{D}_{j,1,<\chi}&=\left\{e^{\min\{V_{p},\chi\}}\tau_j^{-1}|V_{p-1}<\chi, b_{p,j}=+, a_j=a_{i^*}\right\}\\
\mathcal{N}_{j,1,<\chi}&=\left\{e^{V_{p-1}}\tau_j^{-1}|V_{p-1}<\chi, b_{p,j}=+, a_j=a_{i^*}\right\}
\end{align*}

\begin{align*}
\mathcal{D}_{j,0,>\chi}&=\left\{-e^{V_p}\tau_j^{-1}|V_{p}>\chi, b_{p,j}=-, a_j \neq a_{i^*}\right\}\\
\mathcal{N}_{j,0,>\chi}&=\left\{-e^{\max \{V_{p-1},\chi\}}\tau_j^{-1}|V_{p}>\chi, b_{p,j}=-, a_j \neq a_{i^*}\right\}
\end{align*}

\begin{align*}
\mathcal{D}_{j,0,<\chi}&=\left\{-e^{V_{p-1}}\tau_j^{-1}|V_{p-1}<\chi, b_{p,j}=+, a_j \neq a_{i^*}\right\}\\
\mathcal{N}_{j,0,<\chi}&=\left\{-e^{\min \{V_{p},\chi\}}\tau_j^{-1}|V_{p-1}<\chi, b_{p,j}=+, a_j \neq a_{i^*}\right\}
\end{align*}

Note that zeros and poles of $\mathcal{G}_\mathrm{x}$ are positive when $\mathrm{x}$ is $1,>\chi$ or $1,<\chi$ and they are negative when $\mathrm{x}$ is $0,>\chi$ or $0,<\chi$. 

Observe that when sorted from smallest to largest, points $\mathcal{D}_{j,1,>\chi}\setminus \mathcal{N}_{j,1,>\chi}$ and $\mathcal{N}_{j,1,>\chi}\setminus \mathcal{D}_{j,1,>\chi}$ alternate, the smallest point is from $\mathcal{D}_{j,1,>\chi}\setminus \mathcal{N}_{j,1,>\chi}$ and the  largest from $\mathcal{N}_{j,1,>\chi}\setminus \mathcal{D}_{j,1,>\chi}$ with the rest of the points alternating between these two sets.  Now, to assure that zeros and poles of $\mathcal{G}_{1,>\chi}(w)$ alternate, we need to make sure that as we take the union over all $j$s the property is maintained. If $\tau_i=\tau_j$, it is not hard to check that points $\left(\mathcal{D}_{i,1,>\chi}\cup \mathcal{D}_{j,1,>\chi} \right)\setminus \left(\mathcal{N}_{i,1,>\chi} \cup \mathcal{N}_{j,1,>\chi}\right)$ and $\left(\mathcal{N}_{i,1,>\chi}\cup \mathcal{N}_{j,1,>\chi} \right)\setminus \left(\mathcal{D}_{i,1,>\chi} \cup \mathcal{D}_{j,1,>\chi}\right)$ alternate. In the case when $\tau_i>\tau_j$ we can choose to keep all points in 
$\mathcal{D}_{i,1,>\chi}\cup \mathcal{N}_{i,1,>\chi}$ to the left of  all points in $\mathcal{D}_{j,1,>\chi}\cup \mathcal{N}_{j,1,>\chi}$, which can be achieved if Assumption \ref{ap65} (2) holds (when we take into account the separating condition for other three cases as well).

To keep zeros and poles of $\mathcal{G}_{1,<\chi}(w)$ separated from zeros and poles of $\mathcal{G}_{1,>\chi}(w)$ we can choose
\begin{enumerate}
\item that all zeros and poles of $\mathcal{G}_{1,<\chi}(w)$ are to the left of all zeros and poles of $\mathcal{G}_{1,>\chi}$;
which can be guaranteed if Assumption \ref{ap65} (1) where (\ref{sep1}) holds; or alternatively 
\item that all zeros and poles of $\mathcal{G}_{1,<\chi}(w)$ are to the right of all zeros and poles of $\mathcal{G}_{1,>\chi}$; then Assumption \ref{ap65} (1) where (\ref{sep2}) holds can guarantee this case.
\end{enumerate}
To keep zeros and poles of $\mathcal{G}_{0,<\chi}(w)$ separated from zeros and poles of $\mathcal{G}_{0,>\chi}(w)$, we can make an analogous argument.
\end{remark}

\begin{lemma}\label{l65}Suppose that Assumption \ref{ap65} holds. Then for any $\zeta\in \RR$, the equation in $w$
\begin{align}
\mathcal{G}_{\chi}(w)=\zeta^n,\label{ceqn}
\end{align}
has at most one pair of complex conjugate roots.
\end{lemma}

\begin{proof} Let
\begin{align*}
\mathcal{G}_{\chi}(w)=C\frac{\prod_{b_i\in \mathcal{Z}_{1,>\chi}
    \cup\mathcal{Z}_{0,>\chi}\cup\mathcal{Z}_{1,<\chi}
    \cup\mathcal{Z}_{0,<\chi}}(w-b_i)}{\prod_{a_i\in \mathcal{P}_{1,>\chi}
    \cup\mathcal{P}_{0,>\chi}\cup\mathcal{P}_{1,<\chi}
    \cup\mathcal{P}_{0,<\chi}}(w-a_i)},
\end{align*}
where $\{a_1<a_2<\cdots<a_k\}$, $\{b_1<b_2<\cdots<b_k\}$, and $C\neq 0$ is an absolute constant. 

Observe that under Assumptions \ref{ap65} in both cases (i.e. if either (\ref{sep1}) or (\ref{sep2}) holds) poles can be divided in three segments: $\{a_1,\dots,a_{k_1}\}$,$\{a_{k_1+1},\dots,a_{k_2}\}$, and $\{a_{k_2+1},\dots,a_k\}$ so that for each pair of consecutive poles $a_i,a_{i+1}$ from one of the segments there is a unique zero $b_j\in (a_i,a_{i+1})$ where $j\in[k]$. For such a pair of poles $a_i,a_{i+1}$, it is straightforward to check that one of the following two cases occurs
    \begin{itemize}
        \item $\lim_{w\rightarrow a_i+}\mathcal{G}_{\chi}(w)=-\infty$ and $\lim_{w\rightarrow a_{i+1}-}\mathcal{G}_{\chi}(w)=+\infty$; or
        \item $\lim_{w\rightarrow a_i+}\mathcal{G}_{\chi}(w)=+\infty$ and $\lim_{w\rightarrow a_{i+1}-}\mathcal{G}_{\chi}(w)=-\infty$.
    \end{itemize}
By continuity $\mathcal{G}_{\chi}(w)$ is a surjection from $(a_i,a_{i+1})$ onto $(-\infty,\infty)$. Hence for each $\zeta \in \RR$, 
$\mathcal{G}_{\chi}(w)=\zeta^n$ has at least one root in $(a_i,a_{i+1})$. Since there are three segments, note that this will give us at least $k-3$ real roots of (\ref{ceqn}).

The following cases might occur
\begin{enumerate}
    \item $s^n=C$. In this case the equation (\ref{ceqn}) in $w$ has at most $k-1$ roots in the complex plane. But we already have $(k-3)$ real roots, hence in this case (\ref{ceqn}) has at most one pair of complex conjugate roots.
    \item $s^n\neq C$. In this case the equation (\ref{ceqn}) in $w$ has exactly $k$ roots in the complex plane. Again it is straightforward to check that one of the following two cases occurs
    \begin{itemize}
        \item $\lim_{w\rightarrow a_k+}\mathcal{G}_{\chi}(w)=-\infty$ and $\lim_{w\rightarrow a_{1}-}\mathcal{G}_{\chi}(w)=+\infty$; or
        \item $\lim_{w\rightarrow a_k+}\mathcal{G}_{\chi}(w)=+\infty$ and $\lim_{w\rightarrow a_{1}-}\mathcal{G}_{\chi}(w)=-\infty$.
    \end{itemize}
    By continuity $\mathcal{G}_{\chi}(w)$ is a surjection from $(-\infty,a_{1})\cup(a_p,\infty)$ onto $(-\infty,C)\cup(C,\infty)$. Hence for each $s\neq C$, 
$\mathcal{G}_{\chi}(w)=s$ has at least one root in $(-\infty,a_{1})\cup(a_p,\infty)$. Since (\ref{ceqn}) at least $(k-2)$ real roots, we deduce that it has at most one pair of complex conjugate roots.
\end{enumerate}
\end{proof}

\begin{lemma}\label{l66}Suppose Assumptions \ref{ap5} and \ref{ap65}(2)(3) hold. Let $\mathcal{H}(\chi,\kappa)$ be the limit of the rescaled height function of pure dimer coverings on rail yard graphs as $\epsilon\rightarrow 0$, as obtained in Theorem \ref{t61}. Let
\begin{align*}
\mathcal{R}_{\chi}:=\mathcal{P}_{1,>\chi}\cup \mathcal{P}_{0,>\chi}.
\end{align*}
Assume that the equation (\ref{ceqn}) in $w$ with $\zeta=e^{-\kappa}$ has exactly one pair of nonreal conjugate roots and $\mathcal{R}_{\chi}\neq \emptyset$.
\begin{enumerate}
    \item If $\mathcal{P}_{1,>\chi}=\emptyset$, then 
    \begin{align*}
    \mathcal{G}_{\chi}(0)<e^{-n\kappa}.
    \end{align*}
    \item If $\mathcal{P}_{0,>\chi}=\emptyset$, then 
    \begin{align*}
    \mathcal{G}_{\chi}(0)>e^{-n\kappa}.
    \end{align*}
\end{enumerate}
\end{lemma}

\begin{proof}We only prove part (1) here; part (2) can be proved using exactly the same technique. Assume $\mathcal{R}_{\chi}\neq \emptyset$ and $\mathcal{P}_{1,>\chi}=\emptyset$, then $\mathcal{P}_{0,>\chi}\neq \emptyset$. Let
\begin{align}
B_1:=\max\mathcal{P}_{0,>\chi}<0.\label{db1}
\end{align}
\begin{itemize}
\item If $\mathcal{P}_{1,<\chi}\neq  \emptyset$, let
\begin{align}
B_2:=
\min\mathcal{P}_{1,<\chi}>0.\label{db21}
\end{align}
\item If $\mathcal{P}_{1,<\chi}= \emptyset$, let
\begin{align}
B_2:=+\infty.\label{db22}
\end{align}
\end{itemize}
Then by (\ref{dgs1})-(\ref{dgs4}) and (\ref{dgc}) we have
\begin{align}
\lim_{u\rightarrow B_1+}\mathcal{G}_{\chi}(u)=-\infty.\label{lsb1}
\end{align}
Let $K$ be the total number of complex roots (counting multiplicities) in $w$ of $\mathcal{G}_{\chi}(w)=e^{-n\kappa}$. 
From the proof of Lemma \ref{l65}, we see that there are at least $K-2$ real roots of $\mathcal{G}_{\chi}(w)=e^{-n\kappa}$ in $(-\infty,B_1)\cup (B_2,+\infty)$. If $\mathcal{G}_{\chi}(0)\geq e^{-n\kappa}$, by (\ref{lsb1}) and the continuity of $\mathcal{G}_{\chi}(w)$ when $w\in (B_1,B_2)$, we deduce that there is at least one real root of $\mathcal{G}_{\chi}(w)=e^{-n\kappa}$ in $(B_1,0]$, which contradicts the assumption that $\mathcal{G}_{\chi}(w)=e^{-n\kappa}$ has exactly one pair of nonreal conjugate roots. Then part (1) of the lemma follows.
\end{proof}

Then we have the following proposition:
\begin{proposition}\label{p67}Suppose Assumption \ref{ap5} and \ref{ap65} hold. Let $\mathcal{H}(\chi,\kappa)$ be the limit of the rescaled height function of pure dimer coverings on rail yard graphs as $\epsilon\rightarrow 0$, as obtained in Theorem \ref{t61}. Assume that equation (\ref{ceqn}) in $w$ with $\zeta=e^{-\kappa}$ has exactly one pair of nonreal conjugate roots and $\mathcal{R}_{\chi}\neq \emptyset$,  $\mathcal{G}_{\chi}(0)\neq e^{-n\kappa}$ and $\mathcal{G}_{\chi}(\infty)\neq e^{-n\kappa}$. Let $\mathbf{w}_+$ be the unique nonreal root of $\mathcal{G}_{\chi}(w)=e^{-n\kappa}$ in the upper half plane, and the branch of $\mathrm{arg}(\cdot)$ be chosen such that $\mathrm{arg}(\mathbf{w}_{+})\in(0,\pi)$.
Then
\begin{enumerate}[label=(\alph*)]
\item if Assumption \ref{ap65}(2)(3) hold, 
\begin{align}
\frac{\partial \mathcal{H}(\chi,\kappa)}{\partial \kappa}\label{mde}
=2-\frac{2\mathrm{arg}(\mathbf{w}_{+})}{\pi};
\end{align}
\item if Assumption \ref{ap65}(1)(2) hold, 
\begin{align}
\frac{\partial \mathcal{H}(\chi,\kappa)}{\partial \kappa}\label{mde1}
=\begin{cases}
\frac{2\mathrm{arg}(\mathbf{w}_{+})}{\pi}&\mathrm{if}\ \mathbf{1}_{\mathcal{G}_{\chi}(0)> e^{-n\kappa}}+\mathbf{1}_{\mathcal{G}_{\chi}(\infty)> e^{-n\kappa}}\in\{0,2\};\\
2\left[1-\frac{\mathrm{arg}(\mathbf{w}_{+})}{\pi}\right]&\mathrm{otherwise.}
\end{cases}
\end{align}
\end{enumerate}

\end{proposition}

\begin{proof}We first prove (a) here.
By (\ref{slht}) (\ref{dsm2}), we obtain
\begin{align*}
\frac{\partial \mathcal{H}(\chi,\kappa)}{\partial \kappa}
=-\left.\lim_{\epsilon\rightarrow 0+}\frac{1}{\pi}\Im \Theta_{\chi}(\zeta+\mathbf{i}\epsilon)\right|_{\zeta=e^{-\kappa}}
\end{align*}
By (\ref{ctx}), (\ref{ctx1}), (\ref{ctx2}), we obtain
\begin{align}
\frac{\partial \mathcal{H}(\chi,\kappa)}{\partial \kappa}
=2\mathbf{1}_{\mathcal{G}_{\chi}(0)>e^{-n\kappa}}-\frac{2}{\pi}\sum_{\xi\in \mathcal{R}_{\chi}}\left[\mathrm{arg}(w_{\xi,\chi}(e^{-\kappa}))-\mathrm{arg}(\xi)\right]\label{phk}
\end{align}
where the branch of $\mathrm{arg}$ is chosen to have range $(-\pi, \pi]$. Hence we have
\begin{align*}
\mathrm{arg}(\xi)=\begin{cases}0&\mathrm{if}\ \xi>0;\\\pi&\mathrm{otherwise}.\end{cases}
\end{align*}
Under the assumption that $\mathcal{R}_{\chi}\neq \emptyset$ and $\mathcal{G}_{\chi}(0)\neq e^{-n\kappa}$, the following cases might occur:
\begin{enumerate}
\item 
$\mathcal{P}_{0,>\chi}\neq \emptyset$. In this case
the number of negative poles in $\mathcal{R}_{\chi}$ is exactly $|\mathcal{P}_{0,>\chi}|$. From the proof of Lemma \ref{l65} we see that  there are at least $|\mathcal{P}_{0,>\chi}|-1$ negative real roots in $\{w_{\xi,\chi}(e^{-\kappa})\}_{\xi\in \mathcal{R}_{\chi}}$.
Let $B_1$ be defined as in (\ref{db1}).
\begin{enumerate}
\item If $\mathcal{P}_{1,>\chi}\neq \emptyset$, let
\begin{align*}
B_2:=
\min\mathcal{P}_{1,>\chi}>0.
\end{align*}
\item If $\mathcal{P}_{1,>\chi}= \emptyset$ and If $\mathcal{P}_{1,<\chi}\neq  \emptyset$, let $B_2$ be defined as in (\ref{db21}).
\item If $\mathcal{P}_{1,>\chi}= \emptyset$ and If $\mathcal{P}_{1,<\chi}= \emptyset$, let $B_2$ be defined as in (\ref{db22}).
\end{enumerate}
Then by (\ref{dgs1})-(\ref{dgs4}) and (\ref{dgc}) we have (\ref{lsb1})
The following cases might occur
\begin{enumerate}
    \item $\mathcal{G}_{\chi}(0)>e^{-n\kappa}$: then there exists a unique root in $\{w_{\xi,\chi}(e^{-\kappa})\}_{\xi\in \mathcal{R}_{\chi}}\cap (B_1,B_2)$ which is negative; in this case the argument of each negative pole in $\mathcal{R}_{\chi}$ cancels with an argument of a unique negative root in $\{w_{\xi,\chi}(e^{-\kappa})\}_{\xi\in\mathcal{R}_{\chi}}$;
    \item $\mathcal{G}_{\chi}(0)<e^{-n\kappa}$: then there exists no root in $\{w_{\xi,\chi}(e^{-\kappa})\}_{\xi\in \mathcal{R}_{\chi}}\cap (B_1,0)$; in this case there is a unique negative pole in $\mathcal{R}_{\chi}$ whose argument cannot  cancel with an argument of a unique negative root in $\{w_{\xi,\chi}(e^{-\kappa})\}_{\xi\in\mathcal{R}_{\chi}}.$
\end{enumerate}
In either case we have (\ref{mde}) holds.
\item $\mathcal{P}_{1,>\chi}
\neq\emptyset$ and 
$\mathcal{P}_{0,>\chi}= \emptyset$. In this case there are neither negative poles in $\mathcal{R}_{\chi}$ nor negative real roots in $\{w_{\xi,\chi}(e^{-\kappa})\}_{\xi\in \mathcal{R}_{\chi}}$.
By Lemma \ref{l66}, we have $\mathcal{G}_{\chi}(0)>0$. Then (\ref{mde}) follows from (\ref{phk}).
\end{enumerate}

Now we prove (b). Following the same argument as the proof of part(a), we obtain
\begin{align*}
\frac{\partial \mathcal{H}(\chi,\kappa)}{\partial \kappa}
=
2\left[\mathbf{1}_{\mathcal{G}_{\chi}(0)> e^{-n\kappa}}+\mathbf{1}_{\mathcal{G}_{\chi}(\infty)> e^{-n\kappa}}-\frac{\mathrm{arg}(\mathbf{w})}{\pi}\right]
\end{align*}
where $\mathbf{w}$ is some non-real root of $\mathcal{G}_{\chi}(w)=e^{-n\kappa}$. Depending on the value of 
$\mathbf{1}_{\mathcal{G}_{\chi}(0)> e^{-n\kappa}}+\mathbf{1}_{\mathcal{G}_{\chi}(\infty)> e^{-n\kappa}}$, $\mathbf{w}$ may be chosen to be the root in the upper half plane or the lower half plane with appropriate branch of $arg(\cdot)$ such that $\frac{\partial\mathcal{H}}{\partial\kappa}\in[0,2]$, as in (\ref{drh}). Then (\ref{mde1}) follows.
\end{proof}

As we shall see in section \ref{sect:ex}, for pyramid partitions $\mathcal{G}_{\chi}(0)=\mathcal{G}_{\chi}(\infty)=1$.

\begin{definition}Let $\{RYG(l^{(\epsilon)},r^{(\epsilon)},\underline{a}^{(\epsilon)},\underline{b}^{(\epsilon)})\}_{\epsilon>0}$ be a collection of rail-yard graphs satisfying Assumptions \ref{ap5} and \ref{ap65}. Let $\mathcal{H}(\chi,\kappa)$ be the limit of the rescaled height function of pure dimer coverings on rail yard graphs as $\epsilon\rightarrow 0$, as obtained in Theorem \ref{t61}. The liquid region for the limit shape of pure dimer coverings on these rail yard graphs as $\epsilon\rightarrow 0$ is defined to be
\begin{align*}
\mathcal{L}:=\left\{(\chi,\kappa)\in(l^{(0)},r^{(0)})\times \RR:\frac{\partial \mathcal{H}}{\partial \kappa}(\chi,\kappa)\in (0,2)\right\}
\end{align*}
and the frozen region is defined to be 
\begin{align*}
\left\{(\chi,\kappa)\in(l^{(0)},r^{(0)})\times \RR:\frac{\partial \mathcal{H}}{\partial \kappa}(\chi,\kappa)\in \{0,2\}\right\}.
\end{align*}
The frozen boundary is defined to be the boundary separating the frozen region and the liquid region.
\end{definition}

\begin{remark}By Proposition \ref{p67}, we see that if $\mathcal{R}_{\chi}\neq \emptyset$ and $\mathcal{G}_{\chi}(0)\neq \emptyset$, $(\chi,\kappa)\in(l^{(0)},r^{(0)})\times \RR$ is in the liquid region if and only if the following equation
\begin{align}
\mathcal{G}_{\chi}(w)=e^{-n\kappa}\label{gck}
\end{align}
in $w$ has exactly one pair of nonreal conjugate roots in $w$. By Lemma \ref{l65}, we see that the frozen boundary is given by the condition that (\ref{gck}) has double real roots.
\end{remark}

Next we shall find the frozen boundary. The discussion above shows that if $\mathcal{R}_{\chi}\neq \emptyset$ and $\mathcal{G}_{\chi}(0)\neq \emptyset$, $(\chi,\kappa)\in(l^{(0)},r^{(0)})\times \RR$ is on the frozen boundary if and only if $(\chi,\kappa)$ satisfies the following system of equations
\begin{align}\label{fb}
\begin{cases}
\mathcal{G}_{\chi}(w)=e^{-n\kappa}\\
\frac{d\log \mathcal{G}_{\chi}(w)}{dw}=0
\end{cases}.
\end{align}
The second equation in (\ref{fb}) gives
\begin{align}
0&=\sum_{\substack{p\in [m],V_p>\chi\\j\in[n]: b_{p,j}=-,a_j=a_{i^*}}}
\frac{1}{w-e^{V_p}\tau_j^{-1}}-\frac{1}{w-e^{\max\{V_{p-1},\chi\}}\tau_{j}^{-1}}\label{fb1}\\
&+\sum_{\substack{p\in[m],V_{p-1}<\chi\\j\in[n]:b_{p,j}=+,a_j=a_{i^*}}}
\frac{1}{w- e^{V_{p-1}}\tau_j^{-1}}-\frac{1}{w-e^{\min\{V_p,\chi\}}\tau_{j}^{-1}}\notag\\
&+\sum_{\substack{p\in[m],V_p>\chi\\j\in[n]:b_{p,j}=-,a_j \neq a_{i^*}}}\frac{1}{w+e^{\max\{V_{p-1},\chi\}}\tau_j^{-1}}
-\frac{1}{w+e^{V_p}\tau_{j}^{-1}}\notag
\\
&+\sum_{\substack{p\in[m],V_{p-1}<\chi\\j\in[n]:b_{p,j}=+,a_j \neq a_{i^*}}}
\frac{1} {w+e^{\min\{V_p,\chi\}} \tau_{j}^{-1}}-\frac{1}{w+e^{V_{p-1}}\tau_j^{-1}}.\notag
\end{align}

\section{Height Fluctuations and Gaussian Free Field}\label{sect:gff}

In this section, we prove that the fluctuations of height function converge to the pull-back Gaussian free field (GFF) in the upper half plane under a diffeomorphism from the liquid region to the upper half plane. The main theorem proved in this section is Theorem \ref{t77}.

 \subsection{Gaussian free field} 
 
 Let $C_0^{\infty}$ be the space of smooth real-valued functions with compact
 support in the upper half plane $\HH$. The \emph{Gaussian free field} (GFF)
 $\Xi$ on $\HH$ with the zero boundary condition is a collection of Gaussian
 random variables $\{\Xi_{f}\}_{f\in C_0^{\infty}}$ indexed by functions in
 $C_0^{\infty}$, such that the covariance of two Gaussian random variables
 $\Xi_{f_1}$, $\Xi_{f_2}$ is given by
 \begin{equation}
 \mathrm{Cov}(\Xi_{f_1},\Xi_{f_2})=\int_{\HH}\int_{\HH}f_1(z)f_2(w)G_{\HH}(z,w)dzd\ol{z}dwd\ol{w},\label{cvf}
 \end{equation}
 where
 \begin{equation*}
 G_{\HH}(z,w):=-\frac{1}{2\pi}\ln\left|\frac{z-w}{z-\ol{w}}\right|,\qquad z,w\in \HH
 \end{equation*}
 is the Green's function of the Laplacian operator on $\HH$ with the Dirichlet
 boundary condition. The Gaussian free field $\Xi$ can also be considered as a
 random distribution on $C_0^{\infty}$,
 such that for any $f\in C_0^{\infty}$, we have
 \begin{equation*}
 \Xi(f)=\int_{\HH}f(z)\Xi(z)dz:=\xi_f.
 \end{equation*}
Here $\Xi(f)$ is the Gaussian random variable with respect to $f$, which has
mean 0 and variance given by \eqref{cvf} with $f_1$ and $f_2$ replaced by $f$.
See~\cite{SS07} for more about the GFF.

\subsection{$\mathbf{w}_{+}$\label{ss72} as a mapping from $\mathcal{L}$ to $\mathbb{H}$}

By (\ref{dgs1})-(\ref{dgs4}) and (\ref{dgc}), we may write $\mathcal{G}_{\chi}(w)$ as the quotient of two functions $U_{\chi}(w)$ and $R(w)$, such that $U_{\chi}(w)$ depends on $\chi$ and $R(w)$ is independent of $\chi$. More precisely,
\begin{align*}
\mathcal{G}_{\chi}(w)=\frac{U_{\chi}(w)}{R(w)};
\end{align*}
where
\begin{align*}
U_{\chi}(w)=e^{\phi(\chi)}\frac{\prod_{j\in[n],a_j=R}\left(1+e^{-\chi}w \tau_j\right)}{\prod_{j\in[n],a_j=L}({1-e^{-\chi}w \tau_j})}
\end{align*}
and
\begin{small}
\begin{align*}
    \phi(\chi)&=\sum_{p\in[m],V_{p-1}>\chi} \sum_{\substack{j\in[n]\\(a_j,b_{p,j})=(L,-)}}(V_p-V_{p-1})
   +\sum_{p\in[m],V_{p-1}<\chi<V_p} \sum_{\substack{j\in[n]\\(a_j,b_{p,j})=(L,-)}}(V_p-\chi)\\
   & -\sum_{p\in[m],V_{p-1}>\chi} \sum_{\substack{j\in[n]\\(a_j,b_{p,j})=(R,-)}}(V_p-V_{p-1})
   -\sum_{p\in[m],V_{p-1}<\chi<V_p} \sum_{\substack{j\in[n]\\(a_j,b_{p,j})=(R,-)}}(V_p-\chi)
\end{align*}
\end{small}
and
\begin{align}
R(w)=A_1\cdot A_2\label{drw}
\end{align}
\begin{align*}
A_1&=\prod_{j=1}^{n}\prod_{\substack{p\in [m],V_p>\chi\\b_{p,j}=-,a_j=a_{i^*}}}
\left(1-e^{-V_p}w \tau_{j}\right)^{-1}\prod_{\substack{p\in [m],V_{p-1}>\chi\\b_{p,j}=-,a_j=a_{i^*}}}\left(1- e^{-V_{p-1}}w\tau_j\right)\\
&\times\prod_{\substack{p\in [m],V_{p-1}<\chi\\b_{p,j}=+,a_j=a_{i^*}}}
\left(1- e^{-V_{p-1}}w\tau_j\right)^{ 
-1}
\prod_{\substack{p\in [m],V_p<\chi\\b_{p,j}=+,a_j=a_{i^*}}}\left(1-e^{-V_p}w \tau_{j}\right),
\end{align*}

\begin{align*}
A_2&=\prod_{j=1}^{n}\prod_{\substack{p\in [m],V_p>\chi\\b_{p,j}=-,a_j\neq a_{i^*}}}
\left(1+e^{-V_p}w\tau_{j}\right)\prod_{\substack{p\in [m],V_{p-1}>\chi\\b_{p,j}=-,a_j\neq a_{i^*}}}\left({1+ e^{-V_{p-1}}w\tau_j}\right)^{-1}\\
&\times\prod_{\substack{p\in [m],V_p<\chi\\b_{p,j}=+,a_j\neq a_{i^*}}}
\left(1+e^{-V_p}w \tau_{j}\right)^{-1}
\prod_{\substack{p\in [m],V_{p-1}<\chi\\b_{p,j}=+,a_j\neq a_{i^*}}}\left({1+ e^{-V_{p-1}}w\tau_j}\right).
\end{align*}
Observe that $A_1$ and $A_2$ do not depend on $\chi$; it because the first and third terms combined do not depend on $\chi$, as well as the second and fourth combined. Now assume that $a_{i^*}=L$. Then 
\begin{small}
\begin{align}
A_1&=\prod_{j\in[n]:a_j=L}\left[\left(1-e^{-V_0}w\tau_j\right)^{-\mathbf{1}_{b_{1,j}=+}}\left(1-e^{-V_m}w\tau_j\right)^{-\mathbf{1}_{b_{m,j}=-}}\right]\notag\\
&\times\left[\prod_{p=1}^{m-1}\left(1-e^{-V_p}w\tau_j \right)^{-\mathbf{1}_{b_{p+1,j}=+}+\mathbf{1}_{b_{p,j}=+}}\right]\label{da1}
\end{align}
\end{small}
and
\begin{small}
\begin{align}
A_2&=\prod_{j\in[n]:a_j=R}\left[\left(1+e^{-V_0}w\tau_j\right)^{\mathbf{1}_{b_{1,j}=+}}\left(1+e^{-V_m}w\tau_j\right)^{\mathbf{1}_{b_{m,j}=-}}\right]\notag\\
&\times\left[\prod_{p=1}^{m-1}\left(1+e^{-V_p}w\tau_j \right)^{\mathbf{1}_{b_{p+1,j}=+}-\mathbf{1}_{b_{p,j}=+}}\right].\label{da2}
\end{align}
\end{small}

We shall always use $[\cdot]^{\frac{1}{n}}$ to denote the branch which takes positive real values on the positive real line. 
Define
\begin{align*}
U(w_*,z_*)=\left[\frac{\prod_{j\in[n],a_j=R}\left(1+w_* \tau_j\right)}{\prod_{j\in[n],a_j=L}({1-w_*\tau_j})}\right]^{\frac{1}{n}}-z_*.
\end{align*}

Hence we have 
\begin{align}
\mathcal{G}_{\chi}(w)=e^{-n\kappa}\label{cpeq}
\end{align}
if and only if 
\begin{align*}
\begin{cases}
U(w_*,z_*)=0\\
[R(w)]^{\frac{1}{n}}=z\\
(w_*,z_*)=(e^{-\chi}w,e^{-\kappa}z)
\end{cases}.
\end{align*}

\begin{lemma}\label{l71}
Let
\begin{align*}
n_R:=|\{j\in[n]:a_j=R\}|;\qquad n_L:=|\{j\in[n]:a_j=L\}|;
\end{align*}
For any $(\alpha,\theta)$ such that 
\begin{align}
\theta\in(0,\pi),\qquad \alpha\in \left(0,\frac{n_L\pi+(n_R-n_L)\theta}{n}\right),\label{cta}
\end{align}
there exists a unique pair $(w_*,z_*)$ such that $\arg z_*=\alpha$ for some $k\in\ZZ$, $\arg w_*=\theta$ and $U(w_*,z_*)=0$. 
\end{lemma}

\begin{proof}
Note that $n_R+n_L=n$. For $\theta\in (0,\pi)$, define a map $B_{\theta}: [0,\infty)\rightarrow \RR$ by
\begin{align*}
B_{\theta}(\rho):=\frac{1}{n}\left[\sum_{j\in[n]:a_j=R}\arg(1+\rho e^{\mathbf{i}\theta}\tau_j)
-\sum_{j\in[n]:a_j=L}\arg(1-\rho e^{\mathbf{i}\theta}\tau_j)\right],
\end{align*}
where the branch of $\arg(\cdot)$ is chosen such that it has range $(-\pi,\pi]$. It is straightforward to check that $B_{\theta}(\rho)$ is strictly increasing when $\rho\in(0,\infty)$. Moreover,
\begin{align*}
\lim_{\rho\rightarrow0} B_{\theta}(\rho)=0;\qquad
\lim_{\rho\rightarrow\infty} B_{\theta}(\rho)=\frac{n_L\pi}{n}+\frac{(n_R-n_L)\theta}{n}.
\end{align*}
Since $B_{\theta}$ is a bijection from $(0,\infty)$ to $\left(0,\frac{n_L\pi}{n}+\frac{(n_R-n_L)\theta}{n}\right)$, for any $(\alpha,\theta)$ satisfying (\ref{cta}), we can find a unique $\rho>0$, such that $B_{\theta}(\rho)=\alpha$. Let
\begin{align*}
w_*:=\rho e^{\mathbf{i}\theta} \quad \textrm{and} \quad
z_*:=\left[\frac{\prod_{j\in[n],a_j=R}\left(1+w_* \tau_j\right)}{\prod_{j\in[n],a_j=L}({1-w_*\tau_j})}\right]^{\frac{1}{n}}.
\end{align*}
Then the lemma follows.
\end{proof}

\begin{proposition}For each $(\chi,\kappa)\in \mathcal{L}$, let $\mathbf{w}_+(\chi,\kappa)$ be the unique root of (\ref{cpeq}) in the upper half plane $\HH$. Then $\mathbf{w}_+:\mathcal{L}\rightarrow \HH$ is a diffeomorphism.
\end{proposition}

\begin{proof}We first show that $\mathbf{w}_+$ is a bijection. 
For any $w\in\mathbb{H}$, let $z=[R(w)]^{\frac{1}{n}}$. Let
\begin{align*}
\theta:=\arg w\in (0,\pi);\qquad \alpha:=\arg z.
\end{align*}
By (\ref{drw}) (\ref{da1}) (\ref{da2}) we obtain
\begin{align}
\alpha=\frac{1}{n}&\left[\sum_{j\in[n]:a_j=L}\begin{array}{l}-\arg\left(1-e^{-V_0}w\tau_j\right)\mathbf{1}_{b_{1,j}=+}-\arg\left(1-e^{-V_m}w\tau_j\right)\mathbf{1}_{b_{m,j}=-}\\+\sum_{p=1}^{m-1}\arg\left(1-e^{-V_p}\tau_j w\right)\left(-\mathbf{1}_{b_{p+1,j}=+}+\mathbf{1}_{b_{p,j}=+}\right)\end{array}\right.\notag\\
&\left.+\sum_{j\in[n]:a_j=R}\begin{array}{l}\arg\left(1+e^{-V_0}w\tau_j\right)
\mathbf{1}_{b_{1,j}=+}
+\arg\left(1+e^{-V_m}w\tau_j\right)\mathbf{1}_{b_{m,j}=-}\\+\sum_{p=1}^{m-1}\arg\left(1+e^{-V_p}w\tau_j \right)\left(\mathbf{1}_{b_{p+1,j}=+}-\mathbf{1}_{b_{p,j}=+}\right)\end{array}\right].\label{agr}
\end{align}

Then we have 
\begin{align}
\alpha=\frac{1}{n}&\left[
\sum_{j\in[n]:a_j=L}\begin{array}{l}\sum_{p=1}^m\left[\arg\left(1-e^{-V_{p}}w\tau_j\right)-\arg\left(1-e^{-V_{p-1}}w\tau_j\right)\right]\mathbf{1}_{b_{p,j}=+}
\\-\arg\left(1-e^{-V_m}w\tau_j\right)\end{array}\right.\notag\\
&\left.+\sum_{j\in[n]:a_j=R}\begin{array}{l}\sum_{p=1}^{m}\left[\arg\left(1+e^{-V_{p-1}}w\tau_j\right)-
\arg\left(1+e^{-V_{p}}w\tau_j\right)
\right]\mathbf{1}_{b_{p,j}=+}\\+\arg\left(1+e^{-V_m}w\tau_j\right)\end{array}\right]\label{al1}
\end{align}
and
\begin{align}
\alpha=\frac{1}{n}&\left[\sum_{j\in[n]:a_j=L}\begin{array}{l}\sum_{p=1}^{m}\left[
\arg\left(1-e^{-V_{p-1}}w\tau_j\right)
-\arg\left(1-e^{-V_p}w\tau_j\right)\right]\mathbf{1}_{b_{p,j}=-}\\
-\arg\left(1-e^{-V_0}w\tau_j\right)\end{array}\right.\notag\\
&+\left.\sum_{j\in[n]:a_j=R}
\begin{array}{l}
\sum_{p=1}^{m}\left[\arg\left(1+e^{-V_{p}}w\tau_j\right)-
\arg\left(1+e^{-V_{p-1}}w\tau_j\right)
\right]\mathbf{1}_{b_{p,j}=-}\\
+\arg\left(1+e^{-V_0}w\tau_j\right)
\end{array}
\right].\label{al2}
\end{align}
Note that for any $w\in \HH$, $u,v\in[0,\infty]$ and $u<v$, we have
\begin{align*}
-(\pi-\arg w)\leq \arg\left(1-u^{-1}w\right)< \arg\left(1-v^{-1}w\right)\leq 0
\end{align*}
and
\begin{align*}
\arg w\geq\arg\left(1+u^{-1}w\right)> \arg\left(1+v^{-1}w\right)\geq 0.
\end{align*}
Hence from (\ref{al1}), we obtain
\begin{align*}
\alpha&>\frac{1}{n}\left\{\sum_{j\in[n]:a_j=L}
\left[-\arg\left(1-e^{-V_m}w\tau_j\right)\right]+\sum_{j\in[n]:a_j=R}\left[\arg\left(1+e^{-V_m}w\tau_j\right)\right]\right\}\geq 0
\end{align*}
By (\ref{al2}), we obtain
\begin{align*}
\alpha&<\frac{1}{n}\left\{\sum_{j\in[n]:a_j=L}
\left[-\arg\left(1-e^{-V_0}w\tau_j\right)\right]+\sum_{j\in[n]:a_j=R}\left[\arg\left(1+e^{-V_0}w\tau_j\right)\right]\right\}\\
&\leq \frac{n_L\pi+(n_R-n_L)\theta}{n}.
\end{align*}
By Lemma \ref{l71}, we can a unique pair $(w_*,z_*)$ such that $\arg z_*=\alpha$ for some $k\in\ZZ$, $\arg w_*=\theta$ and $U(w_*,z_*)=0$. 
\begin{align*}
\chi:=\log\left(\frac{w}{w_*}\right),\qquad
\kappa:=\log\left(\frac{z}{z_*}\right),
\end{align*}
where the branch of the $\log(\cdot)$ is chosen such that it takes real values on the positive real axis. Then we deduce that $\mathbf{w}_+$ is a bijection. From the process we see that both the mapping $\mathbf{w}_+$ and its inverse are differentiable. Then the proposition follows.
\end{proof}

\subsection{Convergence of height fluctuations to GFF}

Splitting the sum of the RHS of (\ref{fb1}) into those depending on $\chi$ and those independent of $\chi$, we obtain
\begin{small}
\begin{align*}
0&=-\sum_{j\in[n]:a_j=L}\frac{1}{w-e^{\chi}\tau_j^{-1}}
+\sum_{j\in[n]:a_j=R}\frac{1}{w+e^{\chi}\tau_j^{-1}}\\
&+\sum_{\substack{p\in[m-1],j\in [n]\\a_j=L}}\frac{\mathbf{1}_{b_{p+1,j}=+}-\mathbf{1}_{b_{p,j}=+}}{w-e^{V_p}\tau_j^{-1}}
+\sum_{\substack{j\in[n]\\a_j=L}}\left(\frac{\mathbf{1}_{b_{m,j}=-}}{w-e^{V_m}\tau_j^{-1}}
+\frac{\mathbf{1}_{b_{1,j}=+}}{w-e^{V_0}\tau_j^{-1}}
\right)\\
&+\sum_{\substack{p\in[m-1],j\in [n]\\a_j=R}}\frac{\mathbf{1}_{b_{p,j}=+}-\mathbf{1}_{b_{p+1,j}=+}}{w+e^{V_p}\tau_j^{-1}}
-\sum_{\substack{j\in[n]\\a_j=R}}\left(\frac{\mathbf{1}_{b_{m,j}=-}}{w+e^{V_m}\tau_j^{-1}}
+\frac{\mathbf{1}_{b_{1,j}=+}}{w+e^{V_0}\tau_j^{-1}}
\right).
\end{align*}
\end{small}

Let $\mathcal{S}$ be the set of all the zeros and poles of $\mathcal{G}_{\chi}$ that are independent of $\chi$; or equivalently, $\mathcal{S}$ is the set of all the zeros and poles of $R(w)$. More precisely,
\begin{small}
\begin{align*}
\mathcal{S}&=\left\{e^{V_p}\tau_j^{-1}:p\in[m-1],j\in[n],a_j=L,b_{p,j}\neq b_{p+1,j}\right\}
\\
&\cup
\left\{-e^{V_p}\tau_j^{-1}:p\in[m-1],j\in[n],a_j=R,b_{p,j}\neq b_{p+1,j}\right\}\\
&\cup\{e^{V_0}\tau_j^{-1},e^{V_m}\tau_j^{-1}:j\in[n],a_j=L\}
\cup\{-e^{V_0}\tau_j^{-1},-e^{V_m}\tau_j^{-1}:j\in[n],a_j=R\}.
\end{align*}
\end{small}
Then we have the following lemma.

\begin{lemma}
Each $u\in \RR\setminus \mathcal{S}$ is a double root of (\ref{cpeq}) for a unique pair of $(\chi,\kappa)\in\RR^2$.
\end{lemma}

\begin{proof}Define
\begin{align*}
f(s)=\sum_{j\in[n]:a_j=L}\frac{1}{1-\tau_j^{-1}s}
-\sum_{j\in[n]:a_j=R}\frac{1}{1+\tau_j^{-1}s}
\end{align*}
and
\begin{small}
\begin{align*}
g(w)&=\sum_{\substack{p\in[m-1],j\in [n]\\a_j=L}}\frac{\mathbf{1}_{b_{p+1,j}=+}-\mathbf{1}_{b_{p,j}=+}}{w-e^{V_p}\tau_j^{-1}}
+\sum_{\substack{j\in[n]\\a_j=L}}\left(\frac{\mathbf{1}_{b_{m,j}=-}}{w-e^{V_m}\tau_j^{-1}}
+\frac{\mathbf{1}_{b_{1,j}=+}}{w-e^{V_0}\tau_j^{-1}}
\right)\\
&+\sum_{\substack{p\in[m-1],j\in [n]\\a_j=R}}\frac{\mathbf{1}_{b_{p,j}=+}-\mathbf{1}_{b_{p+1,j}=+}}{w+e^{V_p}\tau_j^{-1}}
-\sum_{\substack{j\in[n]\\a_j=R}}\left(\frac{\mathbf{1}_{b_{m,j}=-}}{w+e^{V_m}\tau_j^{-1}}
+\frac{\mathbf{1}_{b_{1,j}=+}}{w+e^{V_0}\tau_j^{-1}}
\right).
\end{align*}
\end{small}
Then $u$ is a double root for (\ref{ceqn}) for some $(\chi,\kappa)\in[r^{(0)},l^{(0)}]\times \RR$ if and only if 
\begin{align}
&e^{\kappa}=\left[\frac{R(u)}{U_{\chi}(u)}\right]^{\frac{1}{n}}\label{drt1},\\
&f(e^{\chi}u^{-1})=ug(u),\label{drt2}
\end{align}
where $[\cdot]^{\frac{1}{n}}$ is the branch that takes positive real value on the positive real axis.
The function $f(s)$ is defined in $\RR\setminus[\{-\tau_j\}_{j\in[n]:a_j=R}\cup\{\tau_j\}_{j\in[n],a_j=L}]$.
Suppose that we enumerate all the points in $\{-\tau_j\}_{j\in[n]:a_j=R}\cup\{\tau_j\}_{j\in[n],a_j=L}$ in the increasing order as follows:
\begin{align*}
-d_{n_L}<-d_{n_L-1}<\ldots<-d_1<0<\alpha_1<\alpha_2<\ldots<\alpha_{n_R}
\end{align*}
Since for all $s\in\RR\setminus[\{-\tau_j\}_{j\in[n]:a_j=R}\cup\{\tau_j\}_{j\in[n],a_j=L}]$,
\begin{align*}
f'(s)=\sum_{j\in[n]:a_j=L}\frac{1}{\tau_j(1-\tau_j^{-1}s)^2}
+\sum_{j\in[n]:a_j=R}\frac{1}{\tau_j(1+\tau_j^{-1}s)^2}>0;
\end{align*}
we obtain 
\begin{enumerate}
    \item $f$ is strictly increasing in each interval $(\alpha_i,\alpha_{i+1})$, for $i\in[n_R-1]$ from $-\infty$ to $\infty$;
    \item $f$ is strictly increasing in each interval $(-d_{j+1},-d_{j})$, for $j\in[n_L-1]$ from $-\infty$ to $\infty$;
    \item $f$ is strictly increasing in the interval $(-d_1,\alpha_1)$ from $-\infty$ to $\infty$;
    \item $f$ is strictly increasing in the interval $(\alpha_{n_R},\infty)$ from $-\infty$ to $0$;
    \item $f$ is strictly increasing in the interval $(-\infty,\alpha_{n_L})$ from $0$ to $\infty$.
\end{enumerate}

Hence for each $u\in \RR$ and for each set 
\begin{align}
\Delta\in &\{(-d_{n_L},-d_{n_L-1}),\ldots,(-d_2,d_1),(d_1,\alpha_1),(\alpha_1,\alpha_2),\ldots,\label{dta1}\\
&(\alpha_{n_R-1},\alpha_{n_R}),(\alpha_{n_R},\infty)\cup(-\infty,-d_{n_L})\},\notag
\end{align}
there is a unique $\chi$ such that (\ref{drt2}) holds and $e^{\chi}u^{-1}\in \Delta$.

For $j\in[n]$, let
\begin{align}
p_{j,L}&=\max\{p\in[0..m]:e^{V_p}\tau_j^{-1}<u, a_j=L\};\label{dpl}\\
p_{j,R}&=\max\{p\in[0..m]:u<-e^{V_p}\tau_j^{-1}, a_j=R\};\label{dpr}
\end{align}
again we take the convention that the minimum (resp.\ maximum) of an empty set is $\infty$ ($-\infty$); and assume for all $j\in[n]$
\begin{align*}
b_{-\infty,j}=-;\qquad
b_{m+1,j}=+.
\end{align*}
From (\ref{agr}), we obtain
\begin{small}
\begin{align}
\lim_{\epsilon\rightarrow 0+}\arg [R(u+\mathbf{i}\epsilon)]^{\frac{1}{n}}=\frac{\pi}{n}\Bigl[\sum_{\substack{j\in[n]\\a_j=R}}\mathbf{1}_{b_{p_{j,R}+1,j}=+}+
\sum_{\substack{j\in[n]\\a_j=L}}\mathbf{1}_{b_{p_{j,L}+1,j}=+}\Bigr],\label{agm5}
\end{align}
\end{small}
where $k\in\ZZ$. Moreover,
\begin{small}
\begin{align}
\lim_{\epsilon\rightarrow 0+}\arg [U_{\chi}(u+\mathbf{i}\epsilon)]^{\frac{1}{n}}&=\lim_{\epsilon\rightarrow 0+}\frac{1}{n}\Bigl[\sum_{\substack{j\in[n]\\a_j=R}}\arg(1+e^{-\chi}(u+\mathbf{i}\epsilon)\tau_j)-\sum_{\substack{j\in[n]\\a_j=L}}\arg(1-e^{-\chi}(u+\mathbf{i}\epsilon)\tau_j)\Bigr]\notag\\
&=\frac{\pi}{n}\Bigl[\sum_{\substack{j\in[n]\\a_j=R}}
\mathbf{1}_{u<-e^{\chi}\tau_j^{-1}}
+\sum_{\substack{j\in[n]\\a_j=L}}\mathbf{1}_{u>e^{\chi}\tau_j^{-1}}\Bigr].\label{agm6}
\end{align}
\end{small}
The following cases might occur:
\begin{enumerate}
    \item $u< 0$: then
    \begin{align}
    \lim_{\epsilon\rightarrow 0+}\arg [R(u+\mathbf{i}\epsilon)]^{\frac{1}{n}}=
\frac{\pi}{n}\left|\{j\in[n]:a_j=R,b_{p_{j,R}+1,j}=+\}\right|\label{agm1}
    \end{align}
    and
    \begin{align}
    \lim_{\epsilon\rightarrow 0+}\arg [U_{\chi}(u+\mathbf{i}\epsilon)]^{\frac{1}{n}}&=\frac{\pi}{n}\left|\{j\in[n]:a_j=R,
    -\tau_j<e^{\chi}u^{-1}\}\right|\label{agm2}
    \end{align}
    It is straightforward to check that there exists a unique $\Delta$ satisfying (\ref{dta1}), such that (\ref{agm1}) and (\ref{agm2}) are equal when $e^{\chi}u^{-1}\in \Delta$.
     \item $u\geq0$: then
    \begin{align}
    \lim_{\epsilon\rightarrow 0+}\arg [R(u+\mathbf{i}\epsilon)]^{\frac{1}{n}}=
\frac{\pi}{n}\left|\{j\in[n]:a_j=L,b_{p_{j,L}+1,j}=+\}\right|\label{agm3}
    \end{align}
    and
    \begin{align}
    \lim_{\epsilon\rightarrow 0+}\arg [U_{\chi}(u+\mathbf{i}\epsilon)]^{\frac{1}{n}}&=\frac{\pi}{n}\left|\{j\in[n]:a_j=L,
   e^{\chi}u^{-1}<\tau_j\}\right|\label{agm4}
    \end{align}
    It is straightforward to check that there exists a unique $\Delta$ satisfying (\ref{dta1}), such that (\ref{agm3}) and (\ref{agm4}) are equal when $e^{\chi}u^{-1}\in \Delta$.
\end{enumerate}
Then we deduce that $u\in \RR\setminus \mathcal{S}$, there exists a unique $\chi$ such that (\ref{drt2}) holds and $(\ref{agm5})$ and (\ref{agm6}) are equal. The condition that (\ref{agm5}) and (\ref{agm6}) are equal is equivalent of saying that the right hand side of (\ref{drt1}) is real and positive. When the right hand side of (\ref{drt1}) is positive, we obtain a unique $\kappa\in \RR$. Then the lemma follows.
\end{proof}

\begin{assumption}\label{ap74}
Let $i,j\in [n]$ and $p_1,p_2\in [m]$. For $a_i=a_j$ and $\tau_i>\tau_j$ holds that 
\begin{align*}
\tau_{i}^{-1}\tau_{j}<e^{V_{p_2}-V_{p_1}}.
\end{align*}
\end{assumption}

\begin{remark}Under Assumption \ref{ap74}, if we order all the points in $\{-\tau_j\}_{j\in[n]:a_j=R}\cup\{\tau_j\}_{j\in[n],a_j=L}$ as follows
\begin{align*}
-d_{n_R}<-d_{n_{R-1}}<\ldots<-d_1<0<\alpha_1<\alpha_2<\ldots<\alpha_{n_L}.
\end{align*}
Then we can order all the points in $\{e^{V_p}\tau_j^{-1}\}_{p\in[0..m],j\in[n],a_j=L}\cup
\{-e^{V_p}\tau_j^{-1}\}_{p\in[0..m],j\in[n],a_j=R}
$ as follows
\begin{align*}
&-d_1^{-1}e^{V_m}<-d_1^{-1}e^{V_{m-1}}<\ldots<-d_1^{-1}e^{V_0}<\\
&-d_2^{-1}e^{V_m}<-d_2^{-1}e^{V_{m-1}}<\ldots<-d_2^{-1}e^{V_0}<\\
&\ldots\\
&-d_{n_R}^{-1}e^{V_m}<-d_{n_R}^{-1}e^{V_{m-1}}<\ldots<-d_{n_R}^{-1}e^{V_0}<\\
&\alpha_{n_L}^{-1}e^{V_0}<\alpha_{n_L}^{-1}e^{V_1}<\ldots<\alpha_{n_L}^{-1}e^{V_m}\\
&\ldots\\
&\alpha_{1}^{-1}e^{V_0}<\alpha_{1}^{-1}e^{V_1}<\ldots<\alpha_{1}^{-1}e^{V_m}.
\end{align*}
\end{remark}

\begin{lemma}\label{l76}Suppose Assumption \ref{ap74} holds. For $u\in \HH\cup \RR$, let $(\chi_u,\kappa_u)\in \RR^2$ such that
\begin{align*}
\mathcal{G}_{\chi_u}(u)=e^{-n\kappa_u}.
\end{align*}
Assume one of the 
following two conditions holds
\begin{enumerate}
\item $u\rightarrow e^{V_p}\tau_j^{-1}\in\mathcal{S}$ for some $p\in[0..m]$, $j\in[n]$ and $a_j=L$;
\item $u\rightarrow -e^{V_p}\tau_j^{-1}\in\mathcal{S}$ for some $p\in[0..m]$, $j\in[n]$
and $a_j=R$;
\end{enumerate}
 then $\chi_u\rightarrow V_p$.
\end{lemma}
\begin{proof}
\begin{enumerate}
\item
We first consider case (1). 
\begin{enumerate}
\item Assume that $u\rightarrow e^{V_p}\tau_j^{-1}\in \mathcal{S}$ for some $p\in[m-1]$, $j\in[n]$ and $a_j=L$. 
Let $\delta>0$ be positive and small. By (\ref{dpl}), under Assumption \ref{ap74} we obtain that for $i\in[n]$, $a_i=L$, 
\begin{itemize}
\item if $u=e^{V_p}\tau_j^{-1}-\delta$
\begin{align*}
p_{i,L}=\begin{cases}-\infty &\mathrm{If}\ \tau_i<\tau_j\\
p-1 &\mathrm{If}\ \tau_i=\tau_j\\
m &\mathrm{If}\ \tau_i>\tau_j
\end{cases}
\end{align*}
\item if $u=e^{V_p}\tau_j^{-1}+\delta$
\begin{align*}
p_{i,L}=\begin{cases}-\infty &\mathrm{If}\ \tau_i<\tau_j\\
p&\mathrm{If}\ \tau_i=\tau_j\\
m &\mathrm{If}\ \tau_i>\tau_j
\end{cases}
\end{align*}
\end{itemize}
By (\ref{agm3}) we have
 \begin{align}
 \lim_{[u\rightarrow e^{V_p}\tau_j^{-1}-]}   \lim_{\epsilon\rightarrow 0+}\arg [R(u+\mathbf{i}\epsilon)]^{\frac{1}{n}}
=\frac{\pi}{n}\left|\{i\in[n]:a_i=L;
   \tau_i>\tau_j\}\right|+\mathbf{1}_{b_{p,i}=+}\label{es5}
    \end{align}
and
 \begin{align}
 \lim_{[u\rightarrow e^{V_p}\tau_j^{-1}+]}   \lim_{\epsilon\rightarrow 0+}\arg [R(u+\mathbf{i}\epsilon)]^{\frac{1}{n}}
=\frac{\pi}{n}\left|\{i\in[n]:a_i=L;
   \tau_i>\tau_j\}\right|+\mathbf{1}_{b_{p+1,i}=+}.\label{es6}
    \end{align}

Note also that
\begin{align*}
\lim_{[u\rightarrow e^{V_p}\tau_j^{-1}-]} ug(u)
=\begin{cases}
+\infty&\mathrm{if}\ \mathbf{1}_{b_{p+1,j}=+}<\mathbf{1}_{b_{p,j}=+}\\
-\infty&\mathrm{if}\ \mathbf{1}_{b_{p+1,j}=+}>\mathbf{1}_{b_{p,j}=+}\\
\mathrm{a\ finite\ real\ number}&\mathrm{if}\ \mathbf{1}_{b_{p+1,j}=+}=\mathbf{1}_{b_{p,j}=+}
\end{cases}
\end{align*}
and
\begin{align*}
\lim_{[u\rightarrow e^{V_p}\tau_j^{-1}+]} ug(u)
=\begin{cases}
-\infty&\mathrm{if}\ \mathbf{1}_{b_{p+1,j}=+}<\mathbf{1}_{b_{p,j}=+}\\
+\infty&\mathrm{if}\ \mathbf{1}_{b_{p+1,j}=+}>\mathbf{1}_{b_{p,j}=+}\\
\mathrm{a\ finite\ real\ number}&\mathrm{if}\ \mathbf{1}_{b_{p+1,j}=+}=\mathbf{1}_{b_{p,j}=+}
\end{cases}
\end{align*}
Since $e^{V_p}\tau_j^{-1}\in\mathcal{S}$, we obtain that $b_{p,j}\neq b_{p+1,j}$. We obtain that when $u\rightarrow e^{V_p}\tau_j^{-1}+$ or 
$u\rightarrow e^{V_p}\tau_j^{-1}-$, by (\ref{drt2}), $e^{\chi}u^{-1}$ approaches some $\tau_k$ for $a_k=L$, $k\in[n]$. Moreover
\begin{enumerate}
    \item If $\mathbf{1}_{b_{p,j}}>\mathbf{1}_{b_{p+1,j}}$, as $u\rightarrow e^{V_p}\tau_j^{-1}-$, $e^{\chi}u^{-1}$ approaches $\tau_k$ from the left;
    \item If $\mathbf{1}_{b_{p,j}}<\mathbf{1}_{b_{p+1,j}}$, as $u\rightarrow e^{V_p}\tau_j^{-1}-$, $e^{\chi}u^{-1}$ approaches $\tau_k$ from the right;
    \item If $\mathbf{1}_{b_{p,j}}>\mathbf{1}_{b_{p+1,j}}$, as $u\rightarrow e^{V_p}\tau_j^{-1}+$, $e^{\chi}u^{-1}$ approaches $\tau_k$ from the right;
    \item If $\mathbf{1}_{b_{p,j}}<\mathbf{1}_{b_{p+1,j}}$, as $u\rightarrow e^{V_p}\tau_j^{-1}+$, $e^{\chi}u^{-1}$ approaches $\tau_k$ from the left.
\end{enumerate}
By
(\ref{agm4}),
\begin{enumerate}
\item If $\mathbf{1}_{b_{p,j}}>\mathbf{1}_{b_{p+1,j}}$
\begin{align}
&&\lim_{[u\rightarrow e^{V_p}\tau_j^{-1}-]}\lim_{\epsilon\rightarrow 0+}\arg [U_{\chi_u}(u+\mathbf{i}\epsilon)]^{\frac{1}{n}}=\frac{\pi}{n}\left|\{i\in[n]:a_i=L;\label{es1}
  \tau_i\geq \tau_k\}\right|\\
&&\lim_{[u\rightarrow e^{V_p}\tau_j^{-1}+]}\lim_{\epsilon\rightarrow 0+}\arg [U_{\chi_u}(u+\mathbf{i}\epsilon)]^{\frac{1}{n}}=\frac{\pi}{n}\left|\{i\in[n]:a_i=L;
  \tau_i>\tau_k\}\right|\label{es2}
\end{align}
\item The case when $u\rightarrow e^{V_0}\tau_j^{-1}$ and $u\rightarrow e^{V_m}\tau_j^{-1}$ for some $j\in[n]$, $a_j=L$ can be proved similarly.
\begin{align}
&&\lim_{[u\rightarrow e^{V_p}\tau_j^{-1}-]}\lim_{\epsilon\rightarrow 0+}\arg [U_{\chi_u}(u+\mathbf{i}\epsilon)]^{\frac{1}{n}}=\frac{\pi}{n}\left|\{i\in[n]:a_i=L;
  \tau_i> \tau_k\}\right|\label{es3}\\
&&\lim_{[u\rightarrow e^{V_p}\tau_j^{-1}+]}\lim_{\epsilon\rightarrow 0+}\arg [U_{\chi_u}(u+\mathbf{i}\epsilon)]^{\frac{1}{n}}=\frac{\pi}{n}\left|\{i\in[n]:a_i=L;
  \tau_i\geq \tau_k\}\right|\label{es4}
\end{align}
\end{enumerate}
In either case to make (\ref{es1})-(\ref{es4}) equal to the corresponding arguments in (\ref{es5}), (\ref{es6}), we must have $\tau_k=\tau_j$.
\item The case $u\rightarrow e^{V_0}\tau_j^{-1}$ or
$u\rightarrow e^{V_0}\tau_j^{-1}$
for some $j\in [n]$ and $a_j=L$ can be proved similarly.
\end{enumerate}
\item Now we consider case (2). 
\begin{enumerate}
\item Assume that $u\rightarrow -e^{V_p}\tau_j^{-1}$ for some $p\in[m-1]$, $j\in[n]$ and $a_j=R$. 
Let $\delta>0$ be positive and small. By (\ref{dpr}), under Assumption \ref{ap74} we obtain that for $i\in[n]$, $a_i=R$, 
\begin{itemize}
\item if $u=-e^{V_p}\tau_j^{-1}-\delta$
\begin{align*}
p_{i,R}=\begin{cases}-\infty &\mathrm{If}\ \tau_i<\tau_j\\
p &\mathrm{If}\ \tau_i=\tau_j\\
m &\mathrm{If}\ \tau_i>\tau_j
\end{cases}
\end{align*}
\item if $u=-e^{V_p}\tau_j^{-1}+\delta$
\begin{align*}
p_{i,R}=\begin{cases}-\infty &\mathrm{If}\ \tau_i<\tau_j\\
p-1&\mathrm{If}\ \tau_i=\tau_j\\
m &\mathrm{If}\ \tau_i>\tau_j
\end{cases}
\end{align*}
\end{itemize}
By (\ref{agm1}) we have
 \begin{align}
 \lim_{[u\rightarrow -e^{V_p}\tau_j^{-1}-]}   \lim_{\epsilon\rightarrow 0+}\arg [R(u+\mathbf{i}\epsilon)]^{\frac{1}{n}}
=\frac{\pi}{n}\left|\{i\in[n]:a_i=R;
   \tau_i>\tau_j\}\right|+\mathbf{1}_{b_{p+1,i}=+}\label{et5}
    \end{align}
and
 \begin{align}
 \lim_{[u\rightarrow -e^{V_p}\tau_j^{-1}+]}   \lim_{\epsilon\rightarrow 0+}\arg [R(u+\mathbf{i}\epsilon)]^{\frac{1}{n}}
=\frac{\pi}{n}\left|\{i\in[n]:a_i=R;
   \tau_i>\tau_j\}\right|+\mathbf{1}_{b_{p,i}=+}\label{et6}
    \end{align}
Note also that
\begin{align*}
\lim_{[u\rightarrow -e^{V_p}\tau_j^{-1}-]} ug(u)
=\begin{cases}
+\infty&\mathrm{if}\ \mathbf{1}_{b_{p+1,j}=+}<\mathbf{1}_{b_{p,j}=+}\\
-\infty&\mathrm{if}\ \mathbf{1}_{b_{p+1,j}=+}>\mathbf{1}_{b_{p,j}=+}\\
\mathrm{a\ finite\ real\ number}&\mathrm{if}\ \mathbf{1}_{b_{p+1,j}=+}=\mathbf{1}_{b_{p,j}=+}
\end{cases}
\end{align*}
and
\begin{align*}
\lim_{[u\rightarrow -e^{V_p}\tau_j^{-1}+]} ug(u)
=\begin{cases}
-\infty&\mathrm{if}\ \mathbf{1}_{b_{p+1,j}=+}<\mathbf{1}_{b_{p,j}=+}\\
+\infty&\mathrm{if}\ \mathbf{1}_{b_{p+1,j}=+}>\mathbf{1}_{b_{p,j}=+}\\
\mathrm{a\ finite\ real\ number}&\mathrm{if}\ \mathbf{1}_{b_{p+1,j}=+}=\mathbf{1}_{b_{p,j}=+}
\end{cases}
\end{align*}
Since $-e^{V_p}\tau_j^{-1}\in\mathcal{S}$, we obtain that $b_{p,j}\neq b_{p+1,j}$. We obtain that when $u\rightarrow -e^{V_p}\tau_j^{-1}+$ or 
$u\rightarrow -e^{V_p}\tau_j^{-1}-$, by (\ref{drt2}), $e^{\chi}u^{-1}$ approaches some $-\tau_k$ for $a_k=R$, $k\in[n]$. Moreover
\begin{enumerate}
    \item If $\mathbf{1}_{b_{p,j}}>\mathbf{1}_{b_{p+1,j}}$, as $u\rightarrow -e^{V_p}\tau_j^{-1}-$, $e^{\chi}u^{-1}+\tau_k$ approaches $0$ from the left;
    \item If $\mathbf{1}_{b_{p,j}}<\mathbf{1}_{b_{p+1,j}}$, as $u\rightarrow -e^{V_p}\tau_j^{-1}-$, $e^{\chi}u^{-1}+\tau_k$ approaches $0$ from the right;
    \item If $\mathbf{1}_{b_{p,j}}>\mathbf{1}_{b_{p+1,j}}$, as $u\rightarrow-e^{V_p}\tau_j^{-1}+$, $e^{\chi}u^{-1}+\tau_k$ approaches $0$ from the right;
    \item If $\mathbf{1}_{b_{p,j}}<\mathbf{1}_{b_{p+1,j}}$, as $u\rightarrow -e^{V_p}\tau_j^{-1}+$, $e^{\chi}u^{-1}+\tau_k$ approaches $0$ from the left.
\end{enumerate}
By
(\ref{agm2}),
\begin{enumerate}
\item If $\mathbf{1}_{b_{p,j}}>\mathbf{1}_{b_{p+1,j}}$
\begin{align}
&&\lim_{[u\rightarrow -e^{V_p}\tau_j^{-1}-]}\lim_{\epsilon\rightarrow 0+}\arg [U_{\chi_u}(u+\mathbf{i}\epsilon)]^{\frac{1}{n}}=\frac{\pi}{n}\left|\{i\in[n]:a_i=R;\label{et1}
  \tau_i> \tau_k\}\right|\\
&&\lim_{[u\rightarrow -e^{V_p}\tau_j^{-1}+]}\lim_{\epsilon\rightarrow 0+}\arg [U_{\chi_u}(u+\mathbf{i}\epsilon)]^{\frac{1}{n}}=\frac{\pi}{n}\left|\{i\in[n]:a_i=R;
  \tau_i\geq\tau_k\}\right|\label{et2}
\end{align}
\item If $\mathbf{1}_{b_{p,j}}<\mathbf{1}_{b_{p+1,j}}$
\begin{align}
&&\lim_{[u\rightarrow -e^{V_p}\tau_j^{-1}-]}\lim_{\epsilon\rightarrow 0+}\arg [U_{\chi_u}(u+\mathbf{i}\epsilon)]^{\frac{1}{n}}=\frac{\pi}{n}\left|\{i\in[n]:a_i=R;
  \tau_i\geq\tau_k\}\right|\label{et3}\\
&&\lim_{[u\rightarrow -e^{V_p}\tau_j^{-1}+]}\lim_{\epsilon\rightarrow 0+}\arg [U_{\chi_u}(u+\mathbf{i}\epsilon)]^{\frac{1}{n}}=\frac{\pi}{n}\left|\{i\in[n]:a_i=R;
  \tau_i> \tau_k\}\right|\label{et4}
\end{align}
\end{enumerate}
In either case to make (\ref{et1})-(\ref{et4}) equal to the corresponding arguments in (\ref{et5}), (\ref{et6}), we must have $\tau_k=\tau_j$.
\item The case when $u\rightarrow -e^{V_0}\tau_j^{-1}$ or $u\rightarrow -e^{V_m}\tau_j^{-1}$ for some $j\in[n]$, $a_j=R$ can be proved similarly.
\end{enumerate}
\end{enumerate}
\end{proof}

\begin{theorem}\label{t77}Let $\{RYG(l^{(\epsilon)},r^{(\epsilon)},\underline{a}^{(\epsilon)},\underline{b}^{(\epsilon)})\}_{\epsilon>0}$ be a sequence of rail-yard graphs satisfying Assumptions \ref{ap5}, \ref{ap65}(2)(3), and \ref{ap74}. Let $\mathbf{w}_+:\mathcal{L}\rightarrow \HH$ be the diffeomorphism from the liquid region to the upper half plane which maps each point $(\chi,\kappa)$ in the liquid region to the unique root of (\ref{ceqn}) in the upper half plane $\HH$. Then as $\epsilon\rightarrow 0$, the height function of pure dimer coverings on $\{RYG(l^{(\epsilon)},r^{(\epsilon)},\underline{a}^{(\epsilon)},\underline{b}^{(\epsilon)})\}_{\epsilon>0}$ in the liquid region converges to the $\mathbf{w}_+$-pullback of GFF in the sense that for any $(\chi,\kappa)\in \mathcal{L}$, $\chi\notin\{V_p\}_{p=0}^{m}$ and positive real number $\alpha$
\begin{align*}
\int_{-\infty}^{\infty}\left(h_{M}\left(\frac{\chi}{\epsilon},\frac{\kappa}{\epsilon}\right)-
\mathbb{E}\left[h_{M}\left(\frac{\chi}{\epsilon},\frac{\kappa}{\epsilon}\right)
\right]\right)
e^{-\alpha\kappa}d\kappa\longrightarrow\int_{(\chi,\kappa)\in\mathcal{L}}e^{-\alpha \kappa} \Xi(\mathbf{w}_+(\chi,\kappa)) d\kappa
\end{align*}
in distribution.
\end{theorem}

\begin{proof}Let $\chi\in [r^{(0)},l^{(0)}]$ and $k$ be a positive integer. By (\ref{hme}) and Assumption \ref{ap5}, we have
\begin{align*}
\int_{-\infty}^{\infty}\left(h_{M}\left(\frac{\chi}{\epsilon},\frac{\kappa}{\epsilon}\right)-
\mathbb{E}\left[h_{M}\left(\frac{\chi}{\epsilon},\frac{\kappa}{\epsilon}\right)
\right]\right)
e^{-n\beta k\kappa}d\kappa=
\frac{2\epsilon\left[\gamma_k(\lambda^{(m)},t;t)-
\mathbb{E}\gamma_k(\lambda^{(m)},t;t)
\right]}{(k\log t)^2}
\end{align*}
where $\chi=2m-\frac{1}{2}$.
By Theorem \ref{t58}, we obtain that for
\begin{align*}
l^{(0)}<\chi_1<\chi_2<\ldots<\chi_s<r^{(0)}
\end{align*}
and positive integers $k_1,\ldots,k_s$
\begin{align*}
\left\{\int_{-\infty}^{\infty}\left(h_{M}\left(\chi_i,\frac{y}{\epsilon}\right)-
\mathbb{E}\left[h_{M}\left(\chi_i,\frac{y}{\epsilon}\right)
\right]\right)
t^{-k_iy}dy\right\}_{i\in[s]}
\end{align*}
converges to the Gaussian vector with covariance
\begin{align*}
I:=\frac{1}{k_ik_jn^2\beta^2(\pi\mathbf{i})^2}\oint_{\mathcal{C}_w}\oint_{\mathcal{C}_z} \frac{\left[\mathcal{G}_{\chi_i}(z)\right]^{k_i\beta}\left[\mathcal{G}_{\chi_j}(w)\right]^{k_j\beta}}{(z-w)^2}dz dw.
\end{align*}


 Under Assumption \ref{ap65}, we deform the integral contour $\mathcal{C}_w$ to $\widetilde{\mathcal{C}}_w$ such that
 \begin{enumerate}
     \item $\widetilde{\mathcal{C}}_w=C_{w,1}\cup C_{w,2}$;
     \item $C_{w,1}$ lies in the upper half plane except two endpoints along the real axis;
     \item $C_{w,2}$ is the reflection of $C_{w,1}$ along the real axis;
     \item $[\mathbf{w}_+]^{-1}(C_{w,1})$ is the vertical line in $\mathcal{L}$ passing through $(\chi_j,0)$.
 \end{enumerate}
Similarly, we deform the integral contour $\mathcal{C}_z$ to $\widetilde{\mathcal{C}}_z$ such that
 \begin{enumerate}
     \item $\widetilde{\mathcal{C}}_z=C_{z,1}\cup C_{z,2}$;
     \item $C_{z,1}$ lies in the upper half plane except two endpoints along the real axis;
     \item $C_{z,2}$ is the reflection of $C_{z,1}$ along the real axis;
     \item $[\mathbf{w}_+]^{-1}(C_{z,1})$ is the vertical line in $\mathcal{L}$ passing through $(\chi_i,0)$.
 \end{enumerate}
 Then making a change of variables from $(z,w)\in \CC^2$ to $((\chi_1,\kappa_1),(\chi_2,\kappa_2))\in \mathcal{L}^2$ by $[\mathbf{w}_+]^{-1}\times [\mathbf{w}_+]^{-1}$ and the corresponding complex conjugates, we obtain
 \begin{align*}
 I&=\frac{1}{k_ik_jn^2\beta^2(\pi\mathbf{i})^2}\oint_{\widetilde{\mathcal{C}}_w}\oint_{\widetilde{\mathcal{C}}_z} \frac{\left[\mathcal{G}_{\chi_i}(z)\right]^{k_i\beta}\left[\mathcal{G}_{\chi_j}(w)\right]^{k_j\beta}}{(z-w)^2}dz dw=\frac{1}{k_ik_jn^2\beta^2(\pi\mathbf{i})^2}\\
&\times \left[\int_{(\chi_j,\kappa_j)\in\mathcal{L}}\int_{(\chi_i,\kappa_i)\in\mathcal{L}} \frac{e^{-n\kappa_ik_i\beta}e^{-n\kappa_jk_j\beta}}{(\mathbf{w}_+(\chi_i,\kappa_i)-\mathbf{w}_+(\chi_j,\kappa_j))^2}\frac{\partial \mathbf{w}_+(\chi_i,\kappa_i)}{\partial \kappa_i}\frac{\partial \mathbf{w}_+(\chi_j,\kappa_j)}{\partial \kappa_j}d\kappa_i d\kappa_j\right.\\
  &-\int_{(\chi_j,\kappa_j)\in\mathcal{L}}\int_{(\chi_i,\kappa_i)\in\mathcal{L}} \frac{e^{-n\kappa_ik_i\beta}e^{-n\kappa_jk_j\beta}}{(\mathbf{w}_+(\chi_i,\kappa_i)-\overline{\mathbf{w}_+(\chi_j,\kappa_j)})^2}\frac{\partial \mathbf{w}_+(\chi_i,\kappa_i)}{\partial \kappa_i}\frac{\partial\overline{ \mathbf{w}_+(\chi_j,\kappa_j)}}{\partial \kappa_j}d\kappa_i d\kappa_j\\
   &-\int_{(\chi_j,\kappa_j)\in\mathcal{L}}\int_{(\chi_i,\kappa_i)\in\mathcal{L}} \frac{e^{-n\kappa_ik_i\beta}e^{-n\kappa_jk_j\beta}}{\overline{(\mathbf{w}_+(\chi_i,\kappa_i)}-\mathbf{w}_+(\chi_j,\kappa_j))^2}\frac{\partial \mathbf{w}_+(\chi_i,\kappa_i)}{\partial \kappa_i}\frac{\partial \overline{\mathbf{w}_+(\chi_j,\kappa_j)}}{\partial \kappa_j}d\kappa_i d\kappa_j\\
   &+\left.\int_{(\chi_j,\kappa_j)\in\mathcal{L}}\int_{(\chi_i,\kappa_i)\in\mathcal{L}} \frac{e^{-n\kappa_ik_i\beta}e^{-n\kappa_jk_j\beta}}{\overline{(\mathbf{w}_+(\chi_i,\kappa_i)}-\overline{\mathbf{w}_+(\chi_j,\kappa_j)})^2}\frac{\partial \overline{\mathbf{w}_+(\chi_i,\kappa_i)}}{\partial \kappa_i}\frac{\partial \overline{\mathbf{w}_+(\chi_j,\kappa_j)}}{\partial \kappa_j}d\kappa_i d\kappa_j\right]\\
 \end{align*}

Integrating by parts, we obtain that 
\begin{align*}
I&=\frac{2}{(\pi\mathbf{i})^2}\int_{(\chi_j,\kappa_j)\in\mathcal{L}}\int_{(\chi_i,\kappa_i)\in\mathcal{L}}
e^{-n\kappa_ik_i\beta}e^{-n\kappa_jk_j\beta}
\log\left|\frac{\mathbf{w}_+(\chi_i,\kappa_i)-\mathbf{w}_+(\chi_j,\kappa_j)}{
\mathbf{w}_+(\chi_i,\kappa_i)-\overline{\mathbf{w}_+(\chi_j,\kappa_j)}
}\right|
d\kappa_i d\kappa_j\\
&=4\mathrm{Cov}\left(
\int_{(\chi_i,\kappa_i)\in\mathcal{L}}e^{-n\kappa_ik_i\beta} \Xi(\mathbf{w}_+(\chi_i,\kappa_i)) d\kappa_i,
\int_{(\chi_j,\kappa_j)\in\mathcal{L}} \Xi(\mathbf{w}_+(\chi_j,\kappa_j))
e^{-n\kappa_jk_j\beta}d\kappa_j
\right).
\end{align*}
Then the proposition follows.
\end{proof}

\section{Pyramid Partitions}\label{sect:ex}

In this section, we discuss specific example of the rail-yard graph, known as pyramid partitions. The limit shape and height fluctuations of perfect matchings on these graphs can be obtained by the technique developed in the paper. 

\subsection{Pyramid partitions}

A fundamental pyramid partition is a heap of square bricks such that
\begin{itemize}
    \item Each square brick is of size $2\times 2$ and has a central line dividing it into two equal-size rectangular parts; hence the direction of the central line determines the direction of the square brick; and
    \item Each square brick lies upon two side-by-side square bricks; and is rotated 90 degrees from the bricks immediately below it; and
    \item there is a unique brick on the top.
\end{itemize}

A pyramid partition is obtained from the fundamental pyramid partition by removing finitely many square bricks, such that if a square brick is removed, then all the square bricks above it are also removed. See the first figure in Introduction. 

Let $s$ be a fixed positive integer which is odd. Let $\Lambda_s$ be the set of pyramid partitions that can be obtained from the fundamental partition where the center of the square brick on the top is $(0,0)$ and where we can only take off bricks that lie inside the strip $-s-1\leq x-y\leq s+1$.

Looking from the top each pyramid partition corresponds to a domino tiling of the square grid. See the second figure in Introduction. From a pyramid partition, we can obtain a pure dimer covering on a rail-yard graph by the following steps:
\begin{enumerate}
\item rotate the pyramid partition clockwise by 45 degrees,  
\item For each blue vertex $v_b$, assume it has 4 incident edges $e_1,e_2,e_3,e_4$. Assume that $e_1$ and $e_2$ (resp.\ $e_3$ and $e_4$) are to the left (resp.\ right) of $v_b$. Split each blue vertex $v_b$ of the dual graph into 3 vertices, $v_{b_1}$, $v_{b_2}$, $v_{b_3}$ such that $v_{b_1}$ and $v_{b_3}$ are blue vertices while $v_{b_2}$ is a red vertex. The red vertex $v_{b_2}$ has exactly two incident edges joining it to $v_{b_1}$ and $v_{b_3}$, respectively. $v_{b_1}$ has 3 incident edges $e_1$ $e_2$ and $(v_{b_1},v_{b_2})$; while $v_{b_3}$ has 3 incident edges $e_3$ $e_4$ and $(v_{b_3},v_{b_2})$.
\item If one of $e_1,e_2$ (resp.\ $e_3,e_4$) is in the dimer covering, while neither $e_3$ nor $e_4$ (resp.\ neither $e_1$ nor $e_2$) are in the dimer covering, make $(v_{b_2},v_{b_3})$ (resp.\ $(v_{b_1},v_{b_2})$) present in the dimer covering and $(v_{b_1},v_{b_2})$ (resp. $(v_{b_2},v_{b_3})$) absent in the dimer covering.
\end{enumerate}
See the third and fourth figures in Introduction for 
for the pure dimer covering on a rail-yard graph corresponding to the pyramid partitions given as the examples.

\begin{proposition}
There is a one-to-one correspondence between pyramid partitions in $\Lambda_s$ and pure dimer coverings on the rail-yard graph such that for $i\in[-s..s-1]$
\begin{align*}
a_i=\begin{cases}
L& i \textrm{ is odd}\\
R& i \textrm{ is even}
\end{cases} 
\quad {and } \quad
b_i=\begin{cases}
+& i<0 \\
-& i\geq 0.
\end{cases} 
\end{align*}
Equivalently, there is a bijection between pyramid partitions in $\Lambda_s$ and sequences of partitions $(\lambda^{(-s)},\lambda^{(-s+1)},\ldots,\lambda^{(0)},\lambda^{(1)},\ldots,\lambda^{(s)})$ such that
\begin{align*}
\emptyset=\lambda^{(-s)}\prec\lambda^{(-s+1)}\prec'\lambda^{(-s+2)}\ldots\prec \lambda^{(0)}\succ'\lambda^{(1)}\succ\lambda^{(2)}\ldots\succ'\lambda^{(s)}=\emptyset.
\end{align*}
\end{proposition}
\begin{proof}
See Lemma 5.9 of \cite{you10} and Proposition 8 of \cite{BCC17}.
\end{proof}

The formula to compute partition function of pyramid partitions was conjectured in \cite{Ken05,Sze08} and proved in \cite{BY09,you10}.

Consider the pure dimer coverings on rail-yard graphs corresponding to pyramid partitions. Then we have $m=2$, $V_1=0$ and $V_0=-V_2$. Assume that the model is periodic with $n=2$.

Recall that $\mathcal{G}_{\chi}$ is defined by (\ref{dgc}). Then the frozen boundary has the following parametric equation (parametrized by $w$):
\begin{align*}
\begin{cases}
\frac{U_{\chi}(w)}{R(w)}=e^{-2\kappa}\\
f(e^{\chi}w^{-1})=wg(w)
\end{cases},
\end{align*}
where
\begin{align*}
f(s):=\frac{1}{1-\tau_1^{-1}s}
-\frac{1}{1+\tau_2^{-1}s},
\end{align*}
\begin{align*}
g(w):&=-\frac{1}{w-e^{V_1}\tau_1^{-1}}
+\frac{1}{w-e^{V_2}\tau_1^{-1}}
+\frac{1}{w-e^{V_0}\tau_1^{-1}}\\
&+\frac{1}{w+e^{V_1}\tau_2^{-1}}
-\frac{1}{w+e^{V_2}\tau_2^{-1}}
-\frac{1}{w+e^{V_0}\tau_2^{-1}},
\end{align*}
and
\begin{align*}
U_{\chi}(w)&=\frac{\left(1+e^{-\chi}w \tau_2\right)}{({1-e^{-\chi}w \tau_1})},
\end{align*}
\begin{align*}
R(w)&=\frac{\left(1+e^{-V_0}w\tau_2\right)\left(1+e^{-V_2}w\tau_2\right)\left(1-e^{-V_1}\tau_1 w\right)}{\left(1-e^{-V_0}w\tau_1\right)\left(1-e^{-V_2}w\tau_1\right)\left(1+e^{-V_1}\tau_2 w\right)}.
\end{align*}
By (\ref{dpl}), (\ref{dpr}), we obtain
\begin{align*}
p_{1,L}&=\max\{p\in\{0,1,2\}:e^{V_p}\tau_1^{-1}<w\};\\
p_{2,R}&=\max\{p\in\{0,1,2\}:w<-e^{V_p}\tau_2^{-1}\}.
\end{align*}
By (\ref{agm1})-(\ref{agm4}), we have
\begin{itemize}
    \item $w< 0$: then
    \begin{align*}
    \lim_{\epsilon\rightarrow 0+}\arg [R(w+\mathbf{i}\epsilon)]^{\frac{1}{2}}=
\frac{\pi}{2}\mathbf{1}_{b_2(p_{2,R},p_{2,R}+1)=+}
    \end{align*}
    and
    \begin{align*}
    \lim_{\epsilon\rightarrow 0+}\arg [U_{\chi}(w+\mathbf{i}\epsilon)]^{\frac{1}{2}}&=\frac{\pi}{2}\mathbf{1}_
    {-\tau_2<e^{\chi}w^{-1}}
    \end{align*}
   \item $w\geq0$: then
    \begin{align*}
    \lim_{\epsilon\rightarrow 0+}\arg [R(w+\mathbf{i}\epsilon)]^{\frac{1}{2}}=
\frac{\pi}{2}\mathbf{1}_{b_1(p_{1,L},p_{1,L}+1)=+}
    \end{align*}
    and
    \begin{align*}
    \lim_{\epsilon\rightarrow 0+}\arg [U_{\chi}(w+\mathbf{i}\epsilon)]^{\frac{1}{2}}&=
   \mathbf{1}_{e^{\chi}w^{-1}<\tau_1}
    \end{align*}
In order to make   
\begin{align*}
\lim_{\epsilon\rightarrow 0+}\arg [R(w+\mathbf{i}\epsilon)]^{\frac{1}{2}}= \lim_{\epsilon\rightarrow 0+}\arg [U_{\chi}(w+\mathbf{i}\epsilon)]^{\frac{1}{2}},
\end{align*}
 we have
 \begin{enumerate}
     \item If $w>e^{V_2}\tau_1^{-1}$, $e^{\chi}w^{-1}\in (0,\tau_1)$;
     \item If $w\in(e^{V_1}\tau_1^{-1},e^{V_2}\tau_1^{-1})$, $e^{\chi}w^{-1}\in (\tau_1,\infty)$;
     \item If $w\in(e^{V_0}\tau_1^{-1},e^{V_1}\tau_1^{-1})$, $e^{\chi}w^{-1}\in (0,\tau_1)$;
     \item If $w\in(0,e^{V_0}\tau_1^{-1})$, $e^{\chi}w^{-1}\in (\tau_1,\infty)$;
     \item If $w<-e^{V_2}\tau_2^{-1}$, $e^{\chi}w^{-1}\in (-\tau_2,0)$;
     \item If $w\in(-e^{V_2}\tau_2^{-1},-e^{V_1}\tau_2^{-1})$, $e^{\chi}w^{-1}\in (-\infty,-\tau_2)$;
     \item If $w\in(-e^{V_1}\tau_2^{-1},-e^{V_0}\tau_2^{-1})$, $e^{\chi}w^{-1}\in (-\tau_2,0)$;
     \item If $w\in(-e^{V_0}\tau_2^{-1},0)$, $e^{\chi}w^{-1}\in (-\infty,-\tau_2)$.
    \end{enumerate}
\end{itemize}
Hence for each $w\in \mathbb{R}\setminus \{\pm e^{V_p}\tau_j^{-1},0\}_{p\in\{0,1,2\},j\in\{1,2\}}$, we can find a unique $\chi$ satisfying (1)-(8) and $f(e^{\chi}w^{-1})=wg(w)$; then knowing $w$ and $\chi$ we can find a unique $\kappa$ by $\frac{U_{\chi}(w)}{R(w)}=e^{-n\kappa}$. See Figure \ref{limshapePP} from Introduction for the frozen boundary of pyramid partitions.

\appendix

\section{}\label{sc:dmp}
Here we recall some facts about Macdonald polynomials and include some known technical results that were used in this paper. 

Let $X=(x_1,\ldots,x_n,\ldots)$ and $Y=(y_1,\ldots,y_n,\ldots)$ be two countable sets of variables. Let $\Lambda_X$ be the algebra of symmetric functions of $X$ over $\CC$. The power symmetric functions $\{p_{\lambda}(X)\}_{\lambda\in \YY}$ form a linear basis for $\Lambda_X$, where
\begin{align*}
p_{\lambda}(X)=\prod_{i\in\NN}p_{\lambda_i}(X) \quad \textrm{and} \quad p_i(X)=\sum_{j\in \NN}x_j^i,\textrm{for }i\in\NN.
\end{align*}

For each fixed pair of parameters $q,t\in(0,1)$ and $\lambda,\mu\in\YY$ define the scalar product $\langle \cdot,\cdot \rangle: \Lambda_X\times \Lambda_X\rightarrow \RR$ as a bilinear map such that:
\begin{align}
\langle p_{\lambda},p_{\mu}  \rangle=\delta_{\lambda\mu}\left[\prod_{i=1}^{l(\lambda)}\frac{1-q^{\lambda_i}}{1-t^{\lambda_i}}\right]\left[\prod_{j=1}^{\infty}j^{m_j(\lambda)}(m_j(\lambda))!\right]\label{dsp}
\end{align}
where $\delta_{\lambda \mu}=1$ if and only if $\lambda=\mu$, and $m_j(\lambda)$ is the number of parts in $\lambda$ equal to $j$. 

Macdonald symmetric functions $P_\lambda(X;q,t)$ and $Q_{\lambda}(X;q,t)$, for the definition see (4.7) and (4.12) in Chapter VI of \cite{IGM15}, form two bases $(P_{\lambda})$ and $(Q_\lambda)$, which are dual to each other with respect to the above scalar product, i.e.  
\begin{align*}
\langle P_{\lambda}(X;q,t),Q_{\lambda}(X;q,t) \rangle=\delta_{\lambda\mu}.
\end{align*}
Skew Macdonald symmetric functions are defined by the branching rules
\begin{align*}
P_{\lambda}(X,Y;q,t)&=\sum_{\mu\in \YY}P_{\lambda/\mu}(X;q,t)P_{\mu}(Y;q,t),\\
Q_{\lambda}(X,Y;q,t)&=\sum_{\mu\in \YY}Q_{\lambda/\mu}(X;q,t)Q_{\mu}(Y;q,t).
\end{align*}
When $q=t$, 
\begin{align*}
&P_{\lambda}(X;t,t)=Q_{\lambda}(X;t,t)=s_{\lambda}(X),\notag\\
&P_{\lambda/\mu}(X;t,t)=Q_{\lambda/\mu}(X;t,t)=s_{\lambda/\mu}(X).
\end{align*}
It is known that (see Remarks 1.~on Page 346 of \cite{IGM15}) for a single variable $x$
\begin{align*}
P_{\lambda/\mu}(x)=\delta_{\mu\prec\lambda}\psi_{\lambda/\mu}(q,t)x^{|\lambda|-|\mu|},\qquad
Q_{\lambda/\mu}(x)=\delta_{\mu\prec\lambda}\phi_{\lambda/\mu}(q,t)x^{|\lambda|-|\mu|},
\end{align*}
where $\psi_{\lambda/\mu}(q,t)$ and $\phi_{\lambda/\mu}(q,t)$ are independent of $x$ and furthermore
\begin{align*}
\left.\psi_{\lambda/\mu}(q,t)\right|_{q=t}=
\left.\phi_{\lambda/\mu}(q,t)\right|_{q=t}=1.
\end{align*}

\begin{definition}
Let $k\in \mathbb{Z}^{+}$ and $q,t\in \mathbb{R}^{+}$ be parameters. Let $D_{-k,X}$ be an operator acting on symmetric functions $\Lambda_X$.
For any analytic symmetric function $F(X)$ satisfying
\begin{align*}
F(X) =\sum_{\lambda \in \YY}c_{\lambda}P_{\lambda}(X;q,t),
\end{align*}
where $c_{\lambda}$'s are complex coefficients, define $D_{-k,X}F\in\Lambda_{X}$ to be
\begin{align}
D_{-k,X;q,t}F(X) =\sum_{\lambda\in\YY} c_{\lambda}\left\{(1-t^{-k})\left[\sum_{i=1}^{l(\lambda)}(q^{\lambda_i}t^{-i+1})^k\right]+t^{-kl(\lambda)}\right\}P_{\lambda}(X;q,t).\label{ngt}
\end{align}
\end{definition}

Let $W=(w_1,\ldots,w_k)$ be an ordered set of variables. Define
\begin{align}
D(W;q,t)=\frac{(-1)^{k-1}}{(2\pi\mathbf{i})^{k}}\frac{\sum_{i=1}^k\frac{w_k t^{k-i}}{w_i q^{k-i}}}{\left(1-\frac{tw_2}{qw_1}\right)\ldots\left(1-\frac{tw_k}{qw_{k-1}}\right)}\prod_{i<j}\frac{(1-\frac{w_i}{w_j})(1-\frac{qw_i}{tw_j})}{\left(1-\frac{w_i}{tw_j}\right)\left(1-\frac{qw_i}{w_j}\right)}\prod_{i=1}^k\frac{dw_i}{w_i}.\label{ddf}
\end{align}
Recall that $H(W,X;q,t)$ was defined as in (\ref{dh}).

The following proposition is a slightly more general form of Proposition 4.10 of \cite{GZ16}.
\begin{proposition}\label{pa2}Assume one of the following two conditions holds
\begin{enumerate}
    \item $q\in(0,1)$ and $t\in(0,1)$; or
    \item $q\in (1,\infty)$ and $t\in (1,\infty)$.
\end{enumerate}

Let $f:\CC\rightarrow\CC$ be a function analytic in a neighborhood of 0, and $f(0)\neq 0$. Let $g:\CC\rightarrow\CC$ be a function analytic in a neighborhood of 0, and 
\begin{align*}
g(z)=\frac{f(z)}{f(q^{-1}z)};
\end{align*}
for $z$ in a small neighborhood of $0$. Then
\begin{align}
D_{-k,X;q,t}\left(\prod_{x_i\in X} f(x_i)\right)=\left(\prod_{x_i\in X}f(x_i)\right)\label{ss}
\oint\cdots\oint D(W;q,t)H(W,X;q,t)\left(\prod_{i=1}^k g(w_i)\right)
\end{align}
where the contours of the integral satisfy the following conditions 
\begin{itemize}
\item all the contours are in the neighborhood of $0$ such that both $f$ and $g$ are analytic;
\item each contour enclose $0$ and $\{qx_i\}_{x_i\in X}$;
\item If case (1) holds, $|w_i|\leq |tw_{i+1}|$ for all $i\in[k-1]$;
\item If case (2) holds, $|w_i|\leq \left|\frac{1}{q}w_{i+1}\right|$ for all $i\in[k-1]$;
\end{itemize}
$H(W,X;q,t)$ is given by (\ref{dh}), and $D(W;q,t)$ is given by (\ref{ddf}).
\end{proposition}

\begin{proof}When $X$ consists of finitely many variables and when case (1) holds, the proposition was proved in Proposition 4.10 of \cite{GZ16}. It is straightforward to check the Proposition when case (2) holds by (\ref{pqr}).

When $X$ consists of countably many variables, the identity (\ref{ss}) holds formally, since its projection onto any finitely many variables $(x_1,\ldots,x_n)$ by letting $x_{n+1}=x_{n+2}=\ldots=0$ holds.
\end{proof}

\begin{lemma}\label{la3}
Let $(a,q)_{\infty}=\prod_{r=0}^{\infty}(1-aq^r)$ and
\begin{align}
\Pi(X,Y;q,t):=\frac{(tx_iy_j;q)_{\infty}}{(x_iy_j;q)_{\infty}}, \quad \Pi'(X,Y):=\prod_{i,j}(1+x_iy_j).\label{defPPprime}
\end{align}
Then
\begin{align*}
\sum_{\lambda\in\YY}P_{\lambda}(X;q,t)Q_{\lambda}(Y;q,t)&=\sum_{\lambda\in\YY}P_{\lambda'}(X;q,t)Q_{\lambda'}(Y;q,t)=\Pi(X,Y;q,t);\\
\sum_{\lambda\in \YY}P_{\lambda}(X;q,t)P_{\lambda'}(Y;t,q)&=\sum_{\lambda\in\YY}Q_{\lambda}(X;q,t)Q_{\lambda'}(Y;t,q)=\Pi'(X,Y).
\end{align*}
In particular, when $q=t$ we obtain the Cauchy identities for Schur polynomials.
\begin{align*}
&\sum_{\lambda\in \YY}s_{\lambda}(X)s_{\lambda}(Y)=\prod_{i,j}\frac{1}{1-x_iy_j};\\
&\sum_{\lambda\in \YY}s_{\lambda}(X)s_{\lambda'}(Y)=\prod_{i,j}(1+x_iy_j).
\end{align*}
\end{lemma}
\begin{proof}See (2.5), (4.13), and (5.4) in Section VI of \cite{IGM15}. 
\end{proof}

\begin{lemma}\label{la6} Let $\Pi$, $\Pi'$, and $H$ be as in (\ref{defPPprime}) and (\ref{dh}), then
\begin{align*}
&\Pi(X,Y;q,t)=\exp\left(\sum_{n=1}^{\infty}\frac{1-t^n}{1-q^n}\frac{1}{n}p_n(X)p_n(Y)\right),\\
&\Pi'(X,Y)=\exp\left(\sum_{n=1}^{\infty}\frac{(-1)^{n+1}}{n}p_n(X)p_n(Y)\right),\\
&H(X,Y;q,t)=\exp\left(\sum_{n=1}^{\infty}\frac{1-t^{-n}}{n}p_n(qX^{-1})p_n(Y)\right).
\end{align*}
\end{lemma}

\begin{proof}The first identity follows from Page 310 of \cite{IGM15}. The other two follow from \begin{align*}
\Pi'(X,Y)&=\left[\Pi(-X,Y;0,0)\right]^{-1},\\
H(X,Y;q,t)&=\Pi(qX^{-1},Y;0,t^{-1}).
\end{align*}
\end{proof}

\begin{definition}Let $\mathcal{A}$ be a graded algebra over a field $F$. For $a\in \mathcal{A}$, define $\mathrm{ldeg}(a)$ to be the minimum degree of all the homogeneous components in $a$.
\end{definition}

\begin{lemma}\label{la5}(Proposition 2.3 pf \cite{bcgs13}) Let $\{d_k\}_k$ $\{u_k\}_k$ be two sequences of elements of graded algebras $\mathcal{A}$ and $\mathcal{B}$. Assume $\lim_{k\rightarrow\infty}\mathrm{ldeg}(d_k)=\infty$ and 
$\lim_{k\rightarrow\infty}\mathrm{ldeg}(u_k)=\infty$. For non-negative integer $k$, let $p_k$ be the power sum. Then
\begin{align*}
\left\langle \exp\left(\sum_{k=1}^{\infty}\frac{d_kp_k(Y)}{k}\right), 
\exp\left(\sum_{k=1}^{\infty}\frac{u_kp_k(Y)}{k}\right)
\right\rangle_Y =\exp\left(\sum_{k=1}^{\infty}\left(\frac{1-q^k}{1-t^k}\cdot\frac{d_k u_k}{k}\right)\right),
\end{align*}
where $a_k$, $b_k$ are independent of the variables in $Y$.
\end{lemma}

\begin{definition}Let $F\supset \CC$ be a field. Let $\mathcal{A}$ be a ($\ZZ_{\geq 0}$-)graded algebra over $F$. For each nonnegative integer $n$, let $\mathcal{A}_n$ denote the $n$-th homogeneous component of $\mathcal{A}$. 

The completion $\widehat{\mathcal{A}}$ consists of formal sums $\sum_{n=1}^{\infty}a_n$ where $a_n\in \mathcal{A}_n$.
For two graded algebras $\mathcal{A}$, $\mathcal{A}'$ over $F$, let $\mathcal{A}\otimes_F\mathcal{A}'$ be a graded algebra over $F$ such that for $a\in \mathcal{A}_m$ and $a'\in \mathcal{A}'_n$, $a\otimes a'\in (\mathcal{A}\otimes_F{A}')_{m+n}$. Let $\mathcal{A}\widehat{\otimes}_F\mathcal{A}'$ be the completion of $\mathcal{A}\otimes_{F}\mathcal{A}'$.

If $\mathcal{B}$ is a graded algebra over $\CC$, let $\mathcal{B}_F$ be the graded algebra $\mathcal{B}\otimes_{\CC}F$ over $F$, i.e. the extension of coefficients from $\CC$ to $F$.
Let  $\Lambda_{X}[F]$  denote  the $F$-algebra  of  symmetric  functions  in $X=\{x_1,x_2,...\}$,  with coefficients in $F$. 
\end{definition}

\begin{definition}
Let $\mathcal{A}$ and $\mathcal{A}'$ be graded algebras over $\CC$ and $\{a_{n,j}\}_j$ be a basis for $\mathcal{A}_n$ for each $n\geq 0$.  We say that an element $f\in \mathcal{A}\widehat{\otimes}\mathcal{A}'[F]$ is $\mathcal{A}$-projective if 
\begin{align*}
f=\sum_{n,j}a_{n,j}\otimes \alpha'_{n,j}, \qquad \alpha'_{n,j}\in\mathcal{A}'(F)
\end{align*}
such that $\lim_{n\rightarrow\infty}\min_j\mathrm{ldeg}(\alpha'_{n,j})=\infty.$ This property is independent of the choice of basis.
\end{definition}

\begin{definition}Let $\mathcal{A}$, $\mathcal{B}$ be graded algebras over $\CC$, and let $F\supset \CC$ be a field. Define the Macdonald scalar product to be the bilinear map
\begin{align*}
\left(\mathcal{A}\otimes \Lambda_X\right)[F]\times (\Lambda_X\otimes \mathcal{B})[F]\rightarrow \mathcal{A}\otimes \mathcal{B}[F]
\end{align*}
such that 
\begin{align*}
\langle a\otimes P_{\lambda}, Q_{\mu}\otimes b \rangle_X:=\langle P_{\lambda}, Q_{\mu} \rangle a\otimes b =\delta_{\lambda\mu}a\otimes b. 
\end{align*}
\end{definition}

\begin{definition}Let $Z:=(z_{1},\ldots,z_{k}),$ where $k$ is a positive integer. Let $\mathcal{L}(Z)$ be the field of formal Laurent series in the variables 
\begin{align*}
\left\{\frac{z_{1}}{z_{2}},\frac{z_{2}}{z_{3}},\ldots,\frac{z_{k-1}}{z_{k}},z_{k}\right\}. 
\end{align*}
Let $\oint dZ:\mathcal{L}(Z)\rightarrow \CC$, such that for each Laurent series $f\in \mathcal{L}(Z)$, $\oint f dZ$ is the coefficient of $\frac{1}{z_1\cdot\ldots\cdot z_k}$ in $f$.
\end{definition}

The following lemma about the commutative properties of the residue operator and the Macdonald scalar product was proved in \cite{Ah18}. 

\begin{lemma}\label{A12}(Lemma 3.8 in \cite{Ah18}) Let $\mathcal{A}$, $\mathcal{B}$ be graded algebras over $\CC$, let $f\in \mathcal{A}\widehat{\otimes}\Lambda_X[\mathcal{L}(Z)]$ and $g\in \Lambda_X\widehat{\otimes}\mathcal{B}[L(W)]$. If $f$ is $\Lambda_X$-projective, then
\begin{align*}
\left\langle \oint fdZ,g \right\rangle_X =\oint\langle f,g \rangle_X dZ; \\
\left\langle  f,\oint g dZ \right\rangle_X =\oint\langle f,g \rangle_X dZ.
\end{align*}
\end{lemma}

The following technical lemma is elementary, as proved in \cite{Ah18}.

\begin{lemma}\label{al51} (Lemma 5.7 of \cite{Ah18}) Let $\theta\in (0,\pi)$, and $\xi>0$. Define
\begin{align*}
R_{\epsilon,\theta,\xi}:=\{w\in \CC:\mathrm{dist}(w,[1,\infty))\leq \xi\}\cap \{w\in \CC:|\arg(w-(1-\epsilon))|\leq \theta\}.
\end{align*}
Let $\alpha>0$ and suppose $N(\epsilon)\in\ZZ>0$ such  that $\limsup_{\epsilon\rightarrow 0}\epsilon N(\epsilon)>0$ as $\epsilon\rightarrow 0$.  Then  for  any fixed $\theta\in (0,\pi),
\xi>0$, we have
\begin{align*}
\frac{(z;e^{-\epsilon})_{N(\epsilon)}}{(e^{-\epsilon \alpha} z;e^{-\epsilon})_{N(\epsilon)}}=\left(\frac{1-z}{1-e^{-\epsilon N(\epsilon)}z}\right)^{\alpha}\exp\left(O\left(\frac{\epsilon\min\{|z|,|z|^2\}}{|1-z|}\right)\right)
\end{align*}
uniformly for $z\in \CC\setminus R_{\epsilon,\theta,\xi}$ and $\epsilon$ arbitrarily small.
\end{lemma}

\begin{lemma}\label{lb2}(Corollary A.2 in \cite{GZ16}) Let $d,h,k$ be  positive  integers.  Let $f$, $g_1$,\ldots , $g_d$ be  meromorphic  functions  with  possible poles at ${z_1,...,z_h}$.  Then for $k\geq2$,
\begin{align*}
&\frac{1}{(2\pi\mathbf{i})^k}\oint\ldots\oint \frac{1}{(v_2-v_1)\cdots\ldots\cdot(v_k-v_{k-1})}\prod_{j=1}^{d}\left(\sum_{i=1}^kg_j(v_i)\right)\prod_{i=1}^k f(v_i)dv_i\\
&=\frac{k^{d-1}}{2\pi\mathbf{i}}\oint f(v)^k\prod_{j=1}^dg_j(v)dv.
\end{align*}
where  the  contours  contain $\{z_1,...,z_h\}$ and  on  the  left  side  we  require  that  the $v_i$-contour  is  contained  in the $v_j$-contour whenever $i < j$.
\end{lemma}

\bigskip

\bigskip
\noindent\textbf{Acknowledgements.} This material is based upon work supported by the National Science Foundation under Grant No. DMS-1928930 while ZL and MV participated in a program hosted by the Mathematical Sciences Research Institute in Berkeley, California, during the Fall 2021 semester. ZL acknowledges support from National Science Foundation DMS1608896 and Simons Foundation grant 638143. We thank Vadim Gorin for helpful discussions.

\bibliography{fpmm}
\bibliographystyle{plain}

\end{document}